\documentclass[amsfonts,11pt]{amsart}
\RequirePackage{amsmath,amssymb,amsthm, amscd, comment,mathtools}
\usepackage[margin=1.2in]{geometry}
\usepackage{amsmath}
\usepackage{amsthm}
\usepackage{amssymb}
\usepackage{amscd}
\usepackage{mathrsfs}
\usepackage{graphicx}
\usepackage{rotating, graphicx}
\usepackage[all,cmtip]{xy}
\usepackage{enumerate}
\usepackage{tabularx}
\usepackage{tikz-cd}
\usepackage[thinlines]{easytable}
\usepackage{multirow}
\usepackage{mathtools}
\usepackage{mathrsfs}
\usepackage{changepage} 
\usepackage{paralist}

\usepackage{kantlipsum}
\usepackage{mathabx}
\usepackage{enumitem}
\usepackage{times}
\usepackage{subfig}
\usepackage{color}
\usepackage[pagebackref,breaklinks,colorlinks,linkcolor=red,anchorcolor=red,citecolor=blue]{hyperref}

\calclayout

\newtheorem{proposition}[subsection]{Proposition}
\newtheorem{corollary}[subsection]{Corollary}
\newtheorem{theorem}[subsection]{Theorem}
\newtheorem{lemma}[subsection]{Lemma}
\newtheorem{conjecture}[subsection]{Conjecture}
\theoremstyle{definition}
\newtheorem{definition}[subsection]{Definition}
\theoremstyle{remark}

\newtheorem{remark}[subsection]{Remark}
\newtheorem{example}[subsection]{Example}

\numberwithin{equation}{subsection}

\DeclareMathAlphabet{\mathbbold}{U}{bbold}{m}{n}

\title[]{The limit and boundary characteristic classes in Borel-Moore motivic homology}

\author{Fangzhou Jin}
\address{School of Mathematical Sciences\\
Key Laboratory of Intelligent Computing and Applications (Ministry of Education)\\
Tongji University\\
Siping Road 1239\\
200092 Shanghai\\
China}
\email{\href{mailto:fangzhoujin@tongji.edu.cn}{fangzhoujin@tongji.edu.cn}}
\urladdr{\url{https://fangzhoujin.github.io/}}

\author{Peng Sun}
\address{School of Mathematics and Computational Science\\
Xiangtan University\\
411105 Xiangtan\\
China}
\email{\href{mailto:sunpeng@xtu.edu.cn}{sunpeng@xtu.edu.cn}}
\urladdr{\url{https://math.xtu.edu.cn/info/1010/4145.htm}}

\author{Enlin Yang}
\address{School of Mathematical Sciences\\
Capital Normal University\\
No.105, XiSanHuan North Road\\
Beijing 100048\\
P.R.China}
\email{\href{mailto: yangenlin@cnu.edu.cn}{yangenlin@cnu.edu.cn}}
\urladdr{\url{https://yelmath.github.io/}}

\date{\number\day-\number\month-\number\year}

\begin{document}

\maketitle

\begin{abstract}

We show that the zero-dimensional part of the pro-Chern-Schwarz-MacPherson class defined by Aluffi can be lifted to the zeroth Suslin homology. The proof uses the pro-characteristic class in the limit Borel-Moore motivic homology, which has a quadratic refinement in the limit Borel-Moore Milnor-Witt homology. In characteristic zero, this construction factors through the group of constructible functions, in a way compatible with the covariant functoriality; in positive characteristic this property fails, and we show that the failure can be measured by the boundary characteristic class in the boundary Borel-Moore motivic homology. We prove a push-forward formula for the boundary characteristic class, and conjecture it to agree with the Swan class defined by Kato-Saito.

\end{abstract}

\setcounter{tocdepth}{1}
\tableofcontents

\noindent
\section{Introduction}
\subsection{}
Inspired by the Chern-Weil theory, Grothendieck shows how to define Chern classes for vector bundles in algebraic geometry with values in the group of algebraic cycles up to rational equivalence (\cite{Gro}), namely the \emph{Chow group}, considered as an algebraic substitute of the cohomology ring of a complex manifold. As a particular case of this general construction, one can define the Chern class of a smooth algebraic variety $X$ as the (total) Chern class of its tangent bundle:
\begin{align}
\label{eq:cX}
c(X)=c(T_X)\in CH_*(X).
\end{align}

\subsection{}
Over the field of complex numbers, the landmarking work of MacPherson (\cite{Mac}) extends the Chern class~\eqref{eq:cX} to every (possibly singular) complex algebraic variety $X$, in a way compatible with proper push-forward maps of constructible functions, and therefore answering affirmatively to a conjecture of Deligne-Grothendieck. MacPherson's original construction takes place in homology and uses transcendental methods, and later algebraic formulas have been found by Gonz\'alez-Sprinberg and Verdier (\cite{Gon}), leading to what is now called the \emph{Chern-Schwarz-MacPherson class} in the Chow group:
\begin{align}
c^{SM}(X)\in CH_*(X).
\end{align}

\subsection{}
In \cite{Alu}, Aluffi gives a refinement of MacPherson's class via an algebraic-geometric construction using resolution of singularities, by defining the \emph{pro-Chern-Schwarz-MacPherson class} (abbreviated as \emph{pro-CSM class})
\begin{align}
lC^{SM}(X)\in lCH_*(X)
\end{align}
in the \emph{pro-Chow group} of $X$, that is, the limit of Chow groups of all possible compactifications of $X$ (\cite[Def. 2.2]{Alu}). This construction is briefly recalled in Theorem~\ref{thm:Alu43} below, and the functoriality of the pro-CSM class is extended to all morphisms of varieties (\cite[Thm. 5.2]{Alu}).

\subsection{}
In the first part of this paper, we provide a category-theoretic approach to the $0$-dimensional component $lC^{SM}_0(X)$ of Aluffi's pro-CSM class $lC^{SM}(X)$ via motivic homotopy theory. 
In a previous work (\cite{JY}), the first and third-named authors defined a \emph{characteristic class} in the Borel-Moore homology as a generalized trace map (\cite[Def. 5.1.3]{JY}) 
\begin{align}
\label{eq:CXintro}
C_X(K)\in \mathbb{E}^{\mathrm{BM}}(X/k)
\end{align}
where $K$ is a constructible motivic spectrum over $X$, and $\mathbb{E}^{\mathrm{BM}}(X/k)$ is the Borel-Moore theory associated to a motivic spectrum $\mathbb{E}$ (see~\ref{num:BM}), 
assuming resolution of singularities or inverting the characteristic (see~\ref{num:RS}). 
The class $C_X(K)$ is defined by analogy with a construction due to Verdier in various categories of sheaves (\cite[Exp. III]{SGA5}, \cite[Def. 9.1.2]{KS}, \cite{AS}).

For example, taking $\mathbb{E}=\mathbf{H}\mathbb{Z}$ to be the motivic Eilenberg-Mac Lane spectrum (\cite{Spi}, see Example~\ref{ex:EM}), we obtain a class $C_X(K)\in CH_0(X)$ in the Chow group of $0$-cycles. In this case, the definition of this class is due to Olsson when $k$ is an algebraically closed field (\cite{Ols}).

\subsection{}
Similar to Aluffi's definition, we define a pro-version of the Borel-Moore theory $l\mathbb{E}^{\mathrm{BM}}(X/k)$ (Definition~\ref{def:limBM}). However, as a substitute of Aluffi's geometric approach, we use the six-functors formalism in motivic homotopy (\cite{Ayo}) to define two pro-versions of the class $C_X(K)$, denoted as $lC_X(K)$ and $lC^{\mathrm{c}}_X(K)$ respectively (Definition~\ref{df:proclass}), which are respectively the homological and compactly supported variants of the pro-characteristic class. Concretely, the class $lC_X(K)$ (respectively $lC^{\mathrm{c}}_X(K)$) is the pro-class associated to the class $C_{\overline{X}}(j_*K)$ (respectively $C_{\overline{X}}(j_!K)$) for every compactification $X\xrightarrow{j}\overline{X}$ of $X$. These two classes are related by the local duality functor (Corollary~\ref{cor:Ccdual}), and satisfy push-forward (\ref{num:Cpf}) and additivity formulas (\ref{num:Cadd}). For the relation with Aluffi's push-forward formula (\cite[Thm. 5.2]{Alu}), see Remark~\ref{rem:Apf}.

\subsection{}
Compared to Aluffi's construction which assumes the existence of weak factorizations (which is known in characteristic $0$ by \cite{AKMW}), we only need a weaker version of resolution of singularities (by blow-ups or by alterations, see~\ref{num:RS}). Indeed, this assumption is used to prove the K\"unneth formulas in motivic homotopy (\cite[Thm. 2.4.6]{JY}), which is enough to define the characteristic class~\eqref{eq:CXintro}. In particular, our construction is unconditional with $\mathbb{Z}[1/p]$-coefficients in positive characteristic $p$. For a discussion on generalizations to higher dimensional cases, see Remark~\ref{rem:CC}.

\subsection{}
One of our main results states that the $0$-dimensional part of Aluffi's pro-CSM class $lC^{SM}_0(X)$ is equal to the compactly supported pro-characteristic class of the unit object $\mathbbold{1}_X$:
\begin{theorem}[\textrm{see Theorem~\ref{th:A=C}}]
For every $k$-scheme $X$, we have 
\begin{align}
lC^{SM}_0(X)=lC^{\mathrm{c}}_X(\mathbbold{1}_X)\in lCH_0(X).
\end{align}
\end{theorem}
The proof uses the properties of the class $lC^{\mathrm{c}}_X(K)$, as well as the computation of the Chern classes of the sheaf of differentials with logarithmic poles as in the proof of \cite[Thm. 3.1]{Sil}. One can also take the spectrum $\mathbb{E}$ in~\eqref{eq:CXintro} to be the Milnor-Witt spectrum (\cite{DF}, see Example~\ref{ex:MW}), in which case the pro-characteristic class $lC^{\mathrm{c}}_X(\mathbbold{1}_X)$ gives a quadratic refinement of the class $lC^{SM}_0(X)$, see~\ref{num:MW}.

Note that in a recent work (\cite{Azo}), Azouri gives independently a construction of the compactly supported pro-characteristic class $lC^{\mathrm{c}}_X(\mathbbold{1}_X)$, with an approach closer to Aluffi's original work.

\subsection{}
Recall that the \emph{Suslin homology} is a homology theory on algebraic varieties analogous to singular homology (see \cite[Prop. 14.18]{MVW} for a precise statement).
The zeroth Suslin homology group can be written as
\begin{align}
H_0^{\operatorname{S}}(X/k)=\operatorname{coker}(\operatorname{Cor}(\mathbb{A}^1,X)\to Z_0(X)).
\end{align}
It is covariant for all $k$-morphisms, and there are maps $H_0^{\operatorname{S}}(X/k)\to lCH_0(X)$ compatible with covariant functoriality on both sides.
We enhance the $lC^{SM}_0(X)$ to a \emph{homological characteristic class} in the $0$-th Suslin homology group:
\begin{theorem}[\textrm{see Theorem~\ref{thm:Constr}}]
Let $k$ be a field of characteristic $0$. There is a unique natural transformation of functors on $k$-varieties
\begin{align}
C^{\operatorname{Hom}}:\operatorname{Cons}(-)\to H_0^{\operatorname{S}}(-/k)
\end{align}
which commutes with push-forwards on both sides, and which lifts the class $lC^{SM}_0$.
\end{theorem}
Note that one cannot expect to lift higher-dimensional parts of the pro-CSM class, see Remark~\ref{rm:higherdim} below.

\subsection{}
In the second part of the paper, we focus on the positive characteristic case. In positive characteristic, the characteristic calss no longer factors through $\operatorname{Cons}(-)$, as shown by the famous Grothendieck–Ogg–Shafarevich formula (\cite[X Thm. 7.1]{SGA5}), where the defect is measured by the \emph{Swan conductor}. In the motivic setting, based on previous work of the first and third-named authors (\cite{JY2}) and inspired by the work of Kato-Saito (\cite{KSa}), we introduce a \emph{boundary characteristic class} $bC^{\mathrm{c}}_X$ (Definition~\ref{def:bCc}) in the \emph{boundary Borel-Moore theory} $bCH_0(X)$ (Definition~\ref{def:bBM}), assuming the existence of smooth compatifications (see~\ref{num:smcomp}). This class measures the difference between $lC^{\mathrm{c}}_X$ and the class provided by the rank function (see~\ref{num:blC}), and satisfies the following push-forward formula:
\begin{theorem}[\textrm{see Theorem~\ref{th:pfbC}}]
Assume that $k$ is a perfect field of characteristic $p>0$ such that all $k$-varieties have smooth compactifications.
Let $f:Y\to X$ be a finite \'etale morphism between smooth connected $k$-schemes. Let $K\in\mathbf{SH}_c(Y)$ be a dualizable object, then we have
\begin{align}
\label{eq:bccpfintro}
bC^{\mathrm{c}}_X(f_*K)=f_*bC^{\mathrm{c}}_Y(K)+\operatorname{rk}(K)\cdot bC^{\mathrm{c}}_X(f_*\mathbbold{1})\in bCH_0(X)[1/p].
\end{align}
\end{theorem}
The proof relies on Proposition~\ref{prop:CXZpf}, whose proof is technical and occupies the whole Section~\ref{sec:pfpbpf}.

\subsection{}
Following the philosophy of uniqueness of the characteristic class (see \cite[Conjecture 6.8]{Sai} and \cite{YZ}), we further conjecture that up to a sign, the class $bC^{\mathrm{c}}_X$ agrees with the Kato-Saito Swan class (\cite{KSa}), see Conjecture~\ref{conj:b=Sw}. Since the Kato-Saito Swan class also satisfies a similar push-forward formula (\cite[Cor. 4.3.4]{KSa}), our formula~\eqref{eq:bccpfintro} provides some evidence for this conjecture.

\subsubsection*{\bf Acknowledgments}
We would like to thank Marc Levine for Remark~\ref{rem:levquad}, Ran Azouri, Lie Fu and Xiping Zhang for very helpful discussions.
The authors are supported by the National Key Research and Development Program of China Grant Nr. 2021YFA1001400. F. Jin is supported by the National Natural Science Foundation of China Grant Nr.12471014, Nr.12101455 and the Fundamental Research Funds for the Central Universities. E. Yang is supported by the NSFC Grant Nr.12271006.

\section{Homology theories in motivic homotopy}

\subsection{}
\label{num:RS}
Throughout the paper, all schemes are assumed separated and of finite type over a field $k$ of exponential characteristic $p$. We consider the following conditions on \textbf{embedded resolution of singularities} over $k$:
    \begin{enumerate}
    \item For every separated integral scheme $X$ of finite type over $k$, there exists a proper birational surjective morphism $X'\to X$ with $X'$ regular;
    \item For every separated integral regular scheme $X$ of finite type over $k$ and every nowhere dense closed subscheme $Z$ of $X$, there exists a proper birational surjective morphism $b:X'\to X$ such that $X'$ is regular, $b$ induces an isomorphism $b^{-1}(X-Z)\simeq X-Z$, and $b^{-1}(Z)$ is a strict normal crossing divisor in $X'$.
    \end{enumerate}

\subsection{}
\label{num:dual}
Let $\mathcal{C}$ be a symmetric monoidal category with unit object $\mathbbold{1}$. If $M$ is a  (strongly) dualizable object in $\mathcal{C}$, denote by $M^\vee$ its (strong) dual, and define the \textbf{categorical Euler characteristic} $\chi(M)$ of $M$ as the composition
\begin{align}
\chi(M):\mathbbold{1}\to M^\vee\otimes M\simeq M\otimes M^\vee\to\mathbbold{1}
\end{align}
considered as an endomorphism of the unit object $\mathbbold{1}$ (\cite[III Def. 7.1]{LMS}). By \cite[III Prop. 7.7]{LMS}, we have a canonical identification $\chi(M)=\chi(M^\vee)$. In other words, the categorical Euler characteristic of a dualizable object agrees with that of its dual.

\subsection{}
For every scheme $X$, we denote by $\mathbf{SH}(X)$ the \textbf{stable motivic homotopy category}, which is endowed with a \emph{six functors formalism} (\cite{Ayo}). The unit object of $\mathbf{SH}(X)$ is denoted as $\mathbbold{1}_X$. If $i:Z\to X$ is a closed immersion with $j:U\to X$ the open complement, there are distinguished triangles of endofunctors of $\mathbf{SH}(X)$ called  \textbf{localization triangles}:
\begin{align}
\label{eq:locseq}
j_!j^!\to1\to i_*i^*,
\end{align}
\begin{align}
\label{eq:locseq1}
i_*i^!\to1\to j_*j^!.
\end{align}
We denote by $\mathbf{SH}_c(X)$ the full subcategory of \textbf{constructible objects} in $\mathbf{SH}(X)$ (\cite[Def. 2.2.3]{Ayo}).

\subsection{}
\label{num:FHM32}
(\cite[Prop. 3.2]{FHM})
If $f:Y\to X$ is a morphism of schemes, $K\in\mathbf{SH}(Y)$ and $L\in\mathbf{SH}(X)$ are two objects, by adjunction we obtain a canonical natural transformation of the form
\begin{align}
\label{eq:FHM32}
f_*K\otimes L
\to
f_*(K\otimes f^*L)
\end{align}
which is an isomorphism if $L$ is dualizable.

\subsection{}
\label{num:BM}
(\cite[2.2.1]{DM})
Let $f:X\to S$ be a morphism of schemes and let $\mathbb{E}\in \mathbf{SH}(S)$ be a motivic spectrum. We define the \textbf{$\mathbb{E}$-homology}
\begin{align}
\mathbb{E}^{}_n(X/S)
=
[\mathbbold{1}_S[n],f_!f^!\mathbb{E}]_{\mathbf{SH}(S)}
\end{align}
and the \textbf{Borel-Moore $\mathbb{E}$-homology}
\begin{align}
\mathbb{E}^{\mathrm{BM}}_n(X/S)
=
[\mathbbold{1}_S[n],f_*f^!\mathbb{E}]_{\mathbf{SH}(S)}.
\end{align}

\subsection{}
The natural transformation of functors $f_!\to f_*$ induces a map
\begin{align}
\label{eq:homBM}
\mathbb{E}^{}_n(X/S)\to\mathbb{E}^{\mathrm{BM}}_n(X/S)
\end{align}
which is an isomorphism when $f$ is proper.

\subsection{}
For any morphism $p:Y\to X$, there is a map
\begin{align}
\label{eq:hompf}
p_*:
\mathbb{E}^{}_n(Y/S)
\to
\mathbb{E}^{}_n(X/S).
\end{align}
In particular, for $X=S$ and $n=0$ we have the \textbf{degree map}
\begin{align}
\label{eq:deghom}
\int_{Y/S}:
\mathbb{E}^{}_0(Y/S)
\to
\mathbb{E}^{}_0(S/S).
\end{align}

\subsection{}
For any proper morphism $p:Y\to X$, there is a map
\begin{align}
\label{eq:BMprop}
p_*:
\mathbb{E}^{\mathrm{BM}}_n(Y/S)
\to
\mathbb{E}^{\mathrm{BM}}_n(X/S).
\end{align}
For any \'etale morphism $p:Y\to X$, there is a map 
\begin{align}
\label{eq:BMet}
p^*:
\mathbb{E}^{\mathrm{BM}}_n(X/S)
\to
\mathbb{E}^{\mathrm{BM}}_n(Y/S)
\end{align}
such that the following diagram commutes:
\begin{align}
\label{eq:HomBMet}
\begin{split}
  \xymatrix@R=10pt{
    \mathbb{E}^{}_n(X/S) \ar[r]^-{\eqref{eq:homBM}} & \mathbb{E}^{\mathrm{BM}}_n(X/S) \ar[d]^-{\eqref{eq:BMet}} \\
    \mathbb{E}^{}_n(Y/S) \ar[r]^-{\eqref{eq:homBM}} \ar[u]^-{\eqref{eq:hompf}} & \mathbb{E}^{\mathrm{BM}}_n(Y/S).
  }
\end{split}
\end{align}

\subsection{}
If $Z\to X$ is a closed immersion with open complement $U$, by the localization triangle~\eqref{eq:locseq1} there is a long exact sequence
\begin{align}
\cdots\to
\mathbb{E}^{\mathrm{BM}}_{n+1}(U/S)
\to
\mathbb{E}^{\mathrm{BM}}_{n}(Z/S)
\to
\mathbb{E}^{\mathrm{BM}}_n(X/S)
\to
\mathbb{E}^{\mathrm{BM}}_n(U/S)
\to\cdots.
\end{align}
By the localization triangle~\eqref{eq:locseq} there is a long exact sequence
\begin{align}
\label{eq:homles}
\cdots\to
\mathbb{E}^{}_{n+1,Z}(X/S)
\to
\mathbb{E}^{}_n(U/S)
\to
\mathbb{E}^{}_n(X/S)
\to
\mathbb{E}^{}_{n,Z}(X/S)
\to\cdots.
\end{align}
Here $\mathbb{E}^{}_{n,Z}(X/S)=[\mathbbold{1}_{S}[n],f_!i_*i^*f^!\mathbb{E}]_{\mathbf{SH}(S)}$ is the \emph{$\mathbb{E}$-homology with support in $Z$}.

\subsection{}
\label{num:comp}
Let $f:X\to S$ be a morphism of schemes. We denote by $\operatorname{Cpt}(X/S)$ the \textbf{category of compactifications} of $f$, such that
\begin{enumerate}
\item
The objects are \textbf{compactifications} of $X$ over $S$, that is, factorizations of $f$ as 
\begin{align}
\label{eq:compfac}
X\xrightarrow{j}\overline{X}\xrightarrow{p}S
\end{align}
with $j$ an open immersion with dense image 
and $p$ a proper morphism; a compactification of the form~\eqref{eq:compfac} is sometimes denoted as $X\xrightarrow{j}\overline{X}$ or simply $\overline{X}$; 
\item
Morphisms from $X\xrightarrow{j}\overline{X}\xrightarrow{p}S$ to $X\xrightarrow{j'}\overline{X}'\xrightarrow{p'}S$ are proper morphisms $\overline{X}\xrightarrow{g}\overline{X}'$ such that $j'=f\circ j$, $p=p'\circ g$, that is, such that there is a commutative diagram
\begin{align}
\begin{split}
  \xymatrix@=10pt{
    X \ar[r]^-{j} \ar[d]_-{j'} & \overline{X} \ar[d]^-{p} \ar[ld]_-{g} \\
    \overline{X}' \ar[r]_-{p'} & S.
  }
\end{split}
\end{align}
\end{enumerate}
The category $\operatorname{Cpt}(X/S)$ is non-empty and cofiltered, see \cite{Con}, \cite[Exp. XVII \S3.2]{SGA4}, \cite[Tag 0ATT]{Stack}.

\subsection{}
Let $X$ be a separated $k$-scheme of finite type.
We denote by $\operatorname{Cpt}^{\operatorname{Sm}}(X/k)$ the \textbf{category of smooth compactifications} of $X$, which is the full subcategory of $\operatorname{Cpt}(X/k)$ consisting of compactifications $\overline{X}$ which are smooth over $k$. 

\subsection{}
\label{num:smcomp}
If $k$ is a perfect field, the statement that the full subcategory $\operatorname{Cpt}^{\operatorname{Sm}}(X/k)$ is cofinal in $\operatorname{Cpt}(X/k)$ is related to resolution of singularities, and such a property is known to hold in the following cases:
\begin{enumerate}
\item
If $k$ is a perfect field and $X$ has dimension at most $3$ (\cite{CP}).
\item
If $k$ is a perfect field which satisfies embedded resolution of singularities in the sense of~\ref{num:RS}. Furthermore, in this case one may replace $\operatorname{Cpt}^{\operatorname{Sm}}(X/k)$ by the full subcategory spanned by objects $\overline{X}\in\operatorname{Cpt}^{\operatorname{Sm}}(X/k)$ such that $\overline{X}-X$ is a simple normal-crossing divisor in $\overline{X}$.
\end{enumerate}

\subsection{}
Using the ideas of \cite{Alu}, one may promote the proper functoriality of Borel-Moore homology~\eqref{eq:BMprop} to all morphisms by looking at the followng theory:
\begin{definition}
\label{def:limBM}
We define the \textbf{limit Borel-Moore $\mathbb{E}$-homology}
\begin{align}
l\mathbb{E}^{\mathrm{BM}}_n(X/S)
=
\operatorname{lim}_{\overline{X}\in\operatorname{Cpt}(X/S)}\mathbb{E}^{\mathrm{BM}}_n(\overline{X}/S)
\end{align}
where the transition maps are given by proper functoriality of Borel-Moore homology~\eqref{eq:BMprop} induced by morphisms of compactifications.
\end{definition}

\subsection{}
\label{num:limBM}
By definition, to determine a class $\alpha$ in $l\mathbb{E}^{\mathrm{BM}}(X/S)$ amounts to determine, for every compactification $\overline{X}$ of $f$, a class $\alpha_{\overline{X}}$ in $\mathbb{E}^{\mathrm{BM}}(\overline{X}/S)$, compatible with proper push-forwards.

\subsection{}
Let $f:X\to S$ be a morphism of schemes. There are canonical maps
\begin{align}
f_!f^!\mathbb{E}
\to
\underset{(X\xrightarrow{j}\overline{X}\xrightarrow{p}S)\in\operatorname{Cpt}(X/S)}{\operatorname{lim}}p_*p^!\mathbb{E}
\to
f_*f^!\mathbb{E}
\end{align}
which, by applying the functor $[\mathbbold{1}_S[n],-]_{\mathbf{SH}(S)}$, give rise to maps
\begin{align}
\label{eq:HomlimBM}
\mathbb{E}^{}_n(X/S)
\to
l\mathbb{E}^{\mathrm{BM}}_n(X/S)
\end{align}
and
\begin{align}
\label{eq:limBMBM}
l\mathbb{E}^{\mathrm{BM}}_n(X/S)
\to
\mathbb{E}^{\mathrm{BM}}_n(X/S).
\end{align}
Concretely, the map~\eqref{eq:limBMBM} is also given by the \'etale functoriality $j^*:\mathbb{E}^{\mathrm{BM}}_n(\overline{X}/S)\to \mathbb{E}^{\mathrm{BM}}_n(X/S)$ in~\eqref{eq:BMet} for every compactification $X\xrightarrow{j}\overline{X}$.

\subsection{}
The composition of the two maps~\eqref{eq:HomlimBM} and~\eqref{eq:limBMBM} agrees with the map~\eqref{eq:homBM}. Both maps are isomorphisms when $f$ is proper.

\subsection{}
Let $f:Y\to X$ be a morphism of schemes. We now define a map
\begin{align}
\label{eq:funclimBM}
f_*:
l\mathbb{E}^{\mathrm{BM}}_n(Y/S)
\to
l\mathbb{E}^{\mathrm{BM}}_n(X/S).
\end{align}
Indeed, by \cite[XVII, Prop. 3.2.6]{SGA4}, for every compactification $\overline{X}$ of $X$, there exists a compactification $\overline{Y}$ of $Y$, a proper morphism $\overline{f}:\overline{Y}\to\overline{X}$ and a commutative diagram
\begin{align}
\label{eq:comprel}
\begin{split}
  \xymatrix@=10pt{
    Y \ar[r]^-{} \ar[d]_-{f} & \overline{Y} \ar[d]^-{\overline{f}} \\
    X \ar[r]^-{} & \overline{X}.
  }
\end{split}
\end{align}
By~\ref{num:limBM}, we define the map~\eqref{eq:funclimBM} such that for every class $\alpha\in l\mathbb{E}^{\mathrm{BM}}_n(Y/S)$, the direct image $(f_*\alpha)\in l\mathbb{E}^{\mathrm{BM}}_n(X/S)$ is the class associated, for every compactification $\overline{X}$ of $X$, to the class $(f_*\alpha)_{\overline{X}}=\overline{f}_{*}\alpha_{\overline{Y}}$, which is well-defined and does not depend on the choice of $\overline{Y}$. 

In particular, for $X=S$ and $n=0$ we have the \textbf{degree map}
\begin{align}
\label{eq:funclimBM}
\int_{Y/S}:
l\mathbb{E}^{\mathrm{BM}}_0(Y/S)
\to
l\mathbb{E}^{\mathrm{BM}}_0(S/S).
\end{align}

\subsection{}
If $p:Y\to X$ is a morphism, we have the following commutative diagram:
\begin{align}
\label{eq:Homlimdiag}
\begin{split}
  \xymatrix@R=10pt{
    \mathbb{E}^{}_n(Y/S) \ar[r]^-{\eqref{eq:HomlimBM}} \ar[d]_-{\eqref{eq:hompf}} & l\mathbb{E}^{\mathrm{BM}}_n(Y/S) \ar[d]^-{\eqref{eq:funclimBM}} \\
    \mathbb{E}^{}_n(X/S) \ar[r]^-{\eqref{eq:HomlimBM}} & l\mathbb{E}^{\mathrm{BM}}_n(X/S).
  }
\end{split}
\end{align}
In particular, for $X=S$ and $n=0$ the following diagram is commutative:
\begin{align}
\begin{split}
  \xymatrix@R=10pt{
    \mathbb{E}^{}_0(Y/S) \ar[r]^-{\eqref{eq:HomlimBM}} \ar[d]_-{\int_{Y/S}} & l\mathbb{E}^{\mathrm{BM}}_0(Y/S) \ar[d]^-{\int_{Y/S}} \\
    \mathbb{E}^{}_0(S/S) \ar[r]^-{\sim} & l\mathbb{E}^{\mathrm{BM}}_0(S/S).
  }
\end{split}
\end{align}

\subsection{}
If $p:Y\to X$ is a proper morphism, we have the following commutative diagram:
\begin{align}
\label{eq:limBMdiag}
\begin{split}
  \xymatrix{
    l\mathbb{E}^{\mathrm{BM}}_n(Y/S) \ar[r]^-{\eqref{eq:limBMBM}} \ar[d]_-{\eqref{eq:funclimBM}} & \mathbb{E}^{\mathrm{BM}}_n(Y/S) \ar[d]^-{\eqref{eq:BMprop}} \\
    l\mathbb{E}^{\mathrm{BM}}_n(X/S) \ar[r]^-{\eqref{eq:limBMBM}} & \mathbb{E}^{\mathrm{BM}}_n(X/S).
  }
\end{split}
\end{align}
The diagrams~\eqref{eq:Homlimdiag} and~\eqref{eq:limBMdiag} are compatible with compositions of morphisms, see \cite[Lemma 2.4]{Alu}.

\begin{example}
\label{ex:EM}
Let $S=\operatorname{Spec}(k)$ be the spectrum of a field of characteristic $p$, and $\mathbb{E}=\mathbf{H}\mathbb{Z}$ the motivic Eilenberg-Mac Lane spectrum (\cite{Spi}).
\begin{enumerate}
\item
The motivic homology $\mathbf{H}\mathbb{Z}^{}_n(X/k)$ is isomorphic to the \textbf{Suslin homology} (\cite[Prop. 14.18]{MVW}):
\begin{align}
\mathbf{H}\mathbb{Z}^{}_n(X/k)
\simeq
H_n^{\operatorname{S}}(X/k).
\end{align}
In particular for $n=0$ we obtain the zeroth Suslin homology group
\begin{align}
\mathbf{H}\mathbb{Z}^{}_0(X/k)
\simeq
H_0^{\operatorname{S}}(X/k)=\operatorname{coker}(\operatorname{Cor}(\mathbb{A}^1,X)\to Z_0(X)).
\end{align}
\item
The Borel-Moore motivic homology $\mathbf{H}\mathbb{Z}^{\mathrm{BM}}_0(X/k)$ agrees with the Chow group of zero cycles on $X$:
\begin{align}
\label{eq:BMCH0}
\mathbf{H}\mathbb{Z}^{\mathrm{BM}}_0(X/k)
\simeq
CH_0(X).
\end{align}
and the map $\mathbf{H}\mathbb{Z}^{}_0(X/k)\to\mathbf{H}\mathbb{Z}^{\mathrm{BM}}_0(X/k)$ in~\eqref{eq:homBM} is identified with the canonical surjection $H_0^{\operatorname{S}}(X/k)\to CH_0(X)$ (\cite[Exercise 2.21]{MVW}).

\item
By~\eqref{eq:BMCH0}, it follows that the limit Borel-Moore motivic homology $l\mathbf{H}\mathbb{Z}^{\mathrm{BM}}_0(X/k)$
agrees with the \textbf{pro-Chow group} of zero-cycles defined by Aluffi (\cite[Def. 2.2]{Alu}):
\begin{align}
\label{eq:limBMHZ}
l\mathbf{H}\mathbb{Z}^{\mathrm{BM}}_0(X/k)
\simeq
lCH_0(X)
=
\operatorname{lim}_{\overline{X}\in\operatorname{Cpt}(X/k)}
CH_{0}(\overline{X}).
\end{align}
and the map $\mathbf{H}\mathbb{Z}^{}_0(X/k)\to l\mathbf{H}\mathbb{Z}^{\mathrm{BM}}_0(X/k)$ in~\eqref{eq:HomlimBM} is identified with the canonical map $H_0^{\operatorname{S}}(X/k)\to lCH_0(X)$.
\end{enumerate}
Note the isomorphisms~\eqref{eq:BMCH0} and~\eqref{eq:limBMHZ} can be extended to (pro-)Chow groups of cycles of any dimension (see for example \cite{Jin}), which we do not discuss here.
\end{example}

\begin{example}
\label{ex:MW}
For $\mathbb{E}=\mathbf{H}\mathbb{Z}_{\mathrm{MW}}$ the Milnor-Witt spectrum and $S$ the spectrum of an infinite perfect field $k$ (\cite{DF}), the Borel-Moore Milnor-Witt homology $\mathbf{H}\mathbb{Z}^{\mathrm{BM}}_{\mathrm{MW}}(X/k)$ agrees with the \emph{Chow-Witt group} of zero cycles on $X$, which is a quadratic refinement of the usual Chow group (see \cite[\S8]{DFJK}):
\begin{align}
\mathbf{H}\mathbb{Z}^{\mathrm{BM}}_{\mathrm{MW},0}(X/k)
\simeq
CH^{\mathrm{MW}}_{0}(X).
\end{align}
The limit Borel-Moore Milnor-Witt homology agrees with the \textbf{pro-Chow-Witt group} of zero-cycles defined in a similar way
\begin{align}
l\mathbf{H}\mathbb{Z}_{\mathrm{MW},0}^{\mathrm{BM}}(X/k)
\simeq
lCH^{\mathrm{MW}}_{0}(X)
=
\operatorname{lim}_{\overline{X}\in\operatorname{Cpt}(X/k)}
CH^{\mathrm{MW}}_{0}(\overline{X}).
\end{align}
\end{example}

\begin{remark}
By \cite{Fel}, the Chow-Witt group of zero cycles is a birational invariant:
if $g:Y\to X$ is a proper birational morphism between smooth proper $k$-schemes, then the map $g_*:CH^{\mathrm{MW}}_0(Y)\to CH^{\mathrm{MW}}_0(X)$ is an isomorphism. It follows that if $k$ is perfect and $\operatorname{Cpt}^{\operatorname{Sm}}(X/k)$ is cofinal in $\operatorname{Cpt}(X/k)$ as in~\ref{num:smcomp}, then we have
\begin{align}
lCH^{\mathrm{MW}}_{0}(X)
=
\operatorname{lim}_{\overline{X}\in\operatorname{Cpt}^{\operatorname{Sm}}(X/k)}
CH^{\mathrm{MW}}_{0}(\overline{X})
\end{align}
where the transition morphisms in the limit on the right-hand side are isomorphisms. Therefore $lCH^{\mathrm{MW}}_{0}(X)\simeq CH^{\mathrm{MW}}_{0}(\overline{X})$ for \emph{any} smooth compactification $\overline{X}$ of $X$. Note that the same property holds if one replaces Chow-Witt groups by Chow groups.
\end{remark}

\subsection{}
We end this section by proving an analogue of \cite[Th. 1.1]{Lau} for motives.
For any scheme $X$ of characteristic $p$, denote by $H(X)$ the quotient of the Grothendieck ring $K_0(\mathbf{DM}_c(X)[1/p])$ by the ideal generated by $[\mathbbold{1}(1)]-[\mathbbold{1}]$. We now analyse the relation between local duality and classes in $H(X)$.
\begin{lemma}
\label{lm:dualcls}
Let $X$ be a scheme. Then for any constructible object $K\in\mathbf{DM}_c(X)[1/p]$, the two objects $K$ and $\mathbb{D}_{X/k}(K)$ have the same class in $H(X)$.
\end{lemma}
\proof
We know that there is a bounded Chow weight structure on $\mathbf{DM}_c(X)[1/p]$, and therefore it suffices to prove the case where $K$ is a direct summand of $g_*\mathbbold{1}_Y(n)$, where $g:Y\to X$ is a proper morphism with $Y$ smooth over a finite purely inseparable extension of $k$. In this case the result follows from a direct computation: indeed, if $K$ is represented by an idempotent endomorphism $u$ of $g_*\mathbbold{1}_Y(n)$, where $d$ is the dimension of $Y$, then $\mathbb{D}_{X/k}(K)$ is isomorphic to the object represented by the idempotent endomorphism $u(d)[2d]$ of $g_*\mathbbold{1}_Y(n+d)[2d]$, which is isomorphic to $K(d)[2d]$ and therefore has the same class as $K$ in $H(X)$.
\endproof

\begin{proposition}
\label{prop:Lau1.1}
Let $f:Y\to X$ be a morphism of schemes. Then for any constructible object $K\in\mathbf{DM}_c(Y)[1/p]$, the two objects $f_*K$ and $f_!K$ have the same class in the group $H(X)$.
\end{proposition}
\proof
Since the six functors commute with Tate twists, if two objects $K$ and $K'$ have the same class in $H(Y)$, then the two objects $f_!K$ and $f_!K'$ have the same class in $H(X)$.
By Lemma~\ref{lm:dualcls}, in the group $H(X)$ we have an equality of classes
\begin{align}
[f_*K]=[\mathbb{D}_{X/k}(f_*K)]=[f_!\mathbb{D}_{X/k}(K)]=[f_!K]
\end{align}
which finishes the proof.
\endproof

\section{Characteristic classes in motivic homology theories}

\subsection{}
In this section, we work over a field $k$ of exponential characteristic $p$.

\subsection{}
\label{num:CX}
Let $f:X\to k$ be the structure morphism with $\delta:X\to X\times_kX$ the diagonal morphism. Let $K\in \mathbf{SH}_c(X)$ be a constructible motivic spectrum. Let $\mathbb{E}\in \mathbf{SH}(k)$ be a motivic spectrum with a map $\mathbbold{1}_k\to \mathbb{E}$. In \cite[Def. 5.1.3]{JY}, after inverting $p$, we constructed a map
\begin{align}
\begin{split}
\label{eq:CX}
C_X(K,\mathbb{E}):\mathbbold{1}_X
&\to
\underline{Hom}(K,K)
\simeq
\delta^!(\mathbb{D}_{X/k}(K)\boxtimes_kK)
\to
\delta^*(\mathbb{D}_{X/k}(K)\boxtimes_kK)\\
&=
\mathbb{D}_{X/k}(K)\otimes_kK
\simeq
K\otimes_k\mathbb{D}_{X/k}(K)
\to
f^!\mathbbold{1}_k
\to
f^!\mathbb{E}
\end{split}
\end{align}
viewed as an element of $\mathbb{E}^{\mathrm{BM}}_0(X/k)[1/p]$. This class $C_X(K)\in\mathbb{E}^{\mathrm{BM}}_0(X/k)[1/p]$ is called the ($\mathbb{E}$-valued) \textbf{characteristic class} of $K$. 

\subsection{}
If $X=k$, then $C_X(K)$ agrees with the categorical Euler characteristic of $K$ (see \ref{num:dual}).

\subsection{}
If $k$ is a perfect field which satisfies embedded resolution of singularities in the sense of~\ref{num:RS}, then the class $C_X(K)$ can be defined in $\mathbb{E}^{\mathrm{BM}}_0(X/k)$ without inverting $p$.

\subsection{}
(\cite[Thm. 4.2.8]{JY})
\label{num:adddt}
For a distinguished triangle $L\to M\to N$ of constructible objects in $\mathbf{SH}_c(X)$, we have
\begin{align}
\label{eq:adddt}
C_X(M)=C_X(L)+C_X(N).
\end{align}

\subsection{}
\label{num:CXprop}
(\cite[Cor. 5.1.8]{JY})
If $p:Y\to X$ is a proper morphism, then 
\begin{align}
\label{eq:CXprop}
p_*C_Y(K)=C_X(p_*K).
\end{align}

\subsection{}
(\cite[Example 5.1.16]{JY})
If $X$ is smooth of relative dimension $n$ over $k$, then 
\begin{align}
\label{eq:JY5116}
C_X(\mathbbold{1}_X)=e(T_{X/k})
\end{align}
where $e(T_{X/k})$ is the Euler class of the tangent bundle of $X$ (\cite[Def. 3.1.2]{DJK}). If $\mathbb{E}$ is an \emph{oriented motivic spectrum} (for example $\mathbf{H}\mathbb{Z}$), then $e(T_{X/k})$ agrees with the top Chern class in the $\mathbb{E}$-cohomology of $X$ (\cite[4.4.3]{DJK}).

\subsection{}
(\cite[Thm. 2.4.9]{BD}, \cite[Thm. 3.1.1]{EK})
For $f:X\to k$, denote by $\mathbb{D}_{X/k}$ the \textbf{local duality functor} or \textbf{Verdier duality functor}
\begin{align}
\label{eq:Verdual}
\begin{split}
\mathbb{D}_{X/k}:\mathbf{SH}_c(X)&\to\mathbf{SH}_c(X)\\
M&\mapsto \underline{Hom}(M,f^!\mathbbold{1}_k).
\end{split}
\end{align}
Then the \textbf{Verdier duality} holds, which states that after inverting $p$, the following canonical map of functors 
is invertible:
\begin{align}
1\to\mathbb{D}_{X/k}\circ\mathbb{D}_{X/k}.
\end{align}

\begin{lemma}
\label{lm:CXdual}
If $K\in \mathbf{SH}_c(X)$ is a constructible motivic spectrum, then we have
\begin{align}
C_X(K)=C_X(\mathbb{D}_{X/k}(K)).
\end{align}
\end{lemma}
\proof
We use a variant of the construction in \cite[Construction 2.6]{LZ}: one define the $2$-category of \emph{cohomological correspondences} $\mathcal{C}_k$ over $k$, where objects are pairs $(X,K)$ with $X$ a scheme and $K\in \mathbf{SH}(X)$, and morphisms and compositions of morphisms are described in loc. cit.. The $2$-category $\mathcal{C}_k$ is symmetric monoidal, and by \cite[Lemma 2.14, Thm. 2.16]{LZ} and the K\"unneth formulas in motivic homotopy (\cite[Thm. 2.4.6]{JY}), any pair $(X,K)$ with $K\in \mathbf{SH}_c(X)$ constructible is a dualizable object in $\mathcal{C}_k$ whose dual is $\mathbb{D}_{X/k}(K)$. By \cite[Construction 1.6]{LZ}, the class $C_X(K)$ can be interpreted as the categorical Euler characteristic of the pair $(X,K)$ in $\mathcal{C}_k$. The result then follows from the fact that the categorical Euler characteristic of a dualizable object in a symmetric monoidal category agrees with that of its dual (see \ref{num:dual}).
\endproof

\subsection{}
We now lift the characteristic class to limit Borel-Moore homology theories (Definition~\ref{def:limBM}) using the category of compactifications in~\ref{num:comp}:
\begin{definition}
\label{df:proclass}
Let $K\in \mathbf{SH}_c(X)$ be a constructible motivic spectrum.
To every compactification $X\xrightarrow{j}\overline{X}\in\operatorname{Cpt}(X/k)$, we associate two elements 
\begin{align}
(C_{\overline{X}}(j_*K)) \textrm{ and } (C_{\overline{X}}(j_!K))
\end{align}
in $\mathbb{E}^{\mathrm{BM}}_0(\overline{X}/k)$. By~\ref{num:limBM} and~\eqref{eq:CXprop}, the formation of these two families of elements associated to every compactification of $f$ determines two elements
\begin{align}
lC_X(K) \textrm{ and } lC^{\mathrm{c}}_X(K)
\end{align}
in $l\mathbb{E}^{\mathrm{BM}}_0(X/k)$, which we call the \textbf{pro-characteristic class} and \textbf{compactly supported pro-characteristic class} (or \emph{limit characteristic class} and \emph{compactly supported limit characteristic class}) respectively. 
\end{definition}

\subsection{}
Since the local duality functor~\eqref{eq:Verdual} exchanges the functors $f_*$ and $f_!$, we deduce from Lemma~\ref{lm:CXdual} the following relation between the two classes $lC_X$ and $lC^{\mathrm{c}}_X$:
\begin{corollary}
\label{cor:Ccdual}
If $K\in \mathbf{SH}_c(X)$ is a constructible motivic spectrum, then we have
\begin{align}
lC_X(K)=lC^{\mathrm{c}}_X(\mathbb{D}_{X/k}(K)).
\end{align}
\end{corollary}

\subsection{}
\label{num:Cpf}
By~\ref{num:limBM} and~\eqref{eq:CXprop}, we have the following property: if $f:Y\to X$ is a morphism of schemes, the map $f_*:l\mathbb{E}^{\mathrm{BM}}_n(Y/S)\to l\mathbb{E}^{\mathrm{BM}}_n(X/S)$ in~\eqref{eq:funclimBM} satisfies
\begin{align}
\label{eq:Cpf}
f_*lC_Y(K)=lC_X(f_*K)
\end{align}
and
\begin{align}
\label{eq:Ccpf}
f_*lC^{\mathrm{c}}_Y(K)=lC^{\mathrm{c}}_X(f_!K).
\end{align}
In particular, for $f:X\to k$, the map $f_*:l\mathbb{E}^{\mathrm{BM}}_0(X/k)\xrightarrow{\eqref{eq:funclimBM}}l\mathbb{E}^{\mathrm{BM}}_0(k/k)=[\mathbbold{1}_k,\mathbb{E}]$ sends $lC_X(K)$ (resp. $lC^{\mathrm{c}}_X(K)$) to the $\mathbb{E}$-valued Euler characteristic of $f_*K$ (resp. the $\mathbb{E}$-valued Euler characteristic of $f_!K$), see \cite[5.3.1]{JY}.

\begin{remark}
\label{rem:Apf}
The formulas~\eqref{eq:Cpf} and~\eqref{eq:Ccpf} are formulas related to push-forward maps, expressed in terms of (motivic) sheaves. Over a field of characteristic $0$, we will translate these relations geometrically in terms of constructible functions and give a refinement of \cite[Thm. 5.2]{Alu}, see~\ref{num:consf} below.
\end{remark}

\subsection{}
The canonical map
$
l\mathbb{E}^{\mathrm{BM}}_0(X/k)
\xrightarrow{\eqref{eq:limBMBM}}
\mathbb{E}^{\mathrm{BM}}_0(X/k)
$
sends both classes $lC_X(K)$ and $lC^{\mathrm{c}}_X(K)$ to the class $C_X(K)$ in~\eqref{eq:CX}. This follows from~\ref{num:limBM} and the fact that the class $C_X(K)$ is compatible with pullbacks by \'etale morphisms (\cite[Rem. 3.1.7]{JY}). In particular, If $X$ is proper over $k$, then 
under the isomorphism~\eqref{eq:limBMBM} we have identifications
\begin{align}
lC_X(K)=lC^{\mathrm{c}}_X(K)=C_X(K)\in\mathbb{E}^{\mathrm{BM}}_0(X/k).
\end{align}

\subsection{}
\label{num:Cadd}
We deduce from~\ref{num:adddt} the following property: for a distinguished triangle $K\to\mathbb{G}\to\mathbb{H}$ of constructible objects in $\mathbf{SH}_c(X)$,
we have
\begin{align}
lC_X(\mathbb{G})=lC_X(K)+lC_X(\mathbb{H}),
\end{align}
\begin{align}
\label{eq:addCc}
lC^{\mathrm{c}}_X(\mathbb{G})=lC^{\mathrm{c}}_X(K)+lC^{\mathrm{c}}_X(\mathbb{H}).
\end{align}

\subsection{}
We now compare the characteristic classes with or without compact support $lC_X(K)$ and $lC^{\mathrm{c}}_X(K)$ in the case $\mathbb{E}=\mathbf{H}\mathbb{Z}[1/p]$.

\begin{corollary}
\label{cor:CXdual}
For $\mathbb{E}=\mathbf{H}\mathbb{Z}[1/p]$, if $f:Y\to X$ is a morphism between separated $k$-schemes of finite type, and if $K\in\mathbf{SH}_c(Y)$ is a constructible motivic spectrum, then we have
\begin{align}
\label{eq:CX*!}
C_X(f_*K,\mathbf{H}\mathbb{Z}[1/p])=C_X(f_!K,\mathbf{H}\mathbb{Z}[1/p])\in CH_0(X)[1/p].
\end{align}
\end{corollary}
\proof
We know that there is a premotivic adjunction $\mathbf{SH}\to\mathbf{DM}$ realizing $\mathbf{DM}(k)[1/p]$ as the category of modules over the commutative ring spectrum $\mathbf{H}\mathbb{Z}[1/p]$, and therefore we may assume that $K\in\mathbf{DM}_c(Y)[1/p]$. By \cite[5.1.12]{JY} the map $C_X(-,\mathbf{H}\mathbb{Z}[1/p])$ are invariant under Tate twists, and therefore factors through the group $H(X)$. The result then follows from Proposition~\ref{prop:Lau1.1}.
\endproof

\subsection{}
For \'etale sheaves property~\eqref{eq:CX*!} follows from \cite[Th. 1.1]{Lau}, and reflects intrinsic nature of the \emph{characteristic cycle}, see \cite[Lemma 5.13]{Sai}. For a topological explanation in characteristic $0$, see \cite[Note 13 to Chapter 4, p.141-142]{Fult}. 

\begin{corollary}
\label{cor:CXcpt}
For $\mathbb{E}=\mathbf{H}\mathbb{Z}[1/p]$ and any object $K\in\mathbf{SH}_c(X)$, we have 
\begin{align}
\label{eq:lCX*!}
lC_X(K,\mathbf{H}\mathbb{Z}[1/p])
=
lC^{\mathrm{c}}_X(K,\mathbf{H}\mathbb{Z}[1/p])
\in
lCH_0(X)[1/p].
\end{align}
\end{corollary}

\subsection{}
Note that in the quadratic setting, formulas~\eqref{eq:CX*!} and~\eqref{eq:lCX*!} do not hold for $\mathbb{E}=\mathbf{H}_{\mathrm{MW}}\mathbb{Z}[1/p]$, see Remark~\ref{rem:levquad} below.

\subsection{}
We now introduce a variant of the characteristic class.
Let $U\xrightarrow{j}X\xrightarrow{f}k$ be two composable morphisms of schemes, with $j$ an open immersion. Let $\delta_X:X\to X\times_kX$ be the diagonal morphism. Let $K\in\mathbf{SH}_c(U)$ be a constructible motivic spectrum.  
Consider the composition
\begin{align}
\label{eq:CXU}
\begin{split}
C_{X,U}(K):
\mathbbold{1}_X
&\to
\underline{Hom}(j_!K,j_!K)
\simeq
\delta_X^!(\mathbb{D}(j_!K)\boxtimes_kj_!K)
\to
\delta_X^*(\mathbb{D}(j_!K)\boxtimes_kj_!K)\\
&\simeq
j_!(\mathbb{D}(K)\otimes_kK)
\to
j_!j^!f^!\mathbbold{1}_k
\end{split}
\end{align}
(see also \cite[Def. 2.1.8]{AS}). 
If $\mathbb{E}\in \mathbf{SH}(k)$ is a motivic spectrum with a map $\mathbbold{1}_k\to \mathbb{E}$, we denote by $C_{X,U}(K,\mathbb{E})$ the composition 
\begin{align}
\label{eq:CXUHom}
C_{X,U}(K,\mathbb{E}):
\mathbbold{1}_X
\xrightarrow{C_{X,U}(K)}
j_!j^!f^!\mathbbold{1}_k
\xrightarrow{}
j_!j^!f^!\mathbb{E}.
\end{align}
\subsection{}
\label{num:propequal}
If $U=X$, then $C_{X,U}^{}(K,\mathbb{E})=C_U(K,\mathbb{E})\in\mathbb{E}^{\mathrm{BM}}_0(U/k)$. 


\begin{lemma}
\label{lm:addhom}
For a distinguished triangle $L\to M\to N$ of constructible objects in $\mathbf{SH}_c(X)$, the following equality holds: 
\begin{align}
C_{X,U}(M,\mathbb{E})=C_{X,U}(L,\mathbb{E})+C_{X,U}(N,\mathbb{E}).
\end{align}
\end{lemma}
\proof
The claim is a form of additivity of Euler characteristics, which can be proved in the style of \cite{May}, \cite{GPS}. We proceed as in the proof of \cite[Prop. 4.2.6]{JY}: using the language of higher categories, there exist objects $u,u',v,w$ and a commutative diagram of the form
\begin{align}
\begin{split}
\xymatrix@=8pt{
    & \mathbbold{1}_X \ar^-{}[d] \ar^-{}[rd] \ar_-{}[ld] &\\
   \underline{Hom}(j_!L,j_!L)\oplus\underline{Hom}(j_!N,j_!N) \ar^-{\wr}[dd] \ar^-{}[rd]
   & v \ar^-{}[r] \ar^-{}[l] 
   & \underline{Hom}(j_!M,j_!M) \ar^-{}[ld] \ar_-{\wr}[dd]\\
   & w \ar@{.>}^-{}[d]
   &\\
   \delta_X^!(\mathbb{D}(j_!L)\boxtimes_kj_!L\oplus\mathbb{D}(j_!N)\boxtimes_kj_!N) \ar^-{}[r]  \ar^-{}[d]
   & \delta_X^!u \ar^-{}[d]
   & \delta_X^!(\mathbb{D}(j_!M)\boxtimes_kj_!M) \ar^-{}[l] \ar^-{}[d]\\
   \delta_X^*(\mathbb{D}(j_!L)\boxtimes_kj_!L\oplus\mathbb{D}(j_!N)\boxtimes_kj_!N) \ar^-{}[r]  \ar^-{\wr}[d]
   & \delta_X^*u \ar^-{}[d]
   & \delta_X^*(\mathbb{D}(j_!M)\boxtimes_kj_!M) \ar^-{}[l] \ar_-{\wr}[d]\\
   j_!(\mathbb{D}(L)\otimes_kL\oplus\mathbb{D}(N)\otimes_kN) \ar^-{}[r]  \ar^-{}[rd]
   & u' \ar^-{}[d]
   & j_!(\mathbb{D}(M)\otimes_kM) \ar^-{}[l] \ar^-{}[ld]\\
   & j_!f^!\mathbbold{1}_k &
 }
 \end{split}
\end{align}
which finishes the proof.
\endproof

\begin{lemma}
\label{lm:CHomfun}
Consider a commutative diagram of schemes
\begin{align}
\label{eq:comprell}
\begin{split}
  \xymatrix@=10pt{
    U \ar[r]^-{k} \ar[d]_-{f} & X \ar[d]^-{g} \\
    V \ar[r]^-{j} & Y.
  }
\end{split}
\end{align}
where both horizontal maps are open immersions and $g$ is proper. Let $K\in\mathbf{SH}_c(U)$ be a constructible motivic spectrum. Then the following diagram is commutative:
\begin{align}
\begin{split}
  \xymatrix@=10pt{
    \mathbbold{1}_{Y} \ar[r]^-{} \ar[d]_-{C_{Y,V}(f_!K)} & g_*\mathbbold{1}_{X} \ar[d]^-{g_*C_{X,U}(K)}  \\ 
    j_!j^!q^!\mathbbold{1}_k & g_*k_!f^!j^!q^!\mathbbold{1}_k. \ar[l]^-{} \\ 
  }
\end{split}
\end{align}

\end{lemma}
\proof
Let $q:Y\to k$ be the structure morphism.
The result follows from the following commutative diagram:
\begin{align}
\begin{split}
  \xymatrix@=10pt{
    \mathbbold{1}_{Y} \ar[r]^-{} \ar[d]_-{} & g_*\mathbbold{1}_{X} \ar[d]_-{}  \\ 
    \underline{Hom}(j_!f_!K,j_!f_!K) \ar[d]^-{\wr} & g_*\underline{Hom}(k_!K,k_!K) \ar[l]^-{} \ar[d]^-{\wr} \\ 
    \delta_Y^!(\mathbb{D}(j_!f_!K)\boxtimes_kj_!f_!K) \ar[d]^-{} & g_*\delta_X^!(\mathbb{D}(k_!K)\boxtimes_kk_!K) \ar[l]^-{} \ar[d]^-{} \\ 
    \delta_Y^*(\mathbb{D}(j_!f_!K)\boxtimes_kj_!f_!K) \ar[r]^-{} \ar[d]^-{\wr} & g_*\delta_X^*(\mathbb{D}(k_!K)\boxtimes_kk_!K) \ar[d]^-{\wr} \\ 
    j_!(\mathbb{D}(f_!K)\otimes_kf_!K) \ar[r]^-{} \ar[d]^-{} & g_*k_!(\mathbb{D}(K)\otimes_kK) \ar[d]^-{} \\ 
    j_!j^!q^!\mathbbold{1}_k & g_*k_!f^!j^!q^!\mathbbold{1}_k \ar[l]^-{} \\ 
  }
\end{split}
\end{align}
\endproof

\begin{corollary}
\label{cor:compim}
Let $j:U\to X$ be an open immersion and let $K\in\mathbf{SH}_c(U)$ be a constructible motivic spectrum. Then the image of $C_{X,U}(K,\mathbb{E})$ in $\mathbb{E}^{\mathrm{BM}}_0(X/k)$ is $C_{X}(j_!K,\mathbb{E})$. 

\end{corollary}
\proof
Apply Lemma~\ref{lm:CHomfun} to the commutative diagram
\begin{align}
\label{eq:diagmorcomp}
\begin{split}
  \xymatrix@=6pt{
    U \ar[r]^-{j} \ar[d]_-{j} & X \ar@{=}[d]^-{} \\ 
     X \ar@{=}[r]^-{} & X.
  }
\end{split}
\end{align}
and use~\ref{num:propequal}.
\endproof

\subsection{}
\label{num:rellog}
If $X$ is smooth and proper of dimension $n$, $D=X-U$ is a simple normal-crossing divisor in $X$ and $K=\mathbbold{1}_U$, then the class $C_{X,U}(\mathbbold{1}_U,\mathbf{H}\mathbb{Z})$ agrees, up to a sign, with the \emph{relative top Chern class} of the sheaf of differential forms with logarithmic poles $\Omega^1_X(\operatorname{log}D)$ defined in \cite{Sai2}. Note that the construction in loc. cit. also works in the motivic setting, which we briefly recall here.

Let $D=\cup_{i\in I}D_i$ be the components of $D$, and for $J\subset I$ let $D_J=\cap_{i\in J}D_i$. The sheaf $\mathcal{E}=\Omega^1_X(\operatorname{log}D)$ is a locally free $\mathcal{O}_X$-module of rank $n$, and there are morphisms of $\mathcal{O}_{D_i}$-modules called \textbf{residue morphisms} $\operatorname{res}_i:\mathcal{E}_{|D_i}\to\mathcal{O}_{D_i}$ such that for all $J\subset I$, the map $\oplus_{i\in J}\operatorname{res}_i:\mathcal{E}_{D_J}\to\mathcal{O}_{D_J}^J$ is surjective (\cite[3.7.4]{Del}). Let $V=\operatorname{Spec}_X(\operatorname{Sym}^*(\mathcal{E}^\vee))$ be the vector bundle on $X$ associated to the dual sheaf of $\mathcal{E}$. The map $\operatorname{res}_i$ induces a morphism $r_i:V_{D_i}\to \mathbb{A}^1_{1_{D_i}}$, and let $\Delta_i=r_i^{-1}(D_i)$ be the closed subscheme of $V_{D_i}$, where $1_{D_i}$ is the image of the $1$-section of $\mathbb{A}^1_{D_i}$. Then $\Delta=\cup_{i\in I}\Delta_i$ is contained in the complement of the image of the $0$-section of $V$, and the $0$-section $0:X\to V$ factors as
\begin{align}
X\xrightarrow{l}V-\Delta\xrightarrow{k}V
\end{align}
with $l$ a regular closed immersion. We have a map
\begin{align}
\label{eq:fdl0}
\operatorname{Th}(-T_{V/X})
\to 
0_*0^*\operatorname{Th}(-T_{V/X})
\simeq 
0_*\operatorname{Th}(-N_l)
\underset{\sim}{\xrightarrow{\eta_l}}
0_*l^!\mathbbold{1}_{V-\Delta}
\to
k_!\mathbbold{1}_{V-\Delta}
\end{align}
where 
\begin{enumerate}
\item
$\operatorname{Th}$ is the Thom transformation (\cite[2.1.4]{DJK});
\item
$T_{V/X}$ is the relative tangent bundle of the projection $\pi:V\to X$;
\item
$N_l$ is the normal bundle of the regular immersion $l$;
\item
$\eta_l:\operatorname{Th}(-N_l)\to l^!\mathbbold{1}_{V-\Delta}$ is the fundamental class of $l$ (\cite[Thm. 3.3.2]{DJK}).
\end{enumerate}
We define the relative top Chern class of the sheaf of differential forms with logarithmic poles as the composition
\begin{align}
\begin{split}
c_{X,D}(\Omega^1_{X/k}(\operatorname{log}D)):
\mathbbold{1}_X
&\to
\pi_*\mathbbold{1}_V
\simeq
\pi_*(\pi^!\mathbbold{1}_X\otimes\operatorname{Th}(-T_{V/X}))\\
&\xrightarrow{\eqref{eq:fdl0}}
\pi_*(\pi^!\mathbbold{1}_X\otimes k_!\mathbbold{1}_{V-\Delta})
\simeq
\pi_*\pi^*j_!\operatorname{Th}(T_{V-\Delta/X})\\
&\simeq
j_!\operatorname{Th}(T_{V-\Delta/X})
\to
j_!(f^*\mathbf{H}\mathbb{Z}\otimes\operatorname{Th}(T_{V-\Delta/X}))\\
&\simeq
j_!f^*\mathbf{H}\mathbb{Z}(n)[2n]
\simeq
j_!f^!\mathbf{H}\mathbb{Z}.
 \end{split}
\end{align}
Here $1\simeq\pi_*\pi^*$ is the homotopy invariance, and $f^*\mathbf{H}\mathbb{Z}\otimes\operatorname{Th}(T_{V-\Delta/X})\simeq f^*\mathbf{H}\mathbb{Z}(n)[2n]$ follows from the orientation on $\mathbf{H}\mathbb{Z}$. Then we have $C_{X,U}(\mathbbold{1}_U,\mathbf{H}\mathbb{Z})=(-1)^nc_{X,D}(\Omega^1_{X/k}(\operatorname{log}D))$.

\subsection{}
If $X$ is a compactification of $U$, then the class $C_{X,U}(K,\mathbb{E})$ in~\eqref{eq:CXUHom} defines a class in the $\mathbb{E}$-homology
\begin{align}
C_{X,U}(K,\mathbb{E})\in \mathbb{E}_0(U/k)[1/p].
\end{align}

\begin{corollary}
\label{cor:Homind}
For $X$ a compactification of $U$, the class $C_{X,U}(K,\mathbb{E})$ is independent of the choice of the compactification $j:U\to X\in\operatorname{Cpt}(X/k)$.
\end{corollary}
\proof
If $p:X\to X'$ is a morphism of compactifications of $U$, we have the following commutative diagram
\begin{align}
\label{eq:diagmorcomp}
\begin{split}
  \xymatrix@=6pt{
    U \ar[r]^-{}  \ar@{=}[d]^-{} & X \ar[d]^-{p}\\ 
     U \ar[r]^-{} & X'.
  }
\end{split}
\end{align}
Lemma~\ref{lm:CHomfun} applied to diagram~\eqref{eq:diagmorcomp} shows that $C_{X,U}^{}(K,\mathbb{E})=C_{X',U}^{}(K,\mathbb{E})\in\mathbb{E}_0(U/k)[1/p]$, and the result follows.
\endproof

\begin{remark}
If we replace $j_!K$ by $j_*K$ in~\eqref{eq:CXU}, the composition
\begin{align}
\begin{split}
\mathbbold{1}_X
&\to
\underline{Hom}(j_*K,j_*K)
\simeq
\delta^!(\mathbb{D}(j_*K)\boxtimes_kj_*K)
\to
\delta^*(\mathbb{D}(j_*K)\boxtimes_kj_*K)\\
&\simeq
j_!(\mathbb{D}(K)\otimes_kK)
\to
j_!f^!\mathbbold{1}_k
\to
j_!f^!\mathbb{E}
\end{split}
\end{align}
is also independent of the compactification $X$. Similar to the discussion in Corollary~\ref{cor:Ccdual}, Corollary~\ref{cor:CXdual} and Remark~\ref{rem:levquad}, this class is dual to the class $C_{X,U}(K,\mathbb{E})$; for $\mathbb{E}=\mathbf{H}\mathbb{Z}[1/p]$ the two agree, while for $\mathbb{E}=\mathbf{H}_{\mathrm{MW}}\mathbb{Z}[1/p]$ the two are different.

\end{remark}


\begin{definition}
\label{def:CHom}
By Corollary~\ref{cor:Homind}, the class $C_{X,U}(K,\mathbb{E})\in\mathbb{E}^{}_0(U/k)[1/p]$ is independent of $X$, and we denote it by $C^{\operatorname{Hom}}_U(K,\mathbb{E})$ or simply $C^{\operatorname{Hom}}_U(K)$, called the ($\mathbb{E}$-valued) \textbf{homological characteristic class} of $K$. 
\end{definition}

\subsection{}
Lemma~\ref{lm:addhom} implies that the homological characteristic class induces a map
\begin{align}
\label{eq:CHom}
\begin{split}
C^{\operatorname{Hom}}_U:K_0(\mathbf{SH}_c(U))&\to\mathbb{E}^{}_0(U/k)[1/p]\\
K&\mapsto C^{\operatorname{Hom}}_U(K,\mathbb{E}).
\end{split}
\end{align}

\subsection{}
This construction applies to $\mathbb{E}=\mathbf{H}\mathbb{Z}$ as well as $\mathbf{H}^{\mathrm{MW}}\mathbb{Z}$, the latter being a quadratic refinement of the former.

\subsection{}
If $k$ is a perfect field which satisfies embedded resolution of singularities in the sense of~\ref{num:RS}, then the class $C^{\operatorname{Hom}}_U(K)$ can be defined in $\mathbb{E}_0^{}(U/k)$ without inverting $p$.

\begin{corollary}
\label{cor:CHomfun}
Let $f:U\to V$ be a morphism and let $K\in\mathbf{SH}_c(U)$ be a constructible motivic spectrum. Then the map $f_*:\mathbb{E}_0(U/k)\to\mathbb{E}_0(V/k)$ in~\eqref{eq:hompf} sends $C_{U}^{\operatorname{Hom}}(K)$ to $C_{V}^{\operatorname{Hom}}(f_!K)$.
\end{corollary}
\proof
Let $j:V\to\overline{V}$ be a compactification. As in~\eqref{eq:comprel}, one may choose a commutative diagram 
\begin{align}
\label{eq:comprell}
\begin{split}
  \xymatrix@=10pt{
    U \ar[r]^-{k} \ar[d]_-{f} & \overline{U} \ar[d]^-{\overline{f}} \\
    V \ar[r]^-{j} & \overline{V}.
  }
\end{split}
\end{align}
where $k$ is a compactification and $\overline{f}$ is proper. The result then follows by applying Lemma~\ref{lm:CHomfun} to the diagram~\eqref{eq:comprell}.
\endproof


\begin{corollary}
\label{cor:CHomlift}
Let $U$ be a scheme and let $K\in\mathbf{SH}_c(U)$ be a constructible motivic spectrum. Then the map $\mathbb{E}^{}_0(U/k)\to l\mathbb{E}^{BM}_0(U/k)$ in~\eqref{eq:HomlimBM} sends $C_{U}^{\operatorname{Hom}}(K)$ to $lC^{\mathrm{c}}_{U}(K)$. 

\end{corollary}
\proof 
Follows from~\ref{num:limBM}, Definition~\ref{df:proclass} and Corollary~\ref{cor:compim}.
\endproof


\subsection{}
We now summarize the results in this section.
Let $U$ be a scheme and let $K\in\mathbf{SH}_c(U)$ be a constructible motivic spectrum. There are three characteristic classes in different homology theories:
\begin{enumerate}
\item
(Definition~\ref{def:CHom}) $C_{U}^{\operatorname{Hom}}(K)\in\mathbb{E}^{}_0(U/k)[1/p]$;
\item
(Definition~\ref{df:proclass}) $lC^{\mathrm{c}}_{U}(K)\in l\mathbb{E}^{BM}_0(U/k)[1/p]$;
\item
(\ref{num:CX}) $C_{U}^{}(K)\in\mathbb{E}^{BM}_0(U/k)[1/p]$.
\end{enumerate}
If $k$ is a perfect field which satisfies embedded resolution of singularities in the sense of~\ref{num:RS}, then these classes can be defined without inverting $p$.
These classes are compatible with the caonical maps as illustrated in the following diagram:
\begin{align}
\begin{split}
  \xymatrix@R=4pt{
    \mathbb{E}^{}_0(U/k) \ar[r]^-{\eqref{eq:HomlimBM}}  & l\mathbb{E}^{BM}_0(U/k) \ar[r]^-{\eqref{eq:limBMBM}} & \mathbb{E}^{BM}_0(U/k)\\
    C_{U}^{\operatorname{Hom}}(K) \ar@{(-}[u] \ar@{|->}[r]^-{} & lC^{\mathrm{c}}_{U}(K) \ar@{(-}[u] \ar@{|->}[r]^-{} & C_{U}^{}(K). \ar@{(-}[u]
  }
\end{split}
\end{align}
These classes satisfy the following functorialities: if $f:U\to V$ is a morphism of schemes, then
\begin{enumerate}
\item
(Corollary~\ref{cor:CHomfun}) the map $f_*:\mathbb{E}_0(U/k)\xrightarrow{\eqref{eq:hompf}}\mathbb{E}_0(V/k)$ sends $C_{U}^{\operatorname{Hom}}(K)$ to $C_{V}^{\operatorname{Hom}}(f_!K)$;
\item
(\ref{num:Cpf}) the map $f_*:l\mathbb{E}^{BM}_0(U/k)\xrightarrow{\eqref{eq:funclimBM}}l\mathbb{E}^{BM}_0(V/k)$ sends $lC^{\mathrm{c}}_{U}(K)$ to $lC^{\mathrm{c}}_{V}(f_!K)$;
\item
(\ref{num:CXprop})
if $f$ is proper, the map $f_*:\mathbb{E}^{BM}_0(U/k)\xrightarrow{\eqref{eq:BMprop}}\mathbb{E}^{BM}_0(V/k)$ sends $C^{}_{U}(K)$ to $C^{}_{V}(f_!K)$.
\end{enumerate}

\section{Relation with the pro-CSM class}
\subsection{}
In this section, unless otherwise mentioned, we focus on the case $\mathbb{E}=\mathbf{H}_{}\mathbb{Z}$.
Recall Aluffi's construction of the pro-CSM class:
\begin{theorem}[\textrm{\cite[Prop. 4.3]{Alu}}]
\label{thm:Alu43}
Assume that $k$ is a perfect field which satisfies resolution of singularities by weak factorizations (as in \cite[Thm. 0.1.1]{AKMW}, which is stronger than embedded resolution of singularities in the sense of~\ref{num:RS}). Then there is a unique way to associate a class 
$lC^{SM}_*(X)\in lCH_*(X)$ for every $k$-scheme $X$, 
such that
\begin{enumerate}
\item
(Additivity for stratifications)
If $X$ has a finite stratification into locally closed subschemes $X=\cup_{i=1}^nU_i$, with each $U_i$ smooth and irreducible, then we have
\begin{align}
\label{eq:addcsm}
lC^{SM}_*(X)=\sum_{i=1}^n\iota_{i*}lC^{SM}_*(U_i)\in lCH_*(X),
\end{align}
where $\iota_{i}:U_i\to X$ is the inclusion.
\item
If $X=\overline{X}-D\xrightarrow{j}\overline{X}$ with $\overline{X}$ smooth and proper 
and $D$ a simple normal-crossing divisor in $\overline{X}$, then the image of $lC^{SM}_*(X)$ in $CH_*(\overline{X})$ is given by the top Chern class of the dual of the sheaf of differential forms with logarithmic poles along $D$ (\cite[Def. 3.1]{Del}, \cite[Def. 1.2]{Sai1}):
\begin{align}
j_*lC^{SM}_*(X)=c(\Omega^1_{\overline{X}/k}(\operatorname{log}D)^\vee)\cap[\overline{X}]\in lCH_*(\overline{X})\simeq CH_*(\overline{X}).
\end{align}
\end{enumerate}
\end{theorem}

\subsection{}
Consider the zero-dimensional part $lC^{SM}_0(X)\in lCH_0(X)$ of the pro-CSM class $lC^{SM}(X)$. We now relate this class with our compactly supported pro-characteristic class $lC^{\mathrm{c}}_X(\mathbbold{1}_X,\mathbf{H}\mathbb{Z})$ defined in Definition~\ref{df:proclass}:
\begin{theorem}
\label{th:A=C}
Assume that the conditions in Theorem~\ref{thm:Alu43} are satisfied.
Let $X$ be a scheme. Then via the isomorphism $lCH_0(X)\overset{\eqref{eq:limBMHZ}}{\simeq} l\mathbf{H}\mathbb{Z}^{\mathrm{BM}}(X/k)$, we have 
\begin{align}
lC^{SM}_0(X)=lC^{\mathrm{c}}_X(\mathbbold{1}_X,\mathbf{H}\mathbb{Z})\in lCH_0(X).
\end{align}
\end{theorem}

\subsection{}
As preparation for the proof, we now extend~\eqref{eq:JY5116} to simple normal-crossing divisors (see also \cite[Cor. 3.23]{Tsu}):
\begin{lemma}
Let $X$ be a smooth $k$-scheme of dimension $n$, and let $D$ be a simple normal-crossing divisor in $X$. Denote by $c_{l,X}^D$ the localized Chern class (\cite[\S 18.1]{Ful}, \cite{Blo}). Then
\begin{align}
\label{eq:CD1D}
C_D(\mathbbold{1}_D,\mathbf{H}\mathbb{Z})
=
(-1)^{n-1}\sum_{l=1}^{n}c_{n-l}(\Omega^1_{X/k})_{|D}
\cdot
c_{l,X}^D(\Omega^1_{X/k}\to\Omega^1_{X/k}(\operatorname{log}D))\cap[X]
\in CH_0(D).
\end{align}
where $\Omega^1_{X/k}\to\Omega^1_{X/k}(\operatorname{log}D)$ is a two-term complex with $\Omega^1_{X/k}(\operatorname{log}D)$ placed in degree $0$. 

\end{lemma}
\proof
We use the computation of the localized Chern classes of the complex $\Omega^1_{X/k}\to\Omega^1_{X/k}(\operatorname{log}D)$ as in the proof of \cite[Thm. 3.1]{Sil}. Let $D_1,\cdots,D_m$ be the irreducible components of $D$. We have a short exact sequence of sheaves (\cite[3.3.2.2]{Del})
\begin{align}
0\to \Omega^1_{X/k}
\to \Omega^1_{X/k}(\operatorname{log}D)
\to \oplus_{i=1}^m\mathcal{O}_{D_i}\to 0
\end{align}
from which we deduce an equality of total localized Chern classes
\begin{align}
\label{eq:eqtotc}
c_X^D(\Omega^1_{X/k}\to\Omega^1_{X/k}(\operatorname{log}D))
=
\prod_{i=1}^mc_X^D(\mathcal{O}_{D_i}).
\end{align}
From the short exact sequence
\begin{align}
0\to 
\mathcal{O}_X(-D_i)
\to
\mathcal{O}_X
\to
\mathcal{O}_{D_i}\to0
\end{align}
we deduce 
\begin{align}
1=c_X^D(\mathcal{O}_X)=c_X^D(\mathcal{O}_{D_i})c_X^D(\mathcal{O}_X(-D_i))=c_X^D(\mathcal{O}_{D_i})(1+c_{1,X}^D(-D_i)),
\end{align}
and therefore $c_{l,X}^D(\mathcal{O}_{D_i})=c_{1,X}^D(D_i)^l$. Since taking the dual of a bundle changes its $s$-th (localized) Chern class by a sign $(-1)^s$, we are reduced to show that
\begin{align}
\begin{split}
\label{eq:silclaim}
-C_{D}(\mathbbold{1}_D,\mathbf{H}\mathbb{Z})
=
\sum_{l=1}^n
\sum_{j_1+\cdots+j_m=l}
c_{n-l}(T_{X/k})_{|D}\cdot c_{1,X}^D(-D_1)^{j_1}\cdots c_{1,X}^D(-D_m)^{j_m}\cap[X]
\in CH_0(D).
\end{split}
\end{align}
For each $i$, we have the short exact sequence
\begin{align}
0\to
N_{D_i}X
\to
T_{X|D_i}
\to
T_{D_i}
\to0
\end{align}
and the isomorphism $N_{D_i}X\simeq\mathcal{O}_X(D_i)_{|D_i}$, we deduce
\begin{align}
c(T_{X})_{|D_i}
=
c(N_{D_i}X)c(T_{D_i})
=
c(\mathcal{O}_X(D_i))_{|D_i}c(T_{D_i})
=
(1+c_1(D_i))c(T_{D_i})
\end{align}
and consequently
\begin{align}
\label{eq:sil2}
c_d(T_{X})_{|D_i}
=
c_d(T_{D_i})
-
c_{d-1}(T_{D_i})c_1(-D_i)_{|D_i}.
\end{align}
We now prove~\eqref{eq:silclaim} by induction on the number of irreducible components $m$.
\begin{itemize}
\item
If $m=1$, $D=D_1$ is smooth, in which case we have
\begin{align}
\begin{split}
&-\sum_{l=1}^nc_{n-l}(T_{X/k})c_{1,X}^{D_1}(-D_1)^{l}\cap[X]
=
\sum_{l=1}^nc_{n-l}(T_{X/k})_{|D_1}c_1(-D_1)_{|D_1}^{l-1}\cap[D_1]\\
\overset{\eqref{eq:sil2}}{=}
&\sum_{l=1}^nc_{n-l}(T_{D_1})c_1(-D_1)_{|D_1}^{l-1}\cap[D_1]
-
\sum_{l=1}^{n-1}c_{n-l-1}(T_{D_1})c_1(-D_1)_{|D_1}^l\cap[D_1]\\
=
&\sum_{l=0}^{n-1}c_{n-l-1}(T_{D_1})c_1(-D_1)_{|D_1}^{l}\cap[D_1]
-
\sum_{l=1}^{n-1}c_{n-l-1}(T_{D_1})c_1(-D_1)_{|D_1}^l\cap[D_1]\\
=
&c_{n-1}(T_{D_1})\cap[D_1]
\overset{\eqref{eq:JY5116}}{=}
C_{D_1}(\mathbbold{1}_{D_1},\mathbf{H}\mathbb{Z}[1/p])\in CH_0(D_1)[1/p]
\end{split}
\end{align}
which proves~\eqref{eq:silclaim}. 
\item
Assume that~\eqref{eq:silclaim} holds when $D$ has at most $m-1$ branches. We split the right-hand side of~\eqref{eq:silclaim} into sums according to the last non-zero index $j_q$:
\begin{align}
\label{eq:sil21}
\begin{split}
&\sum_{l=1}^n
\sum_{j_1+\cdots+j_m=l}
c_{n-l}(T_{X/k})_{|D}c_{1,X}^{D}(-D_1)^{j_1}\cdots c_{1,X}^{D}(-D_m)^{j_m}\cap[X]\\
=
&\sum_{q=1}^m\sum_{l=1}^n
\sum_{\substack{j_1+\cdots+j_q=l\\ j_q\geqslant1}}
c_{n-l}(T_{X/k})_{|D}c_{1,X}^{D}(-D_1)^{j_1}\cdots c_{1,X}^{D}(-D_q)^{j_q}\cap[X].
\end{split}
\end{align}
For each $q$, denote by $\iota_q:D_q\to D$ the inclusion. Let $E_q=D_1\cup\cdots\cup D_{q-1}$ and denote by $\psi_q:E_q\cap D_q\to D_q$ the inclusion of the simple normal-crossing divisor with $q-1$ branches. 
We proceed as in the case $m=1$:
\begin{align}
\begin{split}
\label{eq:sil23}
&-\sum_{l=1}^n
\sum_{\substack{j_1+\cdots+j_q=l\\ j_q\geqslant1}}
c_{n-l}(T_{X/k})_{|D}c_{1,X}^{D}(-D_1)^{j_1}\cdots c_{1,X}^{D}(-D_q)^{j_q}\cap[X]\\
=
&\sum_{l=1}^n
\sum_{\substack{j_1+\cdots+j_q=l\\ j_q\geqslant1}}
c_{n-l}(T_{X/k})_{|D_q}c_1(-D_1)^{j_1}_{|D_q}\cdots c_1(-D_q)^{j_q-1}_{|D_q}\cap[D_q]\\
\overset{\eqref{eq:sil2}}{=}
&\sum_{l=1}^n
\sum_{\substack{j_1+\cdots+j_q=l\\ j_q\geqslant1}}
c_{n-l}(T_{D_m/k})c_1(-D_1)^{j_1}_{|D_q}\cdots c_1(-D_q)^{j_q-1}_{|D_q}\cap[D_q]\\
-
&\sum_{l=1}^{n-1}
\sum_{\substack{j_1+\cdots+j_q=l\\ j_q\geqslant1}}
c_{n-l-1}(T_{D_m/k})c_1(-D_1)^{j_1}_{|D_q}\cdots c_1(-D_q)^{j_q}_{|D_q}\cap[D_q]\\
=
&\sum_{l=0}^{n-1}
\sum_{\substack{j_1+\cdots+j_q=l}}
c_{n-l-1}(T_{D_m/k})c_1(-D_1)^{j_1}_{|D_q}\cdots c_1(-D_q)^{j_q}_{|D_q}\cap[D_q]\\
-
&\sum_{l=1}^{n-1}
\sum_{\substack{j_1+\cdots+j_q=l\\ j_q\geqslant1}}
c_{n-l-1}(T_{D_q/k})c_1(-D_1)^{j_1}_{|D_q}\cdots c_1(-D_q)^{j_q}_{|D_q}\cap[D_q]\\
=
&c_{n-1}(T_{D_q/k})\cap[D_q]\\
+
&\sum_{l=1}^{n-1}
\sum_{\substack{j_1+\cdots j_{q-1}=l}}
c_{n-l-1}(T_{D_q/k})c_1(-D_1)^{j_1}_{|D_q}\cdots c_1(-D_{q-1})^{j_{q-1}}_{|D_q}\cap[D_q]\\
=
&\iota_{q*}C_{D_q}(\mathbbold{1}_{D_q},\mathbf{H}\mathbb{Z})
-
\iota_{q*}\psi_{q*}C_{E_q\cap D_q}(\mathbbold{1}_{E_q\cap D_q},\mathbf{H}\mathbb{Z})
\in CH_0(D)
\end{split}
\end{align}
where in the last equality we use the induction hypothesis for the inclusion $\psi_q$. 
Combining~\eqref{eq:sil21} and~\eqref{eq:sil23}, we are reduced to show that
\begin{align}
\begin{split}
C_D(\mathbbold{1}_D,\mathbf{H}\mathbb{Z})
=
\sum_{q=1}^m
(\iota_{q*}C_{D_q}(\mathbbold{1}_{D_q},\mathbf{H}\mathbb{Z})
-
\iota_{q*}\psi_{q*}C_{E_q\cap D_q}(\mathbbold{1}_{E_q\cap D_q},\mathbf{H}\mathbb{Z}))
\in CH_0(D).
\end{split}
\end{align}
Using a telescopic sum, this amounts to show that 
\begin{align}
\label{eq:teles}
\begin{split}
&\xi_{q*}C_{E_q}(\mathbbold{1}_{E_q},\mathbf{H}\mathbb{Z})
-
\xi_{q-1,*}C_{E_{q-1}}(\mathbbold{1}_{E_{q-1}},\mathbf{H}\mathbb{Z})\\
=
&\iota_{q*}C_{D_q}(\mathbbold{1}_{D_q},\mathbf{H}\mathbb{Z})
-
\iota_{q*}\psi_{q*}C_{E_q\cap D_q}(\mathbbold{1}_{E_q\cap D_q},\mathbf{H}\mathbb{Z})
\in CH_0(D)
\end{split}
\end{align}
where $\xi_q:E_q\to D$ is the inclusion. Formula~\eqref{eq:teles} follows from additivity~\eqref{eq:adddt}, the localization triangle~\eqref{eq:locseq} and a standard Mayer-Vietoris argument, which finishes the proof.
\end{itemize}
\endproof

\begin{corollary}
\label{cor:Cj!}
Let $X$ be a smooth $k$-scheme of dimension $n$, and let $i:D\to X$ be the inclusion of a simple normal-crossing divisor, with open complement $j:U\to X$. Then
\begin{align}
\label{eq:Cj!log}
C_{X}(j_!\mathbbold{1}_U,\mathbf{H}\mathbb{Z})
=
c_n(\Omega^1_{X/k}(\operatorname{log}D)^\vee)\cap[X]
\in CH_0(X).
\end{align}
If $k$ is a perfect field which satisfies embedded resolution of singularities in the sense of~\ref{num:RS}, then~\eqref{eq:Cj!log} holds for $C_{X}(j_!\mathbbold{1}_U,\mathbf{H}\mathbb{Z})\in CH_0(X)$ without inverting $p$.
\end{corollary}
\proof
Formula~\eqref{eq:CD1D} and the Whitney sum formula imply that
\begin{align}
\begin{split}
-i_*C_D(\mathbbold{1}_D,\mathbf{H}\mathbb{Z})
&=
(-1)^{n}\sum_{j=0}^{n-1}c_j(\Omega^1_{X/k})c_{n-j}(\Omega^1_{X/k}\to\Omega^1_{X/k}(\operatorname{log}D))\cap[X]\\
&=
(-1)^{n}(c_n(\Omega^1_{X/k}(\operatorname{log}D))-c_n(\Omega^1_{X/k}))\cap[X]\\
&=
c_n(\Omega^1_{X/k}(\operatorname{log}D)^\vee)\cap[X]-c_n(T_{X/k})\cap[X]\in CH_0(X).
\end{split}
\end{align}
It follows that 
\begin{align}
\begin{split}
&C_{X}(j_!\mathbbold{1}_U,\mathbf{H}\mathbb{Z})
=
C_X(\mathbbold{1}_X,\mathbf{H}\mathbb{Z})-C_X(i_*\mathbbold{1}_D,\mathbf{H}\mathbb{Z})\\
=
&c_n(T_{X/k})\cap[X]-i_*C_D(\mathbbold{1}_D,\mathbf{H}\mathbb{Z})
=
c_n(\Omega^1_{X/k}(\operatorname{log}D)^\vee)\cap[X]\in CH_0(X)
\end{split}
\end{align}
which finishes the proof.
\endproof

\noindent
\emph{Proof of Theorem~\ref{th:A=C}.}
The class $lC^{\mathrm{c}}_X(\mathbbold{1}_X,\mathbf{H}\mathbb{Z})$ also satisfies the additivity for stratifications by additivity~\eqref{eq:adddt} 
and the localization triangle~\eqref{eq:locseq}.
By Aluffi's Theorem~\ref{thm:Alu43}, it suffices to show that if $X=\overline{X}-D\xrightarrow{j}\overline{X}$ with $\overline{X}$ smooth and proper of dimension $n$ and $D$ a simple-normal crossing divisor in $\overline{X}$, then one has
\begin{align}
C_{\overline{X}}(j_!\mathbbold{1}_X,\mathbf{H}\mathbb{Z})
=
c_n(\Omega^1_{\overline{X}/k}(\operatorname{log}D)^\vee)\cap[\overline{X}]
\in CH_0(\overline{X}).
\end{align}
This is a special case of~\eqref{eq:Cj!log}, which finishes the proof.
\endproof

\begin{remark}
\label{rem:CC}
In the case of \'etale sheaves, the analogue of the CSM-class in every dimension has been constructed in \cite[Def. 6.7]{Sai} using the theory of characteristic cycles (see also \cite[Note 87$_1$]{ReS}). One may expect that a similar construction can be performed in motivic homotopy theory, which extends our reconstruction of the $0$-dimensional part of Aluffi's pro-CSM class to higher dimensions. However in positive characteristic it is unreasonble to expect push-forward formulas like~\ref{num:Cpf} unless for the $0$-dimensional part, see \cite[Example 6.10]{Sai}.
\end{remark}

\subsection{}
\label{num:MW}
By Theorem~\ref{th:A=C} and Example~\ref{ex:MW}, for $\mathbb{E}=\mathbf{H}_{\mathrm{MW}}\mathbb{Z}$, our compactly supported pro-characteristic class
\begin{align}
lC^{\mathrm{c}}_X(\mathbbold{1}_X)
\in
lCH^{\mathrm{MW}}_0(X/k)
\end{align}
provides a quadratic refinement of the class $lC^{SM}_0(X)\in lCH_0(X)$.

\begin{remark}
\label{rem:levquad}
M. Levine communicated to us the following remark: 
for $\mathbb{E}=\mathbf{H}_{\mathrm{MW}}\mathbb{Z}$, the direct analogue of~\eqref{eq:Cj!log} between $C_{\overline{X}}(j_!\mathbbold{1}_X,\mathbf{H}_{\mathrm{MW}}\mathbb{Z})$ and the Euler class $e(\Omega^1_{\overline{X}/k}(\operatorname{log}D)^\vee,\mathbf{H}_{\mathrm{MW}}\mathbb{Z})$, as well as Corollary~\ref{cor:CXdual}, fail to hold in general. This can already be seen at the level of the categorical Euler characteristic: for example, for the compactification $\mathbb{G}_m\to\mathbb{P}^1$, we have an isomorphism of sheaves $\Omega^1_{\mathbb{P}^1}(\operatorname{log}D)\simeq \mathcal{O}_{\mathbb{P}^1}$, while the homological Euler characteristic and the compactly supported Euler characteristic of $\mathbb{G}_m$ can be computed respectively: $\chi(\mathbb{G}_m)=1-\langle-1\rangle$, $\chi^{\mathrm{c}}(\mathbb{G}_m)=\langle-1\rangle-1$ (this can be seen, for example, from the fact that the motive of $\mathbb{G}_m$ is isomorphic to $\mathbbold{1}\oplus\mathbbold{1}(1)[1]$, see \cite[5.1.12]{JY}).
\end{remark}

\subsection{}
Recall that we have defined the homological characteristic class $C^{\operatorname{Hom}}_U(K,\mathbf{H}\mathbb{Z})$ in Definition~\ref{def:CHom}. By Corollary~\ref{cor:CHomlift} and Theorem~\ref{th:A=C}, we can relate them to Aluffi's zeroth pro-Chern-Schwartz-MacPherson class $lC^{SM}_0(X)$:
\begin{corollary}
\label{cor:HomCSM}
Assume that the conditions in Theorem~\ref{thm:Alu43} are satisfied.
Let $X$ be a scheme. Then the map $H_0^{\operatorname{S}}(X/k)\to lCH_0(X)$ in Example~\ref{ex:EM} sends the class $C^{\operatorname{Hom}}_X(\mathbbold{1}_X,\mathbf{H}\mathbb{Z})$ to the class $lC^{SM}_0(X)$. 

In other words, the class $C_{U}^{\operatorname{Hom}}(K,\mathbf{H}\mathbb{Z})$ is a lift of $lC^{SM}_0(X)$ to the Suslin homology $H_0^{\operatorname{S}}(X/k)$.
\end{corollary}

\begin{remark}
\label{rm:higherdim}
By Corollary~\ref{cor:CHomlift}, the class $C_{U}^{\operatorname{Hom}}(K)$ is a lift of $lC^{\mathrm{c}}_{U}(K)$ to $\mathbb{E}^{}(U/k)$.
One may wonder if the total (pro)-CSM class lifts to motivic homology of higher dimension. This is false in general, as the following example shows (communicated to us by X. Zhang): the fundamental class of $\mathbb{A}^1$ in $CH_1(\mathbb{A}^1)$, which is a component of its CSM class, is a non-zero class, while by $\mathbb{A}^1$-invariance $\mathbf{H}\mathbb{Z}_{2n,n}(\mathbb{A}^1/k)$ vanishes for every $n\neq0$, so no such lift can exist.
\end{remark}

\subsection{}
\label{num:consf}
In the rest of this section, we assume that the base field $k$ has characteristic $0$. For a scheme $X$, the group of ($\mathbb{Z}$-valued) \textbf{constructible functions} $\operatorname{Cons}(X)$ is the subgroup of functions $X\to\mathbb{Z}$ spanned by linear combinations of functions of the form $1_Z$, where $1_Z$ is the characteristic function of a closed subscheme $Z$ of $X$. A constructible function is \textbf{locally constant} if its value is constant on each connected component.

\subsection{}
By \cite{Rio}, if $F$ is any field of characteristic $0$, then any constructible object $K\in\mathbf{DM}_{c}(F)$ is dualizable. Its categorical Euler characteristic $\chi(K)$ is then an integer.

Let $X$ be a scheme and let $K\in\mathbf{DM}_{c}(X)$ be a constructible object. The \textbf{rank function} of $K$ is the function $r(K)$ on $X$ sending a point $x:\operatorname{Spec}(k(x))\to X$ to the integer $\chi(x^*K)$. If let $K\in\mathbf{SH}_{c}(X)$ is a constructible object, we define its rank function as that of its image via the canonical functor $\mathbf{SH}_{c}(X)\to\mathbf{DM}_{c}(X)$.

We now focus on constructible objects in $\mathbf{DM}_{c}(X)$ instead of $\mathbf{SH}_{c}(X)$.

\begin{lemma}
\label{lm:duallc}
If $K\in\mathbf{DM}_{c}(X)$ is a dualizable object, then $r(K)$ is a locally constant function on $X$.
\end{lemma}
\proof
We know that the categorical Euler characteristic $\chi(K)$ is defined in the group $\operatorname{End}_{\mathbf{DM}_{c}(X)}(\mathbbold{1}_X)\simeq\mathbb{Z}^{\pi_0(X)}$, which can be identified as a locally constant function on $X$. For any point $x:\operatorname{Spec}(k(x))\to X$, the functor $x^*:\mathbf{DM}_{c}(X)\to\mathbf{DM}_{c}(k(x))$ is symmetric monoidal, and therefore $\chi(x^*K)=x^*\chi(K)$. It follows that the value of the rank function $r(K)$ is constant on each connected component of $X$, and the result follows.
\endproof

\begin{corollary}
\label{cor:rkcons}
For any constructible object $K\in\mathbf{DM}_{c}(X)$, the rank function $r(K)$ is a constructible function.
\end{corollary}
\proof
Follows from Lemma~\ref{lm:duallc} and the fact that for any constructible object $K\in\mathbf{DM}_{c}(X)$, one can find a stratification of $X$ such that the restriction of $K$ on each stratum is dualizable.
\endproof

\subsection{}
By additivity of traces (\cite{May}) and Corollary~\ref{cor:rkcons}, the rank function induces a morphism of abelian groups
\begin{align}
\label{eq:rkK0}
r:K_0(\mathbf{DM}_{c}(X))\to \operatorname{Cons}(X).
\end{align}
If $j:U\to X$ is an immersion, then $r[j_!\mathbbold{1}_U]=1_U$ is the characteristic function of $U$.

\begin{proposition}
\label{prop:alu61}
The following diagram commutes:
\begin{align}
\begin{split}
  \xymatrix@=10pt{
    K_0(\mathbf{DM}_{c}(X)) \ar[r]^-{\eqref{eq:rkK0}} \ar[rd]_-{lC^{\mathrm{c}}_X} & \operatorname{Cons}(X) \ar[d]^-{lC^{SM}_0} \\
    & lCH_0(X).
  }
\end{split}
\end{align}
\end{proposition}

\proof
By Theorem~\ref{th:A=C} and the theory of weight structures, we only need to prove the following statement. Let $j:X\to\overline{X}$ be an open immersion with $\overline{X}$ smooth and projective of dimension $m$, such that the complement $D$ of $j$ is a simple normal-crossing divisor. Let $p:W\to X$ be a smooth projective morphism with $W$ of dimension $d$, and let $\Gamma\in CH_d(W\times_XW)\simeq End(j_!p_*\mathbbold{1}_W)$. Then we have
\begin{align}
\label{Alugen}
C_{\overline{X}}(j_!p_*\mathbbold{1}_W,\Gamma)
=
(\Gamma\cdot\Delta_{W/X})\cdot C_{\overline{X}}(j_!\mathbbold{1}_X).
\end{align}
We may choose a Cartesian square
\begin{align}
\begin{split}
  \xymatrix@=10pt{
    W \ar[r]^-{k} \ar[d]_-{p} & \overline{W} \ar[d]^-{q} \\
    X \ar[r]^-{j} & \overline{X}.
  }
\end{split}
\end{align}
where $\overline{W}$ is smooth, $q$ is proper, and the complement $E$ of $k$ is the support of a simple normal-crossing divisor. 

Let $\overline{\Gamma}\in CH_d(\overline{W}\times_{\overline{X}}\overline{W})\simeq End(q_*\mathbbold{1}_{\overline{W}})$ be the closure of $\Gamma$. 
The pairing 
\begin{align}
\begin{split}
CH_d(\overline{W}\times_{\overline{X}}\overline{W})\times CH_0(\overline{W})&\to CH_0(\overline{W})\\
(\alpha,\beta)&\mapsto \pi_{2,*}(\alpha\cdot\pi_1^*\beta)
\end{split}
\end{align}
gives rise to a cycle $\overline{\Gamma}\cdot C_{\overline{X}}(k_!\mathbbold{1}_W)\in CH_0(\overline{W})$.
Then it suffices to show the following equalities:
\begin{align}
\label{eq:CjG1}
C_{\overline{X}}(j_!p_*\mathbbold{1}_W,\Gamma)
=
q_*(\overline{\Gamma}\cdot C_{\overline{X}}(k_!\mathbbold{1}_W))
\end{align}
and
\begin{align}
\label{eq:CjG2}
q_*(\overline{\Gamma}\cdot C_{\overline{X}}(k_!\mathbbold{1}_W))
=
(\Gamma\cdot\Delta_{W/X})\cdot C_{\overline{X}}(j_!\mathbbold{1}_X).
\end{align}
For formula~\eqref{eq:CjG1}, note that the cycle $q_*(\overline{\Gamma}\cdot C_{\overline{X}}(k_!\mathbbold{1}_W))$ agrees with the composition
\begin{align}
\mathbbold{1}_{\overline{X}}
\xrightarrow{}
q_*\mathbbold{1}_{\overline{W}}
\xrightarrow{\overline{\Gamma}}
q_*\mathbbold{1}_{\overline{W}}
\xrightarrow{C_{\overline{X}}(k_!\mathbbold{1}_W)}
q_*\mathcal{K}_{\overline{W}}
\xrightarrow{}
\mathcal{K}_{\overline{X}}
\end{align}
and the result follows from the following commutative diagram
\begin{align}
\begin{split}
  \xymatrix@=10pt{
    \mathbbold{1}_{\overline{X}} \ar[r]^-{} \ar[d]_-{} & \underline{Hom}(q_*k_!\mathbbold{1}_{W},q_*k_!\mathbbold{1}_{W}) \ar[r]^-{} \ar[d]^-{} & q_*k_*\mathcal{K}_{W}\otimes q_*k_!\mathbbold{1}_{W} \ar@/^2pc/[dd] \\
    q_*k_!\mathbbold{1}_{W} \ar[r]^-{\Gamma} \ar[d]_-{} & q_*k_!\mathbbold{1}_{W} \ar[d]_-{}  \ar[r]^-{C_{W}(\mathbbold{1}_W)} & q_*k_!\mathcal{K}_{W} \ar[d]_-{} \ar[u]_-{}\\
    q_*\mathbbold{1}_{\overline{W}} \ar[r]^-{\overline{\Gamma}} & q_*\mathbbold{1}_{\overline{W}} \ar[r]^-{C_{\overline{W}}(k_!\mathbbold{1}_W)} & q_*\mathcal{K}_{\overline{W}}
  }
\end{split}
\end{align}
where
\begin{enumerate}
\item
The upper left vertical map is the composition $\mathbbold{1}_{\overline{X}}\to q_*\underline{Hom}(k_!\mathbbold{1}_{W},k_!\mathbbold{1}_{W})\simeq q_*k_!\mathbbold{1}_{W}$.
\item
The commutativity of the upper right square follows from the proof of \cite[Prop. 5.1.15]{JY}.
\end{enumerate}
Formula~\eqref{eq:CjG1} is a generalization of \cite[Claim 6.1]{Alu}, and can be proved using a similar method: by Corollary~\ref{cor:Cj!}, it suffices to show that
\begin{align}
q_*(\overline{\Gamma}\cdot c_d(\Omega^1_{\overline{W}/k}(\operatorname{log}E))\cap[\overline{W}])
=
(\Gamma\cdot\Delta_{W/X})\cdot c_d(\Omega^1_{\overline{X}/k}(\operatorname{log}D))\cap[\overline{X}].
\end{align}
As in \cite[\S6]{Alu}, we apply MacPherson's graph construction (\cite[\S18]{Ful}) to the differential map 
\begin{align}
dq:q^*\Omega^1_{\overline{X}/k}(\operatorname{log}D)\to\Omega^1_{\overline{W}/k}(\operatorname{log}E)
\end{align}
which induces a rational map
\begin{align}
\xymatrix{
\gamma: \overline{W}\times\mathbb{P}^1 \ar@{.>}[r] & G=\operatorname{Gr}_m(q^*\Omega^1_{\overline{X}/k}(\operatorname{log}D)\oplus\Omega^1_{\overline{W}/k}(\operatorname{log}E))
}
\end{align}
sending a point $(x,(\lambda,1))$ with $\lambda\ne0$ to the graph of $\frac{dg}{\lambda}$ at $x$. The indeterminacies of $\lambda$ are contained in $E\times0\subset\overline{W}\times0$, and blowing up the ideal of indeterminacies gives a morphism
\begin{align}
\widetilde{\gamma}:\widetilde{\overline{W}\times\mathbb{P}^1}\to G.
\end{align}
Let $V=(\overline{W}\times_{\overline{X}}\overline{W})\times_{\overline{W}}(\widetilde{\overline{W}\times\mathbb{P}^1})$. We have a commutative diagram
\begin{align}
\begin{split}
  \xymatrix@=10pt{
    V \ar[r]^-{s} \ar[dd]_-{r} & \widetilde{\overline{W}\times\mathbb{P}^1} \ar[rd]^-{\widetilde{\gamma}} \ar[d]_-{\pi} &\\
    & \overline{W}\times\mathbb{P}^1 \ar@{.>}[r]^-{\gamma\ \ \ \ } \ar[d]_-{\rho} & G \ar[ld]^-{}\\
    \overline{W}\times_{\overline{X}}\overline{W} \ar[r]^-{\pi_1} \ar[d]_-{\pi_2}& \overline{W} \ar[d]^-{q} & \\
    \overline{W} \ar[r]^-{q} & \overline{X} &
  }
\end{split}
\end{align}
Then $[\pi^{-1}(X\times\infty)]=[\pi^{-1}(X\times0)]$ as rational equivalence classes of divisors, and the divisor $\pi^{-1}(X\times0)$ consists of the proper transform $\widetilde{X}$ of $X\times0$, as well as the components of the exceptional divisor. Let $\mathcal{Q}$ be the universal quotient bundle over $G$. Then the arguments of \cite[Lemma 6.2]{Alu} show that
\begin{align}
q_*\pi_{2,*}\left(\overline{\Gamma}\cdot r_*s^*(c_d(\widetilde{\gamma}^*\mathcal{Q})\cap[\pi^{-1}(X\times\infty)])\right)
=
q_*(\overline{\Gamma}\cdot c_d(\Omega^1_{\overline{W}/k}(\operatorname{log}E))\cap[\overline{W}])
\end{align}
\begin{align}
q_*\pi_{2,*}\left(\overline{\Gamma}\cdot r_*s^*(c_d(\widetilde{\gamma}^*\mathcal{Q})\cap[\widetilde{X}])\right)
=
(\Gamma\cdot\Delta_{W/X})\cdot c_d(\Omega^1_{\overline{X}/k}(\operatorname{log}D))\cap[\overline{X}].
\end{align}
Using the projection formula, the result then follows from \cite[Claim 6.4]{Alu}, which states that the components of the exceptional divisor have no contribution to the sum.
\endproof

\subsection{}
We now study \emph{push-forward maps} on constructible functions. Let $f:Y\to X$ be a morphism of schemes. If $Z$ is a closed subscheme of $Y$, using the map~\eqref{eq:deghom}, we define a function $f_*1_Z$ on $X$ which sends a point $x$ to the integer
\begin{align}
f_*1_Z(x)
=
\int_{Z\cap f^{-1}(x)/k(x)}
C^{\operatorname{Hom}}_{Z\cap f^{-1}(x)}(\mathbbold{1}_{Z\cap f^{-1}(x)},\mathbf{H}\mathbb{Z})
\in 
H_0^{\operatorname{S}}(k(x)/k(x))\simeq\mathbb{Z}.
\end{align}
Here the scheme $Z\cap f^{-1}(x)$ is the scheme-theoretic intersection which fits into a Cartesian square
\begin{align}
\begin{split}
  \xymatrix@=10pt{
    Z\cap f^{-1}(x) \ar[r]^-{} \ar[d]_-{} & f^{-1}(x) \ar[r]^-{} \ar[d]^-{} & \operatorname{Spec}(k(x)) \ar[d]^-{x} \\
    Z \ar[r]^-{} & Y \ar[r]^-{} & X.
  }
\end{split}
\end{align}
By Corollary~\ref{cor:HomCSM}, our definition agrees with Aluffi's in \cite[3.4]{Alu}. By \cite[Lemma 5.8]{Alu}, the function $f_*1_Z$ is a constructible function on $X$. By linearity we obtain a push-forward map on constructible functions
\begin{align}
\label{eq:pfcons}
f_*:\operatorname{Cons}(Y)\to\operatorname{Cons}(X).
\end{align}
By \cite[Thm. 5.1]{Alu}, the push-forward map is compatible with compositions: if $g:X\to W$ is another morphism of schemes, then $g_*\circ f_*=(g\circ f)_*$.

\subsection{}
The following lemma shows the compatibility between the rank function and push-forward maps:
\begin{lemma}
\label{lm:rkpf}
Let $f:Y\to X$ be a morphism of schemes. Then the following diagram commputes:
\begin{align}
\begin{split}
  \xymatrix@=10pt{
    K_0(\mathbf{DM}_{c}(Y)) \ar[r]^-{f_!} \ar[d]_-{\eqref{eq:rkK0}}^-{r}  & K_0(\mathbf{DM}_{c}(X)) \ar[d]^-{\eqref{eq:rkK0}}_-{r} \\
    \operatorname{Cons}(Y) \ar[r]_-{\eqref{eq:pfcons}}^-{f_*} & \operatorname{Cons}(X).
  }
\end{split}
\end{align}
\end{lemma}
\proof
Let $K\in\mathbf{DM}_{c}(Y)$ be a constructible object. We need to show that $f_*(r(K))=r(f_!K)$. As explained in the proof of Corollary~\ref{cor:rkcons}, we may assume that $K=j_!L$ where $j:U\to Y$ is an immersion and $L\in\mathbf{DM}_{c}(Y)$ is a dualizable object. By shringking $U$, we may assume that $U$ is connected, and the composition $U\to Y\to X$ is quasi-projective. Therefore replacing $Y$ by $U$ and $K$ by $L$, it suffices to show that $f_*(r(K))=r(f_!K)$ under the additional assumption that $K$ is dualizable object, $Y$ is connected, and $f$ is either an immersion or a smooth projective morphism.

If $f$ is an immersion, then $f_*:\operatorname{Cons}(Y)\to\operatorname{Cons}(X)$ is simply the inclusion of constructible functions, and $r(f_!K)$ is simply $r(K)$ viewed as a constructible function on $X$, and the result follows.

If $f$ is smooth projective, we may also assume that $X$ is connected. The object $f_!K$ is also dualizable by Atiyah duality, and we have $\chi(f_!K)=\chi(K)\cdot\chi(f_!\mathbbold{1})$. Therefore we may assume that $K=\mathbbold{1}_Y$ is the unit object. Let $x:\operatorname{Spec}(k(x))\to X$ be a point of $X$. We need to show that
\begin{align}
(f_*1_Y)(x)=r(f_!\mathbbold{1}_Y)(x).
\end{align}
We have the Cartesian square
\begin{align}
\begin{split}
  \xymatrix@=10pt{
   Y_x \ar[r]^-{f_x} \ar[d]_-{} & \operatorname{Spec}(k(x)) \ar[d]^-{x} \\
   Y \ar[r]^-{f} & X.
  }
\end{split}
\end{align}
By Corollary~\ref{cor:CHomfun} and proper base change we have
\begin{align}
(f_*1_Y)(x)
=
\int_{Y_x/k(x)}
C^{\operatorname{Hom}}_{Y_x}(\mathbbold{1}_{Y_x})
=
C^{\operatorname{Hom}}_{k(x)}(f_{x!}\mathbbold{1}_{Y_x})
=
\chi(f_{x!}\mathbbold{1}_{Y_x}),
\end{align}
\begin{align}
r(f_!\mathbbold{1}_Y)(x)
=
\chi(x^*f_!\mathbbold{1}_Y)
=
\chi(f_{x!}\mathbbold{1}_{Y_x})
\end{align}
which finishes the proof.
\endproof

\begin{remark}
The rank function can also be defined in positive characteristic, by letting the value of $r(K)$ on a point $x$ to be the Euler characteristic of the \'etale realization of the stalk (see \cite{Ill}). Lemma~\ref{lm:duallc} follows from the fact that the Euler characteristic of a locally constant sheaf is constant along specializations, and Corollary~\ref{cor:rkcons} also follows. However as pointed out in \cite[5.2, 6.10]{Alu}, one cannot define functorial push-forwards of constructible functions as in~\eqref{eq:pfcons} for the characteristic $0$ case: for example, in the case of the Artin-Schreier cover, subschemes may have unexpected multiplicities. 
\end{remark}

\subsection{}
We are now ready to state our theorem in the language of constructible functions. Let $k$ be a field of characteristic $0$. Let $\operatorname{Sch}_k$ be the category of separated $k$-schemes of finite type, together with morphisms of $k$-schemes. We view $\operatorname{Cons}(-)$, $H_0^{\operatorname{S}}(-/k)$ and $lCH_0(-)$ as (covariant) functors $\operatorname{Sch}_k\to Ab$ with functoriality provided by~\eqref{eq:hompf}, \eqref{eq:pfcons} and~\eqref{eq:funclimBM} respectively. In \cite[Thm. 5.2]{Alu}, Aluffi defined a natural transformation of functors 
\begin{align}
lC^{SM}_0(-):\operatorname{Cons}(-)\to lCH_0(-)
\end{align}
sending $1_Z$ to the class $lC^{SM}_0(Z)$ in Theorem~\ref{thm:Alu43}. Our result below gives a lift of this transformation to the zeroth Suslin homology $H_0^{\operatorname{S}}(-/k)$:
\begin{theorem}
\label{thm:Constr}
Let $k$ be a field of characteristic $0$.
There is a unique natural transformation of functors 
\begin{align}
C^{\operatorname{Hom}}:\operatorname{Cons}(-)\to H_0^{\operatorname{S}}(-/k)
\end{align}
such that
\begin{enumerate}
\item
(Additivity for stratifications)
If $X$ has a finite stratification into locally closed subschemes $X=\cup_{i=1}^nU_i$, with each $U_i$ smooth and irreducible, then we have
\begin{align}
\label{eq:addcsm}
C^{\operatorname{Hom}}_X(X)=\sum_{i=1}^nC^{\operatorname{Hom}}_X(U_i)\in H_0^{\operatorname{S}}(X/k).
\end{align}
\item
If $X=\overline{X}-D\xrightarrow{j}\overline{X}$ with $\overline{X}$ smooth and proper of dimension $n$ 
and $D$ a simple normal-crossing divisor in $\overline{X}$, then the class $C^{\operatorname{Hom}}_{\overline{X}}(X)$ is given, up to a sign, by the relative top Chern class of the sheaf of differential forms with logarithmic poles along $D$ (\ref{num:rellog}):
\begin{align}
\label{eq:logpole}
C^{\operatorname{Hom}}_{\overline{X}}(X)=(-1)^nc_{\overline{X},D}(\Omega^1_{\overline{X}/k}(\operatorname{log}D))\in H_0^{\operatorname{S}}(\overline{X}/k)\simeq CH_0^{}(\overline{X}).
\end{align}
\end{enumerate}
Furthermore, the image of $C^{\operatorname{Hom}}$ in $lCH_0(-)$ agrees with the zero-dimensional pro-CSM transformation $lC^{SM}_0(-)$.
\end{theorem}

\proof
For a closed immersion $i:Z\to X$, we define $C^{\operatorname{Hom}}_X(1_Z)=i_*C^{\operatorname{Hom}}_Z(\mathbbold{1}_Z,\mathbf{H}\mathbb{Z})\in H_0^{\operatorname{S}}(X/k)$. This allows to define the transformation $C^{\operatorname{Hom}}$ objectwise. 
The additivity property follows from Lemma~\ref{lm:addhom}, Corollary~\ref{cor:CHomfun} and the localization triangle~\eqref{eq:locseq}. 
By Corollary~\ref{cor:HomCSM}, $C^{\operatorname{Hom}}$ is a lift of $lC^{SM}_0(-)$. 

It remains to see that $C^{\operatorname{Hom}}$ commutes with push-forward maps on both sides. 
Using~\ref{num:rellog}, an argument similar to Proposition~\ref{prop:alu61} gives a commutative diagram
\begin{align}
\begin{split}
  \xymatrix@=10pt{
    K_0(\mathbf{DM}_{c}(X)) \ar[r]^-{\eqref{eq:rkK0}} \ar[rd]_-{C^{\operatorname{Hom}}_X} & \operatorname{Cons}(X) \ar[d]^-{C^{\operatorname{Hom}}} \\
    & H_0^{\operatorname{S}}(-/k).
  }
\end{split}
\end{align}
By Lemma~\ref{lm:rkpf}, since the rank function map~\eqref{eq:rkK0} is surjective, we can lift the statement to $K_0(\mathbf{DM}_{c}(-))$, and the result follows from Corollary~\ref{cor:CHomfun}.
\endproof

\begin{remark}
In positive characteristic the natural transformation $K_0(\mathbf{SH}_c(-))\to H_0^{\operatorname{S}}(-/k)$ in~\eqref{eq:CHom} no longer factors through $\operatorname{Cons}(-)$, as shown by the famous Grothendieck–Ogg–Shafarevich formula (\cite[X Thm. 7.1]{SGA5}, \cite[Thm. 4.2.9]{KS}), where the defect is measured by the \emph{Swan conductor}. In what follows we study characteristic classes related to the Swan conductor.
\end{remark}

\section{Boundary homology theories}

\subsection{}
In the style of \cite[Def. 3.1.1]{KSa}, we give the following definition:
\begin{definition}
\label{def:bBM}
Let $f:X\to S$ be a morphism of schemes and let $\mathbb{E}\in \mathbf{SH}(S)$ be a motivic spectrum. We define the \textbf{boundary Borel-Moore $\mathbb{E}$-homology}
\begin{align}
b\mathbb{E}^{\mathrm{BM}}_n(X/S)
=
\operatorname{lim}_{\overline{X}\in\operatorname{Cpt}(X/S)} \mathbb{E}^{\mathrm{BM}}_n(\overline{X}-X/S).
\end{align}
\end{definition}

\subsection{}
\label{num:bBM}
Similar to \ref{num:limBM}, to determine a class $\alpha$ in $b\mathbb{E}^{\mathrm{BM}}_n(X/S)$ amounts to determine, for every compactification $\overline{X}$ of $f$, a class $\alpha_{X/\overline{X}}$ in $\mathbb{E}^{\mathrm{BM}}_n(\overline{X}-X/S)$, compatible with proper push-forwards.

If $S$ is the spectrum of a perfect field $k$ and $\operatorname{Cpt}^{\operatorname{Sm}}(X/k)$ is cofinal in $\operatorname{Cpt}(X/k)$ as in~\ref{num:smcomp}, then to determine a class $\alpha$ in $b\mathbb{E}^{\mathrm{BM}}(X/k)$ amounts to determine, for every smooth compactification $\overline{X}$ of $X$, a class $\alpha_{X/\overline{X}}$ in $\mathbb{E}^{\mathrm{BM}}(\overline{X}-X/k)$, compatible with proper push-forwards.

\subsection{}
When $X$ is proper over $S$, we have $b\mathbb{E}^{\mathrm{BM}}_n(X/S)=0$.

\subsection{}
For any proper morphism $f:Y\to X$, there is a map 
\begin{align}
\label{eq:funcb}
f_*:b\mathbb{E}^{\mathrm{BM}}_n(Y/S)\to b\mathbb{E}^{\mathrm{BM}}_n(X/S).
\end{align}
Indeed, as in~\eqref{eq:comprel} we may choose a morphism $\overline{f}:\overline{Y}\to\overline{X}$ between compactifications, and since $f$ is proper, $\overline{f}$ induces a proper morphism $\overline{Y}-Y\to\overline{X}-X$. Then the map~\eqref{eq:funcb} is induced by the proper functoriality of Borel-Moore homology~\eqref{eq:BMprop}.

\subsection{}
The proper functoriality~\eqref{eq:BMprop} also induces a canonical map
\begin{align}
\label{eq:btol}
b\mathbb{E}^{\mathrm{BM}}_n(X/S)
\to
l\mathbb{E}^{\mathrm{BM}}_n(X/S)
\end{align}
sending a class on the boundary to the total space.
In particular, composing~\eqref{eq:btol} with the map~\eqref{eq:limBMBM} we obtain a \textbf{degree map}
\begin{align}
\int:
b\mathbb{E}^{\mathrm{BM}}_0(X/S)
\to
\mathbb{E}^{\mathrm{BM}}_0(S/S).
\end{align}

\begin{example}
\label{ex:bCH0}
Let $S=\operatorname{Spec}(k)$ be the spectrum of a field of characteristic $p$, and $\mathbb{E}=\mathbf{H}\mathbb{Z}$ the motivic Eilenberg-Mac Lane spectrum (\cite{Spi}). Then the boundary Borel-Moore motivic homology $b\mathbf{H}\mathbb{Z}^{\mathrm{BM}}_0(X/k)$ agrees with the \textbf{boundary Chow group of $0$-cycles} on $X$ defined in \cite[Def. 3.1.1]{KSa}:
\begin{align}
b\mathbb{E}^{\mathrm{BM}}_0(X/k)
\simeq 
bCH_0(X)
=
\operatorname{lim}_{\overline{X}\in\operatorname{Cpt}(X/k)}
CH^{}_{0}(\overline{X}-X).
\end{align}
\end{example}

\subsection{}
We now introduce a homological variant of Definition~\ref{def:bBM}.
Let $X\to S$ be a morphism of schemes. Let $p:\overline{X}\to\overline{X}'$ be a morphism of compactifications of $X$ over $S$, with a commutative diagram
\begin{align}
\begin{split}
  \xymatrix@=10pt{
   & \overline{X}  \ar[d]^-{g} & \overline{X}-X \ar[d]^-{r} \ar[l]_-{i} \\
  X \ar[r]^-{j'} \ar[ru]^-{j} & \overline{X}' & \overline{X}'-X \ar[l]_-{i'}
  }
\end{split}
\end{align}
where one may choose a scheme structure on the closed subset $\overline{X}-X$ of $\overline{X}$ such that the square on the right is Cartesian. Let $\pi:\overline{X}'\to S$ be the structure morphism, and let $\mathbb{F}\in \mathbf{SH}(X)$ be a motivic spectrum. Then we have
\begin{align}
\label{eq:bHomind}
\pi_*i'_*i'^!j'_!\mathbb{F}
\simeq
\pi_*i'_*i'^!g_!j_!\mathbb{F}
\simeq
\pi_*i'_*r_!i^!j_!\mathbb{F}
\simeq
\pi_*g_*i_!i^!j_!\mathbb{F}.
\end{align}
\begin{definition}
\label{def:bHom}
Let $f:X\to S$ be a morphism of schemes and let $\mathbb{E}\in \mathbf{SH}(S)$ be a motivic spectrum. Let $X\xrightarrow{j}\overline{X}\xrightarrow{p}S\in \operatorname{Cpt}(X/S)$ be a compactification with complement $i:\overline{X}-X\to\overline{X}$. We define the \textbf{boundary $\mathbb{E}$-homology}
\begin{align}
b\mathbb{E}^{}_n(X/S)
=
[\mathbbold{1}_S[n],p_*i_*i^!j_!j^!p^!\mathbb{E}]_{\mathbf{SH}(S)}
=
[i_*\mathbbold{1}_{\overline{X}-X}[n],j_!j^!p^!\mathbb{E}]_{\mathbf{SH}(\overline{X})}
\end{align}
By~\eqref{eq:bHomind}, the definition is independent of the choice of compactification $\overline{X}$ up to canonical isomorphism.
\end{definition}

\subsection{}
When $X$ is proper over $S$, we have $b\mathbb{E}^{}_n(X/S)=0$.

\subsection{}
The adjunction $i_*i^!\to1$ induces a map
\begin{align}
\label{eq:bhh}
b\mathbb{E}^{}_n(X/S)
\to
\mathbb{E}^{}_n(X/S)
\end{align}
which by the localization triangle~\eqref{eq:locseq1} fits into a long exact sequence 
\begin{align}
\label{eq:lesbBM}
\cdots
\to
\mathbb{E}^{\mathrm{BM}}_{n+1}(X/S)
\to
b\mathbb{E}^{}_n(X/S)
\to
\mathbb{E}^{}_n(X/S)
\to
\mathbb{E}^{\mathrm{BM}}_n(X/S)
\to\cdots.
\end{align}
In other words, boundary homology can be viewed as the difference between homology and Borel-Moore homology.

\subsection{}
The adjunction $j_!j^!\to1$ induces a map
\begin{align}
\label{eq:bHomBM}
b\mathbb{E}^{}_n(X/S)
\to
\mathbb{E}^{\mathrm{BM}}_n(\overline{X}-X/S)
\end{align}
which by the localization triangle~\eqref{eq:locseq} fits into a long exact sequence 
\begin{align}
\label{eq:blessup}
\cdots
\to
\mathbb{E}^{\mathrm{BM}}_{n+1,\overline{X}-X}(\overline{X}/S)
\to
b\mathbb{E}^{}_n(X/S)
\to
\mathbb{E}^{\mathrm{BM}}_n(\overline{X}-X/S)
\to
\mathbb{E}^{\mathrm{BM}}_{n,\overline{X}-X}(\overline{X}/S)
\to\cdots.
\end{align}
Here $\mathbb{E}^{\mathrm{BM}}_{n,\overline{X}-X}(\overline{X}/S)=[\mathbbold{1}_{\overline{X}-X}[n],i^*p^!\mathbb{E}]_{\mathbf{SH}(\overline{X}-X)}$ is the \emph{Borel-Moore $\mathbb{E}$-homology with support in $\overline{X}-X$}.

\subsection{}
Taking limit of the map~\eqref{eq:bHomBM} over all the compactifications in $\operatorname{Cpt}(X/S)$, we obtain a map
\begin{align}
\label{eq:bHombBM}
b\mathbb{E}^{}_n(X/S)
\to
b\mathbb{E}^{\mathrm{BM}}_n(X/S).
\end{align}

\subsection{}
For any proper morphism $f:Y\to X$, there is a map 
\begin{align}
\label{eq:funcbh}
f_*:b\mathbb{E}^{}_n(Y/S)\to b\mathbb{E}^{}_n(X/S).
\end{align}
Indeed, as in~\eqref{eq:comprel} we may choose a morphism $\overline{f}:\overline{Y}\to\overline{X}$ between compactifications, and since $f$ is proper, $\overline{f}$ induces a proper map $h:\overline{Y}-Y\to\overline{X}-X$.  One may choose a scheme structure on the closed subset $\overline{Y}-Y$ of $\overline{Y}$ such that the immersions $j_Y:Y\to\overline{Y}$ and $i_Y:\overline{Y}-Y\to\overline{Y}$ fit into a Cartesian diagram
\begin{align}
\begin{split}
  \xymatrix@=10pt{
    Y \ar[r]^-{j_Y} \ar[d]_-{f} & \overline{Y} \ar[d]^-{\overline{f}} & \overline{Y}-Y \ar[d]^-{h} \ar[l]_-{i_Y} \\
    X \ar[r]^-{j} & \overline{X} & \overline{X}-X. \ar[l]_-{i}
  }
\end{split}
\end{align}
The map~\eqref{eq:funcbh} is induced by the map
\begin{align}
\begin{split}
\overline{f}_*i_{Y*}i_Y^!j_{Y!}j_Y^!\overline{f}^!
&=
i_*h_*i_Y^!j_{Y!}j_Y^!\overline{f}^!
\simeq
i_*i^!\overline{f}_*j_{Y!}j_Y^!\overline{f}^!\\
&=
i_*i^!j_*f_*j_Y^!\overline{f}^!
\simeq
i_*i^!j_*j^!\overline{f}_*\overline{f}^!
\to
i_*i^!j_*j^!
\end{split}
\end{align}
and one can show that the definition is independent of the choice of compactifications.

\subsection{}
The proper functoriality~\eqref{eq:funcbh} is compatible with the maps~\eqref{eq:bhh},~\eqref{eq:bHomBM} and~\eqref{eq:bHombBM}, where the details are left as exercise.

\begin{remark}
The boundary homology in Definition~\ref{def:bHom} is related to Wildshaus's theory of \textbf{boundary motives} (\cite{Wil}):
in the category $\mathbf{DM}(S)$, define
\begin{enumerate}
\item
$\partial M_S(X)=p_*i_*i^!j_!j^!p^!\mathbbold{1}_S$ the \textbf{boudary motive} of $X$ over $S$;
\item
$M_S(X)=p_*j_!j^!p^!\mathbbold{1}_S$ the \textbf{(homological) motive} of $X$ over $S$;
\item
$M^c_S(X)=p_*j_*j^!p^!\mathbbold{1}_S$ the \textbf{motive with compact support} of $X$ over $S$.
\end{enumerate}
Then these definitions are all independent of the choice of compactification $\overline{X}$, and agree with the definitions in the literature (\cite[Def. 14.1, Def. 16.13]{MVW}, \cite[Def. 2.1]{Wil}). By definition we have
\begin{align}
b\mathbf{H}\mathbb{Z}_n(X/S)
=
[\mathbbold{1}_S[n],\partial M_S(X)]_{\mathbf{DM}(S)};
\end{align}
\begin{align}
\mathbf{H}\mathbb{Z}_n(X/S)
=
[\mathbbold{1}_S[n],M_S(X)]_{\mathbf{DM}(S)};
\end{align}
\begin{align}
\mathbf{H}\mathbb{Z}^{\mathrm{BM}}_n(X/S)
=
[\mathbbold{1}_S[n],M^c_S(X)]_{\mathbf{DM}(S)}.
\end{align}
There is a canonical distinguished triangle of functors
\begin{align}
\label{eq:dtbmot}
p_*i_*i^!j_!j^!p^!
\to
p_*j_!j^!p^!
\to
p_*j_*j^!p^!
\xrightarrow{+1}.
\end{align}
which induces a canonical distinguished triangle (see also \cite[Prop. 2.2]{Wil})
\begin{align}
\begin{split}
\partial M_S(X)\to M_S(X)\to M^c_S(X)\xrightarrow{+1}.
\end{split}
\end{align}
Note that there is another distinguished triangle of the form
\begin{align}
\begin{split}
\partial M_S(X)\to M_S(\overline{X}-X)\to M_{S}(\overline{X}/X)
\end{split}
\end{align}
where $M_{S}(\overline{X}/X)=p_*i_*i^*p^!\mathbbold{1}_S$ is the quotient motive of $\overline{X}$  by $X$, or the motive of $\overline{X}$ with supports in $\overline{X}-X$ (\cite[2.3.14]{CD1}).
\end{remark}


\subsection{}
We end this section with some vanishing of homology groups:
\begin{lemma}
\label{lm:van}
Let $X$ be a smooth scheme over a field $k$ of characteristic $p$. Let $Z$ be a closed subscheme with open complement $U$. Assume that $Z$ is smooth over $k$ or $p$ is invertible. Then the following groups vanish:
\begin{align}
\mathbf{H}\mathbb{Z}^{\mathrm{BM}}_0(Z/X),
\mathbf{H}\mathbb{Z}^{\mathrm{BM}}_{-1}(Z/X),
\mathbf{H}\mathbb{Z}^{}_0(U/X),
b\mathbf{H}\mathbb{Z}^{}_0(U/X).
\end{align}

\end{lemma}
\proof
For the first two groups, one uses the \emph{niveau spectral sequence} (\cite[8.4]{DFJK})
\begin{align}
E^1_{p,q}=\oplus_{x\in Z_{(p)}}\mathbf{H}\mathbb{Z}^{\mathrm{BM}}_{p+q}(\kappa(x)/X)\Rightarrow \mathbf{H}\mathbb{Z}^{\mathrm{BM}}_{p+q}(Z/X).
\end{align}
Therefore it suffices to check that $\mathbf{H}\mathbb{Z}^{\mathrm{BM}}_{p+q}(\kappa(x)/X)=0$ when $p+q=0$ or $-1$, which follows from purity (\cite[Prop. 4.3.10]{DJK}) and \cite[Lemma 3.2]{SV}.

For the last two groups, by~\eqref{eq:homles} and~\eqref{eq:blessup} we have the exact sequences
\begin{align}
\mathbf{H}\mathbb{Z}^{}_1(Z/Z)
\to
\mathbf{H}\mathbb{Z}^{}_0(U/X)
\to
\mathbf{H}\mathbb{Z}^{}_0(X/X)
\to
\mathbf{H}\mathbb{Z}^{}_0(Z/Z)
\end{align}
and
\begin{align}
\mathbf{H}\mathbb{Z}^{}_1(Z/Z)
\to
b\mathbf{H}\mathbb{Z}^{}_0(U/X)
\to
\mathbf{H}\mathbb{Z}^{\mathrm{BM}}_0(Z/X)=0.
\end{align}
Since the map $\mathbf{H}\mathbb{Z}^{}_0(X/X)\to\mathbf{H}\mathbb{Z}^{}_0(Z/Z)$ is injective, we see that it suffices to show that $\mathbf{H}\mathbb{Z}^{}_1(Z/Z)=0$. If $Z$ is smooth over $k$, this follows from the Beilinson-Soul\'e conjecture \emph{in weight $0$} (see \cite{Blo1}). Using alterations (\cite{Tem}), the general case follows from the smooth case by a descent spectral sequence (see \cite[Thm. 3.5]{Kel}, \cite[2.4]{Abe}).
\endproof

\section{The localized characteristic class}

\subsection{}
In the rest of this paper, we switch to the positive characteristic picture.
The theory of \emph{localized characteristic class} is introduced in \cite{AS} for \'etale sheaves and generalized to the motivic setting in \cite{JY2}. In loc. cit., emphasis is put on the quadratic aspect, while in this paper we focus on the classical oriented setting for the motivic Eilenberg-Mac Lane spectrum $\mathbf{H}\mathbb{Z}$.

\subsection{}
Let $k$ be a field and let $X$ be a $k$-scheme. Let $\delta_X:X\to X\times X$ be the diagonal morphism, with open complement $\gamma_X:X\times X-X\to X\times X$. Define the functor
\begin{align}
\label{eq:deltaD}
\delta^\Delta_X=\delta^!_X(-\otimes \gamma_{X*}\mathbbold{1}):\mathbf{SH}(X\times X)\to\mathbf{SH}(X).
\end{align}
The localization triangle~\eqref{eq:locseq1} gives a distinguished triangle
\begin{align}
\label{eq:gammadt}
\delta_{X*}\delta_X^!\mathbbold{1}_{X\times X}
\to
\mathbbold{1}_{X\times X}
\to
\gamma_{X*}\mathbbold{1}_{X\times X-X}
\xrightarrow{+1}
\end{align}
and consequently the functor $\delta^\Delta_X$ fits into a distinguished triangle
\begin{align}
\label{eq:deltaDdt}
\delta_X^*(-)\otimes \delta_X^!\mathbbold{1}_{X\times X}
\to
\delta_X^!
\to
\delta^\Delta_X
\xrightarrow{+1}.
\end{align}
If $X$ is smooth over $k$, then by \emph{absolute purity} (\cite[Prop. 4.3.10]{DJK}) there is an isomorphism
\begin{align}
\label{eq:abspur}
\operatorname{Th}(T_{X/k})
\simeq
\delta_X^!\mathbbold{1}_{X\times X}
\end{align}
and the map 
\begin{align}
\label{eq:DJK431}
\delta_X^*(-)\otimes \operatorname{Th}(T_{X/k})
\simeq
\delta_X^*(-)\otimes \delta_X^!\mathbbold{1}_{X\times X}
\to
\delta_X^!
\end{align}
is the \emph{purity transformation} associated to the regular closed immerion $\delta_X$ (\cite[4.3.1]{DJK}). Note that there is an isomorphism
\begin{align}
\label{eq:deltad}
\delta^\Delta_X\delta_{X*}K
=
\delta^!_X(\delta_{X*}K\otimes \gamma_{X*}\mathbbold{1})
\simeq
\delta^!_X\delta_{X*}(K\otimes\delta_{X}^*\gamma_{X*}\mathbbold{1})
\simeq
K\otimes\delta_{X}^*\gamma_{X*}\mathbbold{1}.
\end{align}

\subsection{}
\label{num:deltaD}
We assume that $k$ is a 
field of charactersitic $p$. Let $X$ be a smooth $k$-scheme with structure morphism $f:X\to k$, and let $i:Z\to X$ be a nowhere dense closed immersion with open complement $j:U\to X$. 
Let $\mathbb{E}\in \mathbf{SH}(k)$ be a motivic spectrum with a map $\mathbbold{1}_k\to \mathbb{E}$. 
We consider the following two classes, which are well-defined after inverting $p$:
\begin{enumerate}
\item
Let $K\in\mathbf{SH}_c(X)$ be a constructible object such that $K_{|U}$ is dualizable. 
The map
\begin{align}
\label{eq:BXZ}
\begin{split}
B_X^Z(K,\mathbb{E})&:
i_*\mathbbold{1}_Z
\to
i_*i^*\delta^\Delta_X\delta_{X*}\mathbbold{1}_X
\to
i_*i^*\delta^\Delta_X\delta_{X*}\underline{Hom}(K,K)\\
&\simeq
i_*i^*\delta^\Delta_X\delta_{X*}\delta_X^!(\mathbb{D}_{X/k}(K)\boxtimes_kK)
\to
i_*i^*\delta^\Delta_X(\mathbb{D}_{X/k}(K)\boxtimes_kK)\\
&\simeq
\delta^\Delta_X(\mathbb{D}_{X/k}(K)\boxtimes_kK)
\to
\delta^\Delta_X\delta_{X*}f^!\mathbbold{1}_k
\to
\delta^\Delta_X\delta_{X*}f^!\mathbb{E}.
\end{split}
\end{align}
is viewed as a class $B_X^Z(K,\mathbb{E})\in[i_*\mathbbold{1}_Z,\delta^\Delta_X\delta_{X*}f^!\mathbb{E}]$.
\item
Let $F\in\mathbf{SH}_c(U)$ be a dualizable object. The map
\begin{align}
\label{eq:BXUZ}
\begin{split}
B_{X,U}^Z(F,\mathbb{E})&:
i_*\mathbbold{1}_Z
\to
i_*i^*\delta^\Delta_X\delta_{X*}\mathbbold{1}_X
\to
i_*i^*\delta^\Delta_X\delta_{X*}\underline{Hom}(j_!F,j_!F)\\
&\simeq
i_*i^*\delta^\Delta_X\delta_{X*}\delta_X^!(\mathbb{D}_{X/k}(j_!F)\boxtimes_kj_!F)
\to
i_*i^*\delta^\Delta_X(\mathbb{D}_{X/k}(j_!F)\boxtimes_kj_!F)\\
&\simeq
\delta^\Delta_X(\mathbb{D}_{X/k}(j_!F)\boxtimes_kj_!F)
\to
\delta^\Delta_X\delta_{X*}j_!(\mathbb{D}_{U/k}(F)\otimes_kF)
\to
\delta^\Delta_X\delta_{X*}j_!j^!f^!\mathbbold{1}_k\\
&\to
\delta^\Delta_X\delta_{X*}j_!j^!f^!\mathbb{E}.
\end{split}
\end{align}
is viewed as a class $B_{X,U}^Z(F,\mathbb{E})\in[i_*\mathbbold{1}_Z,\delta^\Delta_X\delta_{X*}j_!j^!f^!\mathbb{E}]$.
\end{enumerate}
In this section, we study the lifting problems for these two classes:
\begin{enumerate}
\item
Whether the class $B_X^Z(K,\mathbb{E})$ in~\eqref{eq:BXZ} can (uniquely) be lifted through the canonical map 
\begin{align}
\mathbb{E}^{\mathrm{BM}}_0(Z/k)
\to
[i_*\mathbbold{1}_Z,\delta^\Delta_X\delta_{X*}f^!\mathbb{E}].
\end{align}
\item
Whether the class $B_{X,U}^Z(F,\mathbb{E})$ in~\eqref{eq:BXUZ} can (uniquely) be lifted through the canonical map 
\begin{align}
b\mathbb{E}_0(U/k)
\to
[i_*\mathbbold{1}_Z,\delta^\Delta_X\delta_{X*}j_!j^!f^!\mathbb{E}].
\end{align}
\end{enumerate}

\subsection{}
If $F\in\mathbf{SH}_c(U)$ is a dualizable object, we have the following commutative diagram:
\begin{align}
\label{eq:CXUcomp}
\begin{split}
  \xymatrix@R=10pt{
   \mathbbold{1}_X \ar[r]^-{C_{X,U}(F,\mathbb{E})}_-{\eqref{eq:CXU}} \ar[d]^-{} & j_!j^!f^!\mathbb{E} \ar[d]^-{} \\
   i_*\mathbbold{1}_Z \ar[r]^-{B_{X,U}^Z(F,\mathbb{E})} \ar[rd]_-{B_X^Z(j_!F,\mathbb{E})} & \delta^\Delta_X\delta_{X*}j_!j^!f^!\mathbb{E} \ar[d]^-{} \\
    & \delta^\Delta_X\delta_{X*}f^!\mathbb{E}.
  }
\end{split}
\end{align}

\begin{lemma}
\label{lm:dualb0}
If $K$ is a dualizable object, then $B_X^Z(K,\mathbb{E})=0$.
\end{lemma}
\proof
This is because the purity transformation~\eqref{eq:DJK431} is an isomorphism when applied to dualizable objects (\cite[5.4]{FHM}), and therefore $\delta^\Delta_X(\mathbb{D}(K)\boxtimes_kK)=0$.
\endproof

\begin{lemma}
\label{lm:addB}
\begin{enumerate}
\item
If $L\to M\to N$ is a distinguished triangle in $\mathbf{SH}_c(X)$ whose restriction to $U$ are dualizable, then one has $B_X^Z(M,\mathbb{E})=B_X^Z(L,\mathbb{E})+B_X^Z(N,\mathbb{E})$.
\item
If $L\to M\to N$ is a distinguished triangle of dualizable objects in $\mathbf{SH}_c(U)$, then one has $B_{X,U}^Z(M,\mathbb{E})=B_{X,U}^Z(L,\mathbb{E})+B_X^Z(N,\mathbb{E})$.
\end{enumerate}
\end{lemma}
\proof
The proof is similar to Lemma~\ref{lm:addhom}. We prove the first claim: using the language of higher categories, there exist objects $u,v,w$ and a commutative diagram of the form
\begin{align}
\begin{split}
\xymatrix@=6pt{
    & i_*i^*\delta^\Delta_X\delta_{X*}\mathbbold{1}_X \ar^-{}[d] \ar^-{}[rd] \ar_-{}[ld] &\\
   i_*i^*\delta^\Delta_X(\mathbb{D}_{X/k}(L)\boxtimes_kL\oplus\mathbb{D}_{X/k}(N)\boxtimes_kN) \ar^-{\wr}[dd] \ar^-{}[rd]
   & i_*i^*\delta^\Delta_X\delta_{X*}v \ar^-{}[r] \ar^-{}[l] 
   & i_*i^*\delta^\Delta_X(\mathbb{D}_{X/k}(M)\boxtimes_kM) \ar^-{}[ld] \ar_-{\wr}[dd]\\
   & i_*i^*\delta^\Delta_X\delta_{X*}w \ar@{.>}^-{}[d]
   &\\
   \delta^\Delta_X(\mathbb{D}_{X/k}(L)\boxtimes_kL\oplus\mathbb{D}_{X/k}(N)\boxtimes_kN) \ar^-{}[r]  \ar^-{}[rd]
   & \delta^\Delta_X u \ar^-{}[d]
   & \delta^\Delta_X(\mathbb{D}_{X/k}(M)\boxtimes_kM) \ar^-{}[l] \ar^-{}[ld]\\
   & \delta^\Delta_X\delta_{X*}f^!\mathbbold{1}_k &
 }
 \end{split}
\end{align}
which finishes the proof. The second claim is similar and is left as an exercise.
\endproof

\subsection{}
We now consider a particular case where $K_{|U}=0$:
\begin{lemma}
\label{lm:iCZ}
Let $F\in\mathbf{SH}_c(Z)$ be a constructible object and let $K=i_*F$.
Then the following diagram is commutative:
\begin{align}
\begin{split}
  \xymatrix@R=8pt{
   i_*\mathbbold{1}_Z \ar[r]^-{C_Z(F,\mathbb{E})} \ar[rd]_-{B_X^Z(i_*F,\mathbb{E})} & f^!\mathbb{E} \ar[d]^-{} \\
    & \delta^\Delta_X\delta_{X*}f^!\mathbb{E}.
  }
\end{split}
\end{align}
In other words, the class $C_Z(F,\mathbb{E})\in\mathbb{E}^{\mathrm{BM}}_0(X/k)[1/p]\simeq[i_*\mathbbold{1}_Z,f^!\mathbb{E}]$ in~\eqref{eq:CX} is a lift of the class $B_X^Z(i_*F,\mathbb{E})$ in~\eqref{eq:BXZ}.
\end{lemma}
\proof
Using the isomorphism $\mathbb{D}_{X/k}(i_*F)\boxtimes_ki_*F\simeq (i\times i)_*(\mathbb{D}_{Z/k}(F)\boxtimes_kF)$, we have the following commutative diagram:
\begin{align*}
\begin{split}
  \xymatrix@=6pt{
     & i_*\delta^!_Z(\mathbb{D}_{Z/k}(F)\boxtimes_kF) \ar[r]^-{} \ar[d]^-{\wr} & i_*i^!f^!\mathbbold{1}_k  \ar[d]_-{} \ar[rd]_-{} &  \\ 
     i_*\mathbbold{1}_Z \ar[ru]^-{} \ar[d]^-{}  & \delta^!_X(i\times i)_*(\mathbb{D}_{Z/k}(F)\boxtimes_kF) \ar[d]_-{} \ar[r]^-{} & \delta_X^\Delta\delta_{X*}i_*i^!f^!\mathbbold{1}_k \ar[rd]_-{} & f^!\mathbbold{1}_k \ar[d]^-{} \\
     i_*i^*\delta^\Delta_X\delta_{X*}\mathbbold{1}_X \ar[r]^-{} \ar[rd]^-{} & i_*i^*\delta^\Delta_X(i\times i)_*(\mathbb{D}_{Z/k}(F)\boxtimes_kF) \ar[r]^-{\sim} \ar[d]^-{\wr} & \delta^\Delta_X(i\times i)_*(\mathbb{D}_{Z/k}(F)\boxtimes_kF)  \ar[u]^-{} \ar[d]^-{\wr}  & \delta^\Delta_X\delta_{X*}f^!\mathbbold{1}_k \\
    & i_*i^*\delta^\Delta_X(\mathbb{D}_{X/k}(i_*F)\boxtimes_ki_*F) \ar[r]^-{\sim} & \delta^\Delta_X(\mathbb{D}_{X/k}(i_*F)\boxtimes_ki_*F) \ar[ru]^-{} &    
  }
  \end{split}
\end{align*}
and the result follows.
\endproof

\subsection{}
In what follows, we focus on the case $\mathbb{E}=\mathbf{H}\mathbb{Z}$. In this case we generalize Lemma~\ref{lm:iCZ} by showing that the class $B_X^Z(K)=B_X^Z(K,\mathbf{H}\mathbb{Z})$ in~\eqref{eq:BXZ} always lifts uniquely to a $0$-cycle in the Chow group of $Z$:
\begin{corollary}
\label{cor:CXZ}
In the setting of~\ref{num:deltaD}, let $K\in\mathbf{SH}_c(X)$ be a constructible object such that $K_{|U}$ is dualizable. 
Then the class
$B_X^Z(K,\mathbf{H}\mathbb{Z})$ in~\eqref{eq:BXZ}
lifts uniquely to a class 
\begin{align}
\label{eq:CXZ}
C_X^Z(K)\in[i_*\mathbbold{1}_Z, f^!\mathbf{H}\mathbb{Z}]
=
\mathbf{H}\mathbb{Z}^{\mathrm{BM}}_0(Z/k)[1/p]\simeq CH_0(Z)[1/p]
\end{align}
called the \emph{\textbf{localized characteristic class}}.
\end{corollary}
\proof
By Lemma~\ref{lm:van} and the distinguished triangle~\eqref{eq:deltaDdt}, the following canonnical map
\begin{align}
[i_*\mathbbold{1}_Z, f^!\mathbf{H}\mathbb{Z}]
\to
[i_*\mathbbold{1}_Z,\delta^\Delta_X\delta_{X*}f^!\mathbf{H}\mathbb{Z}]
\end{align}
is an isomorphism, and the result follows.
\endproof

\subsection{}
If the base field $k$ is perfect and satisfies embedded resolution of singularities in the sense of~\ref{num:RS}, then the localized characteristic class $C_X^Z(K)$ can be defined in $CH_0(Z)$ without inverting $p$.

\subsection{}
\label{num:CXZdual0}
By Lemma~\ref{lm:dualb0}, if $K$ is a dualizable object, then $C_X^Z(K)=0$.

\subsection{}
By Lemma~\ref{lm:addB}, if $L\to M\to N$ is a distinguished triangle in $\mathbf{SH}_c(X)$ whose restriction to $U$ are dualizable, then one has 
\begin{align}
\label{eq:CXZadd}
C_X^Z(M)=C_X^Z(L)+C_X^Z(N)\in CH_0(Z).
\end{align}
In particular by~\ref{num:CXZdual0} and the localization triangle~\eqref{eq:locseq} we have 
\begin{align}
\label{eq:CXZ-}
C_X^Z(i_*\mathbbold{1}_Z)
=
C_X^Z(\mathbbold{1}_X)
-
C_X^Z(j_!\mathbbold{1}_U)
=
-C_X^Z(j_!\mathbbold{1}_U)
\in CH_0(Z).
\end{align}

\subsection{}
By Lemma~\ref{lm:iCZ}, if $K=i_*F$ where $F\in\mathbf{SH}_c(Z)$ is a constructible object, then we have 
\begin{align}
\label{eq:CXZcl}
C_X^Z(i_*F)=C_Z(F,\mathbf{H}\mathbb{Z})\in CH_0(Z).
\end{align}
In particular by~\eqref{eq:CXZ-} we have
\begin{align}
\label{eq:CXdiff}
\begin{split}
C_X(\mathbbold{1}_X,\mathbf{H}\mathbb{Z})
&=
C_X(j_!\mathbbold{1}_U,\mathbf{H}\mathbb{Z})
+
C_X(i_*\mathbbold{1}_Z,\mathbf{H}\mathbb{Z})\\
&=
C_X(j_!\mathbbold{1}_U,\mathbf{H}\mathbb{Z})
+
i_*C_Z(\mathbbold{1}_Z,\mathbf{H}\mathbb{Z})\\
&=
C_X(j_!\mathbbold{1}_U,\mathbf{H}\mathbb{Z})
+
i_*C_X^Z(i_*\mathbbold{1}_Z)\\
&=
C_X(j_!\mathbbold{1}_U,\mathbf{H}\mathbb{Z})
-
i_*C_X^Z(j_!\mathbbold{1}_U)
\in CH_0(X).
\end{split}
\end{align}

\subsection{}
We recall the notion of rank of a constructible sheaf.
Let $X$ be a connected scheme and let $\mathcal{F}\in D^b_c(X_{\textrm{\'et}},\Lambda)$. Let $U\subset X$ be an open subscheme such that $\mathcal{F}_{|U}$ has locally constant cohomology. The \textbf{rank} $\operatorname{rk}(\mathcal{F})\in\mathbb{Z}$ of $\mathcal{F}$ is defined as the Euler characteristic of $\mathcal{F}_{|U}$, that is, the alternating sum of the dimensions of its cohomology groups, an integer independent of the choice of $U$. If $X$ is not necessarily connected, the rank should be viewed as a locally constant $\mathbb{Z}$-valued function on $X$.

If $K\in\mathbf{SH}_c(X)$ is a constructible object, its \textbf{rank} $\operatorname{rk}(K)\in\mathbb{Z}$ is the rank of its \'etale realization.

\subsection{}
We now study the lifting of the class $B_{X,U}^Z(F,\mathbf{H}\mathbb{Z})$ in~\eqref{eq:BXUZ} to $b\mathbf{H}\mathbb{Z}_0(U/k)$. In the setting of~\ref{num:deltaD}, by~\eqref{eq:lesbBM}  and~\eqref{eq:deltaDdt}, we have the following commutative diagram:
\begin{align}
\label{eq:ASdiag}
\begin{split}
  \xymatrix{
      & b\mathbf{H}\mathbb{Z}_0(U/X)=0 \ar@{=}[r] \ar[d]^-{} & \mathbf{H}\mathbb{Z}_0(U/X) \ar[d]_-{} \\ 
      & [i_*\mathbbold{1}_Z,j_!j^!f^!\mathbf{H}\mathbb{Z}] \ar[d]_-{a} \ar[r]^-{c} & [\mathbbold{1}_X,j_!j^!f^!\mathbf{H}\mathbb{Z}] \ar[d]_-{b} \\
     \mathbb{Z}[1/p] \ar[d]_-{\wr} & [i_*\mathbbold{1}_Z,\delta^\Delta_X\delta_{X*}j_!j^!f^!\mathbf{H}\mathbb{Z}] \ar[r]^-{} \ar[d]_-{\partial} & [\mathbbold{1}_X,\delta^\Delta_X\delta_{X*}j_!j^!f^!\mathbf{H}\mathbb{Z}] \ar[d]_-{} \\ 
  \mathbf{H}\mathbb{Z}^{\mathrm{BM}}_{0}(U/X) \ar[r]^-{d} & b\mathbf{H}\mathbb{Z}_{-1}(U/X) \ar[r]^-{} & \mathbf{H}\mathbb{Z}_{-1}(U/X)
  }
  \end{split}
\end{align}
where both columns and the bottom row are exact. By Lemma~\ref{lm:van}, the maps $a$, $b$ and $d$ are injective.

\begin{lemma}
\label{lm:AS524}
\begin{enumerate}
\item
The class
\begin{align}
\label{eq:beta}
\beta(F)=
B_{X,U}^Z(F,\mathbf{H}\mathbb{Z})
-
\operatorname{rk}(F)\cdot B_{X,U}^Z(\mathbbold{1}_U,\mathbf{H}\mathbb{Z})
\in
[i_*\mathbbold{1}_Z,\delta^\Delta_X\delta_{X*}j_!j^!f^!\mathbf{H}\mathbb{Z}]
\end{align}
lifts uniquely to a class
\begin{align}
C_{X,U}^Z(F)\in [i_*\mathbbold{1}_Z,j_!j^!f^!\mathbf{H}\mathbb{Z}].
\end{align}
\item
The image of $C_{X,U}^Z(F)$ under the map $c$ in the diagram~\eqref{eq:ASdiag} is the class
\begin{align}
\alpha(F)=
C_{X,U}(F,\mathbf{H}\mathbb{Z})-\operatorname{rk}(F)\cdot C_{X,U}(\mathbbold{1}_U,\mathbf{H}\mathbb{Z})
\in 
[\mathbbold{1}_X,j_!j^!f^!\mathbf{H}\mathbb{Z}].
\end{align}
\end{enumerate}
\end{lemma}
\proof
By the upper part of the commutative diagram~\eqref{eq:CXUcomp}, the two classes $\alpha(F)$ and $\beta(F)$ have the same image in $[\mathbbold{1}_X,\delta^\Delta_X\delta_{X*}j_!j^!f^!\mathbf{H}\mathbb{Z}]$. Using the commutative diagram~\eqref{eq:ASdiag}, both claims are reduced to show that $\partial\beta(F)=0\in b\mathbf{H}\mathbb{Z}_{-1}(U/X)$. 
Using the commutative diagram~\eqref{eq:ASdiag} again, we know that there exists a unique class $\sigma(F)\in\mathbf{H}\mathbb{Z}^{\mathrm{BM}}_{0}(U/X)\simeq\mathbb{Z}[1/p]$ such that $d(\sigma(F))=\partial\beta(F)$. So we are reduced to show that $\sigma(F)=0$.

Abbes and Saito proved in \cite[Lemma 5.2.4]{AS} that if we replace $F$ by an \'etale sheaf $\mathcal{F}$ of $\mathbb{Z}/\ell\mathbb{Z}$-modules, where $\ell$ is a prime number invertible in $k$, then $\sigma(\mathcal{F})=0$. Since the six functors, and therefore the construction of the class $\sigma(F)$ are compatible with \'etale realizations, we know that the image of $\sigma(F)$ in $\mathbb{Z}/\ell\mathbb{Z}$ is $0$ for every prime $\ell$ different from $p$. If $\sigma(F)\ne0\in\mathbb{Z}[1/p]$, one may take a sufficiently large prime $\ell$ to obtain a contradiction, which proves the claim.
\endproof

\subsection{}
\label{num:CXUZim}
By the lower part of the commutative diagram~\eqref{eq:CXUcomp}, the image of $C_{X,U}^Z(F)$ in $\mathbf{H}\mathbb{Z}^{\mathrm{BM}}_0(Z/k)\simeq CH_0(Z/k)[1/p]$ is $C_X^Z(j_!F)-\operatorname{rk}(F)\cdot C_X^Z(j_!\mathbbold{1}_U)$.

\subsection{}
If $X$ is proper, the class $C_{X,U}^Z(F)$ is defined in $b\mathbf{H}\mathbb{Z}^{}_0(U/k)$.

\subsection{}
We now prove a formula relating the localized characteristic class in~\eqref{eq:CXZ} and the characteristic class in~\eqref{eq:CX}, which is an analog of \cite[Prop. 5.2.3]{AS} in Chow groups. 
\begin{theorem}
\label{th:GOS}
Let $X$ be a smooth $k$-scheme with structure morphism $f:X\to k$, and let $i:Z\to X$ be a nowhere dense closed immersion with open complement $j:U\to X$. Let $K\in\mathbf{SH}_c(X)$ be a constructible object such that $K_{|U}$ is dualizable.
Then the following formula holds:
\begin{align}
\label{eq:GOS}
C_X(K,\mathbf{H}\mathbb{Z})
=
\operatorname{rk}(K)\cdot C_X(\mathbbold{1}_X,\mathbf{H}\mathbb{Z})
+
i_*C_X^Z(K)
\in 
CH_0(X/k)[1/p].
\end{align}
If the base field $k$ is perfect and satisfies embedded resolution of singularities in the sense of~\ref{num:RS}, then formula~\eqref{eq:GOS} holds in $CH_0(X/k)$ without inverting $p$.
\end{theorem}
\proof
By the localization triangle~\eqref{eq:locseq}, additivity~\eqref{eq:adddt} and ~\eqref{eq:CXZadd}, and the special case where $K_{|U}=0$ in~\eqref{eq:CXZcl}, we are reduced to show that if $F\in\mathbf{SH}_c(U)$ is a dualizable object, then
\begin{align}
C_X(j_!F,\mathbf{H}\mathbb{Z})
=
\operatorname{rk}(F)\cdot C_X(\mathbbold{1}_X,\mathbf{H}\mathbb{Z})
+
i_*C_X^Z(j_!F)
\in 
CH_0(X/k)[1/p].
\end{align}
By~\eqref{eq:CXdiff}, we are reduced to show the following equality:
\begin{align}
\label{eq:CZtr}
C_X(j_!F,\mathbf{H}\mathbb{Z})-\operatorname{rk}(F)\cdot C_X(j_!\mathbbold{1}_U,\mathbf{H}\mathbb{Z})
=
i_*C_X^Z(j_!F)-\operatorname{rk}(F)\cdot i_*C_X^Z(j_!\mathbbold{1}_U).
\end{align}
By Lemma~\ref{lm:AS524} and Corollary~\ref{cor:compim}, the image of the class $C_{X,U}^Z(F)\in [i_*\mathbbold{1}_Z,j_!j^!f^!\mathbf{H}\mathbb{Z}]$ in $\mathbf{H}\mathbb{Z}^{\mathrm{BM}}_0(X/k)\simeq CH_0(X/k)[1/p]$ is the left-hand side of~\eqref{eq:CZtr}, which agrees with the right-hand side of~\eqref{eq:CZtr} by~\ref{num:CXUZim}, and the proof is finished.
\endproof

\subsection{}
If $X$ is in addition proper, taking degrees on both sides of~\eqref{eq:GOS} recovers the Grothendieck–Ogg–Shafarevich formula (\cite[X Thm. 7.1]{SGA5}, \cite[Thm. 4.2.9]{KS}).





 


\section{Functoriality of the localized characteristic class}


\subsection{}
Let $X$ be a smooth $k$-scheme and let $Z\to X$ be a nowhere dense closed subscheme, with open complement $j:U\to X$. Recall that if $K\in\mathbf{SH}_c(X)$ is a constructible object such that $K_{|U}$ is dualizable, we defined in Corollary~\ref{cor:CXZ} the \emph{localized characteristic class}
\begin{align}
C_X^Z(K)\in CH_0(Z)[1/p].
\end{align}
In this section, we establish pull-back and push-forward functorialities of the localized characteristic class $C_X^Z(K)$. We introduce the following variant:
\begin{definition}
Let $F\in\mathbf{SH}_c(U)$ be a dualizable object. We define a class
\begin{align}
\label{eq:CcUX}
\begin{split}
C^{\mathrm{c}}_{U/X}(F)
&=
C_{X}^{Z}(j_!F)-\operatorname{rk}(F)\cdot C_{X}^{Z}(j_!\mathbbold{1}_U)\\
&=
C_{X}^{Z}(j_!F)+\operatorname{rk}(F)\cdot C_{Z}(\mathbbold{1},\mathbf{H}\mathbb{Z})
\in CH_0(Z)[1/p].
\end{split}
\end{align}
If $k$ is perfect and satisfies embedded resolution of singularities in the sense of~\ref{num:RS}, then the class $C^{\mathrm{c}}_{U/X}(F)$ is defined in $CH_0(Z)$ without inverting $p$.
\end{definition}

\begin{remark}
By analogy with Definition~\ref{df:proclass}, one can also consider the class
\begin{align}
\label{eq:CUXcoh}
C^{}_{U/X}(K)
=
C_{X}^{Z}(j_*K)-\operatorname{rk}(K)\cdot C_{X}^{Z}(j_*\mathbbold{1})\in CH_0(Z)[1/p].
\end{align}
instead of the class $C^{\mathrm{c}}_{U/X}(K)$ in~\eqref{eq:CcUX}. Similar to Corollary~\ref{cor:Ccdual} and Corollary~\ref{cor:CXdual}, we have the following equality:
\begin{align}
\label{eq:CUXdual}
C^{}_{U/X}(K)
=
C^{\mathrm{c}}_{U/X}(\mathbb{D}(K))
=
C^{\mathrm{c}}_{U/X}(K)
\in 
CH_0(Z)[1/p].
\end{align}
\end{remark}

\subsection{}
\label{num:addCc}
For a distinguished triangle $L\to M\to N$ of dualizable objects in $\mathbf{SH}_c(X)$, we have
\begin{align}
\label{eq:addCc}
C^{\mathrm{c}}_{U/X}(L)=C^{\mathrm{c}}_{U/X}(M)+C^{\mathrm{c}}_{U/X}(N).
\end{align}
Indeed, the additivity follows from the additivity of the localized characteristic class (\cite[Lemma 4.12]{JY2}) and additivity of the rank (\cite{May}).

\subsection{}
By Theorem~\ref{th:GOS}, the image of $C^{\mathrm{c}}_{U/X}(F)$ in $CH_0(X)[1/p]$ is the class $C_{X}(j_!F)-\operatorname{rk}(F)\cdot C_{X}(j_!\mathbbold{1})$.

\subsection{}
\label{num:birset}
Let $f:Y\to X$ be a morphism between smooth connected $k$-schemes. Let $q:X\to k$ be the structure morphism. Let $i:Z\to X$ be a nowhere dense closed immersion with open complement $j:U\to X$. Form the Cartesian diagram
\begin{align}
\label{eq:pCart}
\begin{split}
  \xymatrix@=10pt{
    V \ar[r]^-{l} \ar[d]_-{g} & Y \ar[d]^-{f} & W \ar[d]^-{h} \ar[l]_-{k} \\
    U \ar[r]^-{j} & X & Z. \ar[l]_-{i}
  }
\end{split}
\end{align}
The following proposition is an analogue of \cite[Cor. 4.3.4]{KSa}, and is the main result of this section:
\begin{proposition} 
\label{prop:CXZpf}
In the setting of~\ref{num:birset}, assume that $f$ is proper and $g$ is \'etale (and finite). Let $F\in\mathbf{SH}_c(V)$ be a dualizable object. 
Then the following equality holds in $CH_0(Z)[1/p]$:
\begin{align}
\label{eq:locpf}
C^{\mathrm{c}}_{U/X}(g_*F)
=
h_*C^{\mathrm{c}}_{V/Y}(F)
+
\operatorname{rk}(F)\cdot C^{\mathrm{c}}_{U/X}(g_*\mathbbold{1}_V)
\end{align}
where $h_*:CH_0(W)\to CH_0(Z)$ is the proper push-forward map.
\end{proposition}

\subsection{}
The proof of Proposition~\ref{prop:CXZpf} is technical and will be given later in Section~\ref{sec:pfpbpf}. Note that by the abstract formalism of six functors, the same proof also works for \'etale sheaves.

We now exploit some of its consequences. The first is the following push-forward functoriality of the localized characteristic class:
\begin{corollary}
In the setting of~\ref{num:birset}, assume that $f$ is proper and $g$ is \'etale (and finite). Let $K\in\mathbf{SH}_c(Y)$ be a constructible object such that $K_{|X-Z}$ is dualizable. Then the following equality holds:
\begin{align}
\label{eq:locpfg}
C^{Z}_{X}(f_*K)
=
h_*C^{W}_{Y}(K)
+
\operatorname{rk}(K)\cdot C^{Z}_{X}(f_*\mathbbold{1}_Y)
\in CH_0(Z)[1/p].
\end{align}
\end{corollary}
\proof
By~\eqref{eq:CXprop}, ~\ref{num:CXZdual0}, ~\eqref{eq:CXZadd}, ~\eqref{eq:CXZcl} and~\eqref{eq:locpf} we have
\begin{align}
\begin{split}
&C^{Z}_{X}(f_*K)
=
C^{Z}_{X}(j_!j^*f_*K)
+
C^{Z}_{X}(i_*i^*f_*K)\\
=
&\operatorname{rk}(K)\cdot C^{Z}_{X}(j_!g_*\mathbbold{1}_V)
+
h_*(C^{W}_{Y}(l_!l^*K)-\operatorname{rk}(K)\cdot C^{W}_{Y}(l_!\mathbbold{1}_V))
+
h_*C^{W}_{Y}(k_*k^*K)\\
=
&h_*C^{W}_{Y}(K)
+
\operatorname{rk}(K)\cdot (C^{Z}_{X}(f_*\mathbbold{1}_Y)-C^{Z}_{X}(i_*i^*f_*\mathbbold{1}_Y)+h_*C_W(\mathbbold{1}_W))\\
=
&h_*C^{W}_{Y}(K)
+
\operatorname{rk}(K)\cdot C^{Z}_{X}(f_*\mathbbold{1}_Y)
\end{split}
\end{align}
which finishes the proof.
\endproof

\begin{remark}
In the case where $f$ is \'etale, the geometric ramification in the diagram~\eqref{eq:pCart} is much simpler, and we have the following result: 
\begin{lemma}
\label{lm:pbpfez}
In the setting of~\ref{num:birset}, 
the following hold:
\begin{enumerate}
\item Assume that $f$ is \'etale. 
Let $K\in\mathbf{SH}_c(X)$ be a constructible object such that $K_{|X-Z}$ is dualizable. Then the following equality holds in $CH_0(W)[1/p]$:
\begin{align}
\label{eq:locpbeq1}
h^*C^Z_{X}(K)
=
C^{W}_{Y}(f^*K).
\end{align}
\item Assume that $f$ is finite and \'etale. Let $K\in\mathbf{SH}_c(Y)$ be a constructible object such that $K_{|V}$ is dualizable. 
Then the following equality holds in $CH_0(Z)[1/p]$:
\begin{align}
\label{eq:locpf1}
C^Z_{X}(f_*K)
=
h_*C^{W}_{Y}(K).
\end{align}
\end{enumerate}
\end{lemma}
\proof
Formula~\eqref{eq:locpbeq1} can be proved directly, while~\eqref{eq:locpf1} is a direct consequence of~\eqref{eq:locpfg}. The details are left as exercise.
\endproof

However, if $f$ is not assumed \'etale, one cannot expect simple formulas as ~\eqref{eq:locpbeq1} and~\eqref{eq:locpf1}:
\begin{enumerate}
\item
Assume that $g$ is the identity of $U$ and both $X$ and $Y$ are smooth compactifications of $U$. 
If the equality~\eqref{eq:locpbeq1} were to hold for $K=j_!\mathbbold{1}_U$, then by Lemma~\ref{lm:Chowpf} below we would have 
\begin{align}
\int C_X^Z(j_!\mathbbold{1}_U)=\int C_Y^{Y_Z}(j_{Y!}\mathbbold{1}_U).
\end{align}
By Theorem~\ref{th:GOS}, this would imply that $\chi(\mathbbold{1}_X)=\chi(\mathbbold{1}_Y)$, which is absurd since it is well-known that the \'etale or topological Euler characteristic is \emph{not} a birational invariant. 
\item
Assume that $X$ and $Y$ are smooth and proper.
If~\eqref{eq:locpf1} were to hold for $K=\mathbbold{1}_Y$, then we would have 
\begin{align}
\int C_X^Z(f_*\mathbbold{1}_Y)=0.
\end{align}
By Theorem~\ref{th:GOS}, this would imply that $\chi(f_*\mathbbold{1}_Y)=\operatorname{rk}(f_*\mathbbold{1}_Y)\cdot\chi(\mathbbold{1}_X)$, which contradicts the Grothendieck-Ogg-Shafarevich formula (\cite[X Thm. 7.1]{SGA5}, \cite[Thm. 4.2.9]{KS}) in the case where $f$ has wild ramification at infinity.
\end{enumerate}
\begin{lemma}
\label{lm:Chowpf}
Consider a Cartesian square
\begin{align}
\begin{split}
  \xymatrix@=10pt{
    W \ar[r]^-{} \ar[d]_-{g} & Y \ar[d]^-{f} \\
    Z \ar[r]^-{} & X.
  }
\end{split}
\end{align}
where $f:Y\to X$ be a proper morphism between smooth connected $k$-schemes of same dimension, and $Z\to X$ is a closed immersion. Let $d=[k(Y):k(X)]$ be the degree of the extension between function fields.
Then the map $g_*g^!:CH_n(Z)\to CH_n(Z)$ agrees with the multiplication by $d$.
\end{lemma}
\proof
By the projection formula (\cite[Thm. 6.2]{Ful}), the map $g_*g^!$ agrees with the (refined) cup-product with the class $f_*\eta_Y$, where $\eta_Y$ is the fundamental class of $Y$. But $f_*\eta_Y$ is $d$ times the fundamental class of $X$, which implies the result.
\endproof
\end{remark}

\section{The boundary characteristic class}
In this section, we define the \emph{boundary characteristic class} in the boundary Chow group of $0$-cycles (Example~\ref{ex:bCH0}), by repeating the limiting process in Definition~\ref{df:proclass} for the class $C^{\mathrm{c}}_{U/X}(F)$ in~\eqref{eq:CcUX}. 

\subsection{}
\label{num:bCcset}
Let $k$ be a perfect field and let $X$ be a smooth $k$-scheme. Let $f:\overline{X}\to\overline{X}'$ be a morphism of smooth compactifications in $\operatorname{Cpt}^{\operatorname{Sm}}(X/k)$. We have the following Cartesian diagram:
\begin{align}
\begin{split}
  \xymatrix@=10pt{
    X \ar[r]^-{j} \ar@{=}[d]_-{} & \overline{X} \ar[d]^-{f} & \overline{X}-X \ar[d]^-{h} \ar[l]_-{k} \\
    X \ar[r]^-{j'} & \overline{X}' & \overline{X}'-X. \ar[l]_-{i}
  }
\end{split}
\end{align}
Let $K\in\mathbf{SH}_c(X)$ be a dualizable object. By~\eqref{eq:locpf}, since $C^{\mathrm{c}}_{X/\overline{X}'}(\mathbbold{1}_X)=0$ by definition, the push-forward map $h_*:CH_0(\overline{X}-X)[1/p]\to CH_0(\overline{X}'-X)[1/p]$ sends the class $C^{\mathrm{c}}_{X/\overline{X}}(K)$ to the class $C^{\mathrm{c}}_{X/\overline{X}'}(K)$. By~\ref{num:bBM}, we give the following definition:
\begin{definition}
\label{def:bCc}
Let $k$ be a perfect field and let $X$ be a smooth $k$-scheme. Assume that $\operatorname{Cpt}^{\operatorname{Sm}}(X/k)$ is cofinal in $\operatorname{Cpt}(X/k)$ (see~\ref{num:smcomp}). Let $K\in\mathbf{SH}_c(X)$ be a dualizable object. 
We define the \textbf{boundary characteristic class} $bC^{\mathrm{c}}_X(K)\in bCH_0(X)[1/p]$ to be the element associated to the formation of the classes 
$C^{\mathrm{c}}_{X/\overline{X}}(K)$ in ~\ref{eq:CcUX}: 
\begin{align}
bC^{\mathrm{c}}_X(K)
=
(C^{\mathrm{c}}_{X/\overline{X}}(K))_{\overline{X}\in\operatorname{Cpt}^{\operatorname{Sm}}(X)}
\in
bCH_0(X)[1/p].
\end{align}
If $k$ satisfies embedded resolution of singularities in the sense of~\ref{num:RS}, then the class $bC^{\mathrm{c}}_X(K)$ is defined in $bCH_0(X)$.
\end{definition}

\begin{remark}
If $k$ is perfect and $\operatorname{Cpt}^{\operatorname{Sm}}(X/k)$ is cofinal in $\operatorname{Cpt}(X/k)$, then by \cite[Cor. 3.1.6]{KSa}, $bCH_0(X)\simeq CH^{}_{0}(\overline{X}-X)$ for \emph{any} smooth compactification $\overline{X}$ of $X$, so the class $bC^{\mathrm{c}}_X(K)$ could have been directly defined in the group $CH_0(\overline{X}-X)$ for any such $\overline{X}$. However, formula~\eqref{eq:locpf} further shows that the class $bC^{\mathrm{c}}_X(K)$ is independent of the choice of the smooth compactification $\overline{X}$, and hence defines a meaningful invariant in the boundary Chow group.
\end{remark}

\begin{remark}
Definition~\ref{def:bCc} relies on the existence of smooth compactifications, which is related to resolution of singularities as discussed in~\ref{num:smcomp}. One may think of giving an alternative unconditional definition using alterations (\cite{Tem}), but we do not know how to do so, since Proposition~\ref{prop:CXZpf} in that case is insufficient to show that the class obtained is independent of the choice of alteration.
\end{remark}





\subsection{}
By~\eqref{eq:addCc}, for a distinguished triangle $L\to M\to N$ of dualizable objects in $\mathbf{SH}_c(X)$, we have
\begin{align}
\label{eq:addbCc}
bC^{\mathrm{c}}_X(L)=bC^{\mathrm{c}}_X(M)+bC^{\mathrm{c}}_X(N).
\end{align}
In particular, if we denote by $\mathbf{SH}^{\mathrm{rig}}_c(X)$ the full subcategory of $\mathbf{SH}_c(X)$ consisting of dualizable objects, we obtain a map
\begin{align}
\begin{split}
bC^{\mathrm{c}}_X:
K_0(\mathbf{SH}^{\mathrm{rig}}_c(X))&\to bCH_0(X)[1/p]. 
\end{split}
\end{align}
The same construction applies to \'etale sheaves, and by an analogue of~\eqref{eq:addbCc} in the filtered derived category \'etale sheaves we obtain a map
\begin{align}
\begin{split}
bC^{\mathrm{c}}_X:
K_0(\mathbf{D}^{b,\mathrm{rig}}_{ctf}(X,\Lambda))&\to b\mathbf{H}_{\mathrm{et}}\Lambda^{BM}_0(X/k). 
\end{split}
\end{align}
Since the \'etale realization is compatible with the six functors, 
we obtain a commutative square
\begin{align}
\begin{split}
  \xymatrix@=10pt{
    K_0(\mathbf{SH}^{\mathrm{rig}}_c(X)) \ar[r]^-{} \ar[d]_-{bC^{\mathrm{c}}_X} & K_0(\mathbf{D}^{b,\mathrm{rig}}_{ctf}(X,\Lambda)) \ar[d]^-{bC^{\mathrm{c}}_X} \\
    bCH_0(X)[1/p] \ar[r]^-{} & b\mathbf{H}_{\mathrm{et}}\Lambda^{BM}_0(X/k).
  }
\end{split}
\end{align}
If $k$ satisfies embedded resolution of singularities in the sense of~\ref{num:RS}, then we have a commutative square
\begin{align}
\begin{split}
  \xymatrix@=10pt{
    K_0(\mathbf{SH}^{\mathrm{rig}}_c(X)) \ar[r]^-{} \ar[d]_-{bC^{\mathrm{c}}_X} & K_0(\mathbf{D}^{b,\mathrm{rig}}_{ctf}(X,\Lambda)) \ar[d]^-{bC^{\mathrm{c}}_X} \\
    bCH_0(X) \ar[r]^-{} & b\mathbf{H}_{\mathrm{et}}\Lambda^{BM}_0(X/k).
  }
\end{split}
\end{align}

\subsection{}
\label{num:blC}
The canonical map $b\mathbb{E}^{\mathrm{BM}}_0(X/k)\xrightarrow{}l\mathbb{E}^{\mathrm{BM}}_0(X/k)$ in~\eqref{eq:btol} sends 
$bC^{\mathrm{c}}_X(K)$ to $lC^{\mathrm{c}}_X(K)-\operatorname{rk}(K)\cdot lC^{\mathrm{c}}_X(\mathbbold{1}_X)$. In particular we have
\begin{align}
\int bC^{\mathrm{c}}_X(K)=\chi_{\mathrm{c}}(K)-\operatorname{rk}(K)\cdot \chi_{\mathrm{c}}(\mathbbold{1}_X).
\end{align}
In other words, the degree of the compactly supported boudary characteristic class gives rise to the opposite of the classical Swan conductor (see \cite{AS}, \cite{KSa}, \cite{JY2}).

\subsection{}
By~\eqref{eq:locpf}, we obtain the following formula:
\begin{theorem}[Push-forward formula for the boundary characteristic class]
\label{th:pfbC}
Let $f:Y\to X$ be a finite \'etale morphism between smooth connected $k$-schemes. Assume that $\operatorname{Cpt}^{\operatorname{Sm}}(X/k)$ is cofinal in $\operatorname{Cpt}(X/k)$, and $\operatorname{Cpt}^{\operatorname{Sm}}(Y/k)$ is cofinal in $\operatorname{Cpt}(Y/k)$. Let $K\in\mathbf{SH}_c(Y)$ be a dualizable object. Then $f_*K\in\mathbf{SH}_c(X)$ is again dualizable, and we have
\begin{align}
\label{eq:bccpf}
bC^{\mathrm{c}}_X(f_*K)=f_*bC^{\mathrm{c}}_Y(K)+\operatorname{rk}(K)\cdot bC^{\mathrm{c}}_X(f_*\mathbbold{1})\in bCH_0(X)[1/p].
\end{align}
\end{theorem}

\section{Relation with the Kato-Saito Swan class}

\subsection{}
Throughout this section, we assume that $k$ is a perfect field of charactersitic $p$. We now recall Kato-Saito's construction of  Swan class (\cite[Def. 4.2.6]{KSa}), which is slightly refined thanks to recent results on alterations in \cite{Tem}.

\subsection{}
Let $f:V\to U$ be a finite \'etale morphism between smooth connected schemes over $k$, and let $Y$ be a compactification of $V$. Then by the theory of tame distillations (\cite{Tem}), there exists a commutative diagram
\begin{align}
\label{eq:KS3.4}
\begin{split}
  \xymatrix@=10pt{
    W \ar[r]^-{} \ar[d]_-{g} & Z \ar[rdd]^-{\overline{h}} \ar[d]_-{\overline{g}} & \\
    V \ar[r]^-{} \ar[d]_-{f} & Y & \\
    U \ar[rr]^-{} & & X 
  }
\end{split}
\end{align}
such that 
\begin{enumerate}
\item
all horizontal maps are open immersions, and both quadrangles are Cartesian
\item
$U$ is the complement of a Cartier divisor in $X$
\item
$Z$ is smooth and $W$ is the complement of a simple normal-crossing divisor in $Z$
\item
$\overline{g}$ is proper, surjective, generically finite of degree $p^r$ for some non-negative integer $r$.
\end{enumerate}
By \cite[Prop. 3.2.2]{KSa}, there exists a map 
\begin{align}
\label{eq:KS38}
(-,\Delta_{Z})^{\mathrm{log}}:
CH_d(W\times_UW\setminus W\times_VW)
\to 
CH_0(Z\setminus W)
\end{align}
called the \textbf{log intersection product}, which is the unique map such that the following diagram commutes:
\begin{align}
\begin{split}
  \xymatrix@R=14pt{
    CH_d(\overline{W\times_UW\setminus W\times_VW}) \ar[rd]^-{(-,\Delta_{Z})_{(Z\times Z)^\sim}} \ar[d]_-{} & \\
    CH_d(W\times_UW\setminus W\times_VW) \ar[r]_-{\eqref{eq:KS38}} & CH_0(Z\setminus W).
  }
\end{split}
\end{align}
Here 
\begin{enumerate}
\item
$(Z\times Z)^\sim$ is the log product (\cite[Def. 1.1.1]{KSa})
\item 
$\overline{W\times_UW\setminus W\times_VW}$ is the closure of $W\times_UW\setminus W\times_VW$ in $(Z\times Z)^\sim$
\item
the slant map is the intersection product with the diagonal of $Z$ in $(Z\times Z)^\sim$
\item
the vertical map is the canonical pull-back map.
\end{enumerate}


\begin{proposition}[\textrm{\cite[Thm. 3.2.3]{KSa}}]
Let $f:V\to U$ be a finite \'etale morphism between smooth connected schemes over $k$.
\begin{enumerate}
\item
There exists a unique map
\begin{align}
\label{eq:bKS}
(-,\Delta_{V})^{\mathrm{log}}_{\mathrm{b}}:
CH_d(V\times_UV\setminus \Delta_V)
\to 
bCH_0(V)[1/p]
\end{align}
such that for any commutative diagram of the form~\eqref{eq:KS3.4}, the following diagram commutes:
\begin{align}
\begin{split}
  \xymatrix@=14pt{
    CH_d(V\times_UV\setminus \Delta_V) \ar[r]^-{\eqref{eq:bKS}} \ar[d]_-{(g\times g)^!} & bCH_0(V)[1/p] \ar[r]^-{} & CH_0(Y\setminus V)[1/p] \\
    CH_d(W\times_UW\setminus W\times_VW) \ar[r]^-{\eqref{eq:KS38}} & CH_0(Z\setminus W) \ar[ru]_-{\frac{1}{p^r}\overline{g}_*} &
  }
\end{split}
\end{align}
\item
If $k$ satisfies embedded resolution of singularities in the sense of~\ref{num:RS}, then the map~\eqref{eq:bKS} factors through $bCH_0(V)$.
\end{enumerate}
\end{proposition}

\begin{definition}

\begin{enumerate}
\item
Let $f:V\to U$ be a finite Galois morphism between smooth connected schemes over $k$ of Galois group $G$. 
For $\sigma\in G$, we denote by $\Gamma_\sigma$ its graph and define the class $s_{V/U}(\sigma)\in bCH_0(V)[1/p]$ as
\begin{align}
s_{V/U}(\sigma)=
\begin{cases}
(V\times_UV\setminus\Delta_V,\Delta_{V})^{\mathrm{log}}_{\mathrm{b}} & \sigma=1 \\
-(\Gamma_\sigma,\Delta_{V})^{\mathrm{log}}_{\mathrm{b}} & \sigma\ne1.
\end{cases}
\end{align}
\item
Let $U$ be a smooth connected scheme over $k$ and let $\mathcal{F}$ be a smooth $\overline{\mathbb{F}}_\ell$-sheaf on $U$. Let $f:V\to U$ be a finite Galois morphism between smooth connected schemes over $k$ of Galois group $G$ which trivializes $\mathcal{F}$, such that $\mathcal{F}$ corresponds to a $\overline{F}_\ell$-representation of $G$ denoted as $M$. We denote by $G_{(p)}\subset G$ the subset of elements of order a power of $p$, and define the class 
\begin{align}
\operatorname{Sw}_{V/U}(\mathcal{F})=
\sum_{\sigma\in G_{(p)}}
\left(
\operatorname{dim}_{\mathbb{F}_\ell}M^{\sigma}
-
\frac{\operatorname{dim}_{K_\ell}M^{\sigma^p}/M^{\sigma}}{p-1}
\right)
s_{V/U}(\sigma)
\in
bCH_0(V)[1/p].
\end{align}
\item
Let $U$ be a smooth connected scheme over $k$ and let $\mathcal{F}$ be a smooth $\overline{\mathbb{F}}_\ell$-sheaf on $U$. If $f:V\to U$ is a finite Galois morphism between smooth connected schemes over $k$ of Galois group $G$ which trivializes $\mathcal{F}$, the \textbf{Swan class} $\operatorname{Sw}_U(\mathcal{F})$ as
\begin{align}
\operatorname{Sw}_U(\mathcal{F})=\frac{1}{|G|}f_*\operatorname{Sw}_{V/U}(\mathcal{F}) \in bCH_0(U)[1/p]
\end{align}
which is independent of the choice of $V$.
\end{enumerate}
\end{definition}

\subsection{}
The Kato-Saito construction gives rise to a map
\begin{align}
\label{eq:Swanmap}
\begin{split}
\operatorname{Sw}_U:
K_0(\mathbf{D}^{b,\mathrm{rig}}_{ctf}(U,\Lambda))&\to bCH_0(U)[1/p]. 
\end{split}
\end{align}
If $k$ satisfies embedded resolution of singularities in the sense of~\ref{num:RS}, then the map~\eqref{eq:Swanmap} lifts to $bCH_0(U)$.

\subsection{}
We conjecture that the Kato-Saito Swan class agrees with the boundary characteristic class:
\begin{conjecture}
\label{conj:b=Sw}

Let $U$ be a smooth connected scheme over a perfect field $k$. 
\begin{enumerate}
\item
Assume that the subcategory $\operatorname{Cpt}^{\operatorname{Sm}}(U/k)$ is cofinal in $\operatorname{Cpt}(U/k)$ (see~\ref{num:smcomp}).
Then the following diagram commutes:
\begin{align}
\begin{split}
  \xymatrix@=12pt{
    K_0(\mathbf{SH}^{\mathrm{rig}}_c(U)) \ar[r]^-{} \ar[d]_-{bC^{\mathrm{c}}_U} & K_0(\mathbf{D}^{b,\mathrm{rig}}_{ctf}(U,\Lambda)) \ar[d]^-{bC^{\mathrm{c}}_U} \ar[ld]_-{-\operatorname{Sw}^{\mathrm{}}_U} \\
    bCH_0(U)[1/p] \ar[r]^-{} & b\mathbf{H}_{\mathrm{et}}\Lambda^{BM}(U/k).
  }
\end{split}
\end{align}
\item
Assume that $k$ satisfies embedded resolution of singularities in the sense of~\ref{num:RS}. Then the following diagram commutes:
\begin{align}
\label{eq:diagres}
\begin{split}
  \xymatrix@=12pt{
    K_0(\mathbf{SH}^{\mathrm{rig}}_c(U)) \ar[r]^-{} \ar[d]_-{bC^{\mathrm{c}}_U} & K_0(\mathbf{D}^{b,\mathrm{rig}}_{ctf}(U,\Lambda)) \ar[d]^-{bC^{\mathrm{c}}_U} \ar[ld]_-{-\operatorname{Sw}^{\mathrm{}}_U} \\
    bCH_0(U) \ar[r]^-{} & b\mathbf{H}_{\mathrm{et}}\Lambda^{BM}(U/k).
  }
\end{split}
\end{align}
\end{enumerate}
\end{conjecture}

\subsection{}
The following formulas provide some evidence for these conjectures: if $f$ is a finite \'etale morphism, then both classes satisfy the formula
\begin{align}
\label{eq:bccpfnew}
bC^{\mathrm{c}}_X(f_*\mathcal{F})=f_*bC^{\mathrm{c}}_Y(\mathcal{F})+\operatorname{rk}(\mathcal{F})\cdot bC^{\mathrm{c}}_X(f_*\mathbbold{1})\in b\mathbf{H}_{\mathrm{et}}\Lambda^{BM}_0(X/k).
\end{align}
\begin{align}
\label{eq:swpf}
\operatorname{Sw}_X(f_*\mathcal{F})=f_*\operatorname{Sw}_Y(\mathcal{F})+\operatorname{rk}(\mathcal{F})\cdot \operatorname{Sw}_X(f_*\mathbbold{1})\in bCH_0(X).
\end{align}
Indeed, the first formula is~\eqref{eq:bccpf}, and the second is \cite[Cor. 4.3.4]{KSa}.

\section{Proof of Proposition~\ref{prop:CXZpf}}
\label{sec:pfpbpf}

\subsection{}
\label{num:lmCXUZ}
We now begin the proof of Proposition~\ref{prop:CXZpf}. 
Recall that in Lemma~\ref{lm:AS524}, we introduced a class $C_{X,U}^Z(F)\in [i_*\mathbbold{1}_Z,j_!j^!q^!\mathbf{H}\mathbb{Z}]$. By~\ref{num:CXUZim}, the class $C^{\mathrm{c}}_{U/X}(F)$ is nothing but the image of $C_{X,U}^Z(F)$ in $CH_0(Z)[1/p]$. So in order to prove formula 
~\eqref{eq:locpf} it suffices to prove the following one:
\begin{align}
\label{eq:locpf2}
C_{X,U}^Z(g_*F)
=
h_*C_{Y,V}^W(F)
+
\operatorname{rk}(F)\cdot C_{X,U}^Z(g_*\mathbbold{1}_V).
\end{align}

\subsection{}
Recall that the class $C_{X,U}^Z(F)$ is defined in Lemma~\ref{lm:AS524} by lifting the class 
\begin{align}
\beta(F)=
B_{X,U}^Z(F,\mathbf{H}\mathbb{Z})
-
\operatorname{rk}(F)\cdot B_{X,U}^Z(\mathbbold{1}_U,\mathbf{H}\mathbb{Z})
\in
[i_*\mathbbold{1}_Z,\delta^\Delta_X\delta_{X*}j_!j^!f^!\mathbf{H}\mathbb{Z}]
\end{align}
in~\eqref{eq:beta}. We now construct some transformations of the functor $\delta^\Delta_X$ introduced in~\eqref{eq:deltaD}.

\subsection{}
\label{num:dDfunc}
Let $f:Y\to X$ be a morphism of schemes. Let $\gamma_X:X\times X-X\to X\times X$ and $\gamma_Y:Y\times Y-Y\to Y\times Y$ be the inclusion of the complement of the diagonal. Let $\gamma_1:Y\times Y-Y\times_XY\to Y\times Y-Y$ be the canonical open immersion. Then we have a commutative square
\begin{align}
\begin{split}
  \xymatrix@R=10pt{
    Y\ar[r]^-{\delta_{Y/X}} \ar[rd]_-{f} & Y\times_XY \ar[d]^-{f_1} \ar[r]^-{\delta_1} & Y\times Y \ar[d]^-{f\times f} & Y\times Y-Y\times_XY \ar[l]_-{\gamma_Y\circ\gamma_1} \ar[d]^-{} \\
      & X \ar[r]^-{\delta_X} & X\times X & X\times X-X \ar[l]_-{\gamma_X}
  }
\end{split}
\end{align}
where both squares are Cartesian. By adjunction we have a map (see \cite[1.1.15]{CD1})
\begin{align}
\label{eq:deltaf*}
\delta_X^*(f\times f)_*
\to 
f_*\delta_Y^*.
\end{align}
There is also a map
\begin{align}
\label{eq:deltaf!}
f_*\delta_Y^!
=
f_{1*}\delta_{Y/X*}\delta_{Y/X}^!\delta_1^!
\to
f_{1*}\delta_1^!
\simeq
\delta_X^!(f\times f)_*.
\end{align}
The following square is commutative:
\begin{align}
\label{eq:delta*!com}
\begin{split}
  \xymatrix@R=10pt{
f_*\delta_Y^! \ar[r]^-{\eqref{eq:deltaf*}} \ar[d]^-{} & \delta_X^!(f\times f)_* \ar[d]^-{}  \\
f_*\delta_Y^*  & \delta_X^*(f\times f)_* \ar[l]_-{\eqref{eq:deltaf!}}.
  }
\end{split}
\end{align}

\subsection{}
In the setting of~\ref{num:dDfunc}, if $f$ is smooth, we have a map
\begin{align}
\label{eq:gammaXY}
\gamma_{Y*}\mathbbold{1}
\to
\gamma_{Y*}\gamma_{1*}\mathbbold{1}
\simeq
(f\times f)^*\gamma_{X*}\mathbbold{1}
\end{align}
and we deduce a map
\begin{align}
\label{eq:Deltadtr}
\delta_Y^*\gamma_{Y*}\mathbbold{1}
&\xrightarrow{\eqref{eq:gammaXY}}
\delta_Y^*(f\times f)^*\gamma_{X*}\mathbbold{1}
=
f^*\delta_X^*\gamma_{X*}\mathbbold{1} 
\end{align}
which fits into a morphism of distinguished triangles
\begin{align}
\begin{split}
  \xymatrix@=10pt{
    \delta_Y^!\mathbbold{1} \ar[r]^-{} \ar[d]_-{} & \delta_Y^*\mathbbold{1} \ar[r]^-{} \ar@{=}[d]^-{} & \delta_Y^*\gamma_{Y*}\mathbbold{1} \ar[d]^-{\eqref{eq:Deltadtr}}  \\
    f^*\delta_X^!\mathbbold{1} \ar[r]^-{} & f^*\delta_X^*\mathbbold{1} \ar[r]^-{} & f^*\delta_X^*\gamma_{X*}\mathbbold{1}.
  }
\end{split}
\end{align}
We define a natural transformation
\begin{align}
\label{eq:dDu!}
\delta^\Delta_Y\delta_{Y*}f^!
\to
f^!\delta^\Delta_X\delta_{X*},
\end{align}
by the composition
\begin{align}
\begin{split}
\delta^\Delta_Y\delta_{Y*}f^!(K)
&\overset{\eqref{eq:deltad}}{\simeq}
f^!(K)\otimes\delta_{Y}^*\gamma_{Y*}\mathbbold{1}
\xrightarrow{\eqref{eq:Deltadtr}}
f^!(K)\otimes f^*\delta_X^*\gamma_{X*}\mathbbold{1}\\
&\to
f^!(K\otimes\delta_X^*\gamma_{X*}\mathbbold{1})
\overset{\eqref{eq:deltad}}{\simeq}
f^!\delta^\Delta_X\delta_{X*}K.
\end{split}
\end{align}
Similarly we define the following transformations
\begin{align}
\label{eq:dDu*}
\delta^\Delta_Y\delta_{Y*}f^*
\to 
f^*\delta^\Delta_X\delta_{X*},
\end{align}
\begin{align}
\label{eq:dDl!}
f_!\delta^\Delta_Y\delta_{Y*}
\to
\delta^\Delta_X\delta_{X*}f_!.
\end{align}

If $f$ is smooth and $X$ is smooth, then the object $\delta_X^*\gamma_{X*}\mathbbold{1}$ is dualizable, and we define a natural transformation
\begin{align}
\label{eq:dDl*}
f_*\delta^\Delta_Y\delta_{Y*}
\to
\delta^\Delta_X\delta_{X*}f_*,
\end{align}
by the composition
\begin{align}
\begin{split}
f_*\delta^\Delta_Y\delta_{Y*}(K)
&\overset{\eqref{eq:deltad}}{\simeq}
f_*(K\otimes\delta_{Y}^*\gamma_{Y*}\mathbbold{1})
\xrightarrow{\eqref{eq:Deltadtr}}
f_*(K\otimes f^*\delta_X^*\gamma_{X*}\mathbbold{1})\\
&\simeq
f_*K\otimes\delta_X^*\gamma_{X*}\mathbbold{1}
\overset{\eqref{eq:deltad}}{\simeq}
\delta^\Delta_X\delta_{X*}f_*K.
\end{split}
\end{align}

If $f$ is \'etale morphism between smooth schemes, the map $\delta_Y^!\mathbbold{1}\to f^*\delta_X^!\mathbbold{1}$ is an isomorphism, and therefore by~\eqref{eq:abspur} and the five lemma, the map~\eqref{eq:Deltadtr} is an isomorphism of dualizable objects. It follows that the transformations~\eqref{eq:dDu!},~\eqref{eq:dDu*},~\eqref{eq:dDl!} and~\eqref{eq:dDl*} are all isomorphisms.

\subsection{}
Let $j:U\to X$ be an open immersion of smooth schemes. We have a commutative diagram
\begin{align}
\begin{split}
  \xymatrix{
    U\times U \ar[r]^-{j\times1_U} \ar[d]_-{1_U\times j} \ar[rd]^-{j\times j} & X\times U \ar[r]^-{p_{U}^{XU}} \ar[d]^-{1_X\times j} & U \ar[d]^-{j} \\
    U\times X \ar[r]_-{j\times1_X} \ar[d]_-{p_{U}^{UX}} & X\times X \ar[r]_-{p_{2X}} \ar[d]^-{p_{1X}} & X \ar[d]^-{} \\
    U \ar[r]^-{j} & X \ar[r]^-{} & k
  }
\end{split}
\end{align}
where all the four squares are Cartesian. The compositions $p_{1U}=p_U^{UX}\circ(1_U\times j):U\times U\to U$ and $p_{1U}=p_U^{XU}\circ(j\times1_U):U\times U\to U$ are respectively the projections to the two components. 

Let $\delta_X:X\to X\times X$ be the diagonal with open complement $\gamma_X:X\times X-X\to X\times X$. Let $F,G\in\mathbf{SH}(U)$ be two objects. There is a canonical map
\begin{align}
\label{eq:ppf}
\begin{split}
&(j_*G\boxtimes j_!F)\otimes\gamma_{X*}\mathbbold{1}
=
p_{1X}^*j_*G\otimes p_{2X}^*j_!F\otimes\gamma_{X*}\mathbbold{1}\\
\simeq
&(j\times1_X)_*p_{U}^{UX*}G\otimes (1_X\times j)_!p_{U}^{XU*}F\otimes\gamma_{X*}\mathbbold{1}\\
\simeq
&(1_X\times j)_!((1_X\times j)^*(j\times1_X)_*p_{U}^{UX*}G\otimes p_{U}^{XU*}F\otimes(1_X\times j)^*\gamma_{X*}\mathbbold{1})\\
\simeq
&(1_X\times j)_!((j\times1_U)_*p_{1U}^{*}G\otimes p_{U}^{XU*}F\otimes(1_X\times j)^*\gamma_{X*}\mathbbold{1})\\
\to
&(1_X\times j)_!(j\times1_U)_*(p_{1U}^{*}G\otimes(j\times1_U)^*p_{U}^{XU*}F\otimes(j\times1_U)^*(1_X\times j)^*\gamma_{X*}\mathbbold{1})\\
=
&(1_X\times j)_!(j\times1_U)_*((G\boxtimes_kF)\otimes(j\times j)^*\gamma_{X*}\mathbbold{1}). 
\end{split}
\end{align}

\begin{lemma}
\label{lm:ppf}
The map~\eqref{eq:ppf} is an isomorphism.
\end{lemma}
\proof
It suffices to show that 
the following canonical map is an isomorphism:
\begin{align}
\begin{split}
&(j\times1_U)_*p_{1U}^{*}G\otimes p_{U}^{XU*}F\otimes(1_X\times j)^*\gamma_{X*}\mathbbold{1}\\
\to
&(j\times1_U)_*(p_{1U}^{*}G\otimes p_{2U}^{*}F\otimes(j\times j)^*\gamma_{X*}\mathbbold{1}).
\end{split}
\end{align}
From~\ref{num:FHM32}, the distinguished triangle~\eqref{eq:gammadt} and the five lemma, we see that is suffices to show the following canonical map is an isomorphism:
\begin{align}
\label{eq:ppfint}
\begin{split}
&(j\times1_U)_*p_{1U}^{*}G\otimes p_{U}^{XU*}F\otimes(1_X\times j)^*\delta_{X*}\delta_X^!\mathbbold{1}\\
\to
&(j\times1_U)_*(p_{1U}^{*}G\otimes p_{2U}^{*}F\otimes(j\times j)^*\delta_{X*}\delta_X^!\mathbbold{1}).
\end{split}
\end{align}
We have a commutative diagram
\begin{align}
\begin{split}
  \xymatrix{
    U \ar[r]^-{\delta_U} \ar@{=}[d]_-{} & U\times U \ar[r]^-{p_{2U}} \ar[d]_-{j\times1_U} & U \ar@{=}[d]_-{} \\
    U \ar[r]^-{(j,1_U)} \ar[d]_-{j} & X\times U \ar[r]^-{p_{U}^{XU}} \ar[d]_-{1_X\times j} & U \ar[d]^-{j} \\
    X \ar[r]^-{\delta_X} & X\times X \ar[r]^-{p_{2X}} & X
  }
\end{split}
\end{align}
where all the four squares are Cartesian. We deduce an isomorphism
\begin{align}
\label{eq:ppfiso}
\begin{split}
&(j\times1_U)_*p_{1U}^{*}G\otimes p_{U}^{XU*}F\otimes(1_X\times j)^*\delta_{X*}\delta_X^!\mathbbold{1}\\
\simeq
&(j\times1_U)_*p_{1U}^{*}G\otimes p_{U}^{XU*}F\otimes(j,1_U)_*j^*\delta_X^!\mathbbold{1}\\
\simeq
&(j,1_U)_*((j,1_U)^*(j\times1_U)_*p_{1U}^{*}G\otimes (j,1_U)^*p_{U}^{XU*}F\otimes j^*\delta_X^!\mathbbold{1})\\
\simeq
&(j,1_U)_*(G\otimes F\otimes j^*\delta_X^!\mathbbold{1})
=
(j\times 1_U)_*\delta_{U*}(G\otimes F\otimes j^*\delta_X^!\mathbbold{1})\\
\simeq
&(j\times1_U)_*(p_{1U}^{*}G\otimes p_{2U}^{*}F\otimes \delta_{U*}j^*\delta_X^!\mathbbold{1})\\
\simeq
&(j\times1_U)_*(p_{1U}^{*}G\otimes p_{2U}^{*}F\otimes (j\times j)^*\delta_{X*}\delta_X^!\mathbbold{1}).
\end{split}
\end{align}
One can check that the map~\eqref{eq:ppfint} agrees with the map~\eqref{eq:ppfiso}, and the result follows.
\endproof





\begin{lemma}
\label{lm:pbpfb}
In the setting of~\ref{num:birset}, assume that $f$ is proper and $g$ is \'etale (and finite). Let $F\in\mathbf{SH}_c(V)$ be a dualizable object. 
Then the following diagram commutes:
\begin{align}
\label{eq:diagpfB}
\begin{split}
  \xymatrix@R=10pt{
    i_*\mathbbold{1}_Z \ar[d]^-{} \ar[r]^-{B_{X,U}^Z(g_*F)} & \delta^\Delta_X\delta_{X*}j_!j^!q^!\mathbf{H}\mathbb{Z}  \\
    f_*k_*\mathbbold{1}_W \ar[r]^-{f_*B_{Y,V}^W(F)} & f_*\delta^\Delta_Y\delta_{Y*}l_!l^!f^!q^!\mathbf{H}\mathbb{Z}. \ar[u]^-{}
  }
\end{split}
\end{align}
Here the right vertical map arises from the map
\begin{align}
\label{eq:Deltal!fin}
\begin{split}
f_*\delta^\Delta_Y\delta_{Y*}l_!l^!f^!
&\overset{\eqref{eq:dDl!}}{\simeq}
f_*l_!\delta^\Delta_V\delta_{V*}l^!f^!
\simeq
f_*l_!\delta^\Delta_V\delta_{V*}g^!j^!
\xrightarrow{\eqref{eq:dDu!}}
f_*l_!g^!\delta^\Delta_U\delta_{U*}j^!\\
&\simeq
j_!g_!g^!\delta^\Delta_U\delta_{U*}j^!
\xrightarrow{}
j_!\delta^\Delta_U\delta_{U*}j^!
\overset{\eqref{eq:dDl!}}{\simeq}
\delta^\Delta_X\delta_{X*}j_!j^!.
\end{split}
\end{align}
\end{lemma}
\proof
We are reduced to show the commutativity of the following diagram:
\begin{align}
\label{eq:bigdiag}
\begin{split}
  \xymatrix@=10pt{
    i_*\mathbbold{1}_Z \ar[rr]^-{} \ar[d]_-{} & & f_*k_*\mathbbold{1}_W \ar[d]^-{} \\
    i_*i^*j_*\mathbbold{1}_U \ar[r]^-{} \ar[d]_-{} & i_*i^*f_*l_*\mathbbold{1}_V \ar[r]^-{\sim} \ar[d]_-{} & f_*k_*k^*l_*\mathbbold{1}_V \ar[d]^-{} \\
    i_*i^*\delta_X^\Delta\delta_{X*}j_*\mathbbold{1}_U \ar[r]^-{(a)} \ar[d]_-{} & i_*i^*f_*\delta_Y^\Delta\delta_{Y*}l_*\mathbbold{1}_V \ar[r]^-{\sim} \ar[d]_-{} & f_*k_*k^*\delta_Y^\Delta\delta_{Y*}l_*\mathbbold{1}_V \ar[d]^-{} \\
    i_*i^*\delta_X^\Delta\delta_{X*}j_*\underline{Hom}(g_*F,g_*F) \ar[r]^-{(b)} \ar[d]_-{} & i_*i^*f_*\delta_Y^\Delta\delta_{Y*}l_*\underline{Hom}(F,F) \ar[r]^-{\sim} \ar[d]_-{} & f_*k_*k^*\delta_Y^\Delta\delta_{Y*}l_*\underline{Hom}(F,F) \ar[d]^-{} \\
    i_*i^*\delta_X^\Delta(\mathbb{D}_{X/k}(j_!g_*F)\boxtimes_kj_!g_*F) & i_*i^*f_*\delta_Y^\Delta(\mathbb{D}_{Y/k}(l_!F)\boxtimes_kl_!F)  \ar[r]^-{\sim} \ar[l]_-{(c)} & f_*k_*k^*\delta_Y^\Delta(\mathbb{D}_{Y/k}(l_!F)\boxtimes_kl_!F)  \\
    \delta_X^\Delta(\mathbb{D}_{X/k}(j_!g_*F)\boxtimes_kj_!g_*F) \ar[u]^-{\wr} \ar[d]_-{} & f_*\delta_Y^\Delta(\mathbb{D}_{Y/k}(l_!F)\boxtimes_kl_!F) \ar@{=}[r]^-{}  \ar[u]^-{\wr} \ar[l]_     -{(d)} \ar@{}[d]|{\Gamma} & f_*\delta_Y^\Delta(\mathbb{D}_{Y/k}(l_!F)\boxtimes_kl_!F)  \ar[u]_-{\wr} \ar[d]^-{} \\
    \delta_X^\Delta\delta_{X*}j_!(\mathbb{D}_{U/k}(g_*F)\otimes g_*F)  \ar[d]_-{} \ar[rr]_-{(e)} & & f_*\delta_Y^\Delta\delta_{Y*}l_!(\mathbb{D}_{V/k}(F)\otimes F) \ar[d]^-{}  \\
    \delta_X^\Delta\delta_{X*}j_!j^!q^!\mathbbold{1}_k & & f_*\delta_Y^\Delta\delta_{Y*}l_!l^!f^!q^!\mathbbold{1}_k.  \ar[ll]_-{\eqref{eq:Deltal!fin}}
  }
\end{split}
\end{align}
We explain how the following maps in diagram~\eqref{eq:bigdiag} are defined, while the other unmarked maps are defined using adjunctions:
\begin{itemize}
\item
The map $(a)$ is induced by the composition 
\begin{align}
\label{eq:deltal*}
\begin{split}
\delta_X^\Delta\delta_{X*}j_*
\overset{\eqref{eq:dDl*}}{\simeq}
j_*\delta_U^\Delta\delta_{U*}
\to
j_*g_*g^*\delta_U^\Delta\delta_{U*}
\overset{\eqref{eq:dDu*}}{\simeq}
f_*l_*\delta_V^\Delta\delta_{V*}g^*
\overset{\eqref{eq:dDu*}}{\simeq}
f_*\delta_Y^\Delta\delta_{Y*}l_*g^*.
\end{split}
\end{align}
\item
The map $(b)$ is induced by the composition 
\begin{align}
\begin{split}
&\delta_X^\Delta\delta_{X*}j_*\underline{Hom}(g_*F,g_*F)
\xrightarrow{\eqref{eq:deltal*}}
f_*\delta_Y^\Delta\delta_{Y*}l_*g^*\underline{Hom}(g_*F,g_*F)\\
\simeq
&f_*\delta_Y^\Delta\delta_{Y*}l_*g^*g_*\underline{Hom}(F,g^!g_*F)
\to
f_*\delta_Y^\Delta\delta_{Y*}l_*\underline{Hom}(F,F).
\end{split}
\end{align}
\item
In order to define the maps $(c)$ and $(d)$, we define a map
\begin{align}
\label{eq:hardmap}
\begin{split}
&(\mathbb{D}_{Y/k}(l_!F)\boxtimes_kl_!F)\otimes\gamma_{Y*}\mathbbold{1}\\
\overset{\eqref{eq:ppf}}{\simeq}
&((1_Y\times l)_!(l\times 1_V)_*((\mathbb{D}_{V/k}(F)\boxtimes_kF)\otimes(l\times l)^*\gamma_{Y*}\mathbbold{1})\\
\simeq
&(1_Y\times l)_!(l\times 1_V)_*((\mathbb{D}_{V/k}(F)\boxtimes_kF)\otimes\gamma_{V*}\mathbbold{1})\\
\xrightarrow{\eqref{eq:gammaXY}}
&(1_Y\times l)_!(l\times 1_V)_*((\mathbb{D}_{V/k}(F)\boxtimes_kF)\otimes (g\times g)^*\gamma_{U*}\mathbbold{1})\\
\simeq
&(1_Y\times l)_!(l\times 1_V)_*((\mathbb{D}_{V/k}(F)\boxtimes_kF)\otimes (l\times l)^*(f\times f)^*\gamma_{X*}\mathbbold{1})\\
\overset{\eqref{eq:ppf}}{\simeq}
&(\mathbb{D}_{Y/k}(l_!F)\boxtimes_kl_!F)\otimes (f\times f)^*\gamma_{X*}\mathbbold{1}.
\end{split}
\end{align}
With the notations in~\ref{num:dDfunc}, the maps $(c)$ and $(d)$ are given by the composition 
\begin{align}
\label{eq:comd}
\begin{split}
&f_*\delta_Y^\Delta(\mathbb{D}_{Y/k}(l_!F)\boxtimes_kl_!F)\\
=
&f_*\delta_Y^!(\mathbb{D}_{Y/k}(l_!F)\boxtimes_kl_!F\otimes\gamma_{Y*}\mathbbold{1})\\
\xrightarrow{\eqref{eq:deltaf!}}
&\delta_X^!(f\times f)_*(\mathbb{D}_{Y/k}(l_!F)\boxtimes_kl_!F\otimes\gamma_{Y*}\mathbbold{1})\\
\xrightarrow{\eqref{eq:hardmap}}
&\delta_X^!(f\times f)_*(\mathbb{D}_{Y/k}(l_!F)\boxtimes_kl_!F\otimes(f\times f)^*\gamma_{X*}\mathbbold{1})\\
\simeq
&\delta_X^!((f\times f)_*(\mathbb{D}_{Y/k}(l_!F)\boxtimes_kl_!F)\otimes\gamma_{X*}\mathbbold{1})\\
=
&\delta_X^\Delta(f\times f)_*(\mathbb{D}_{Y/k}(l_!F)\boxtimes_kl_!F)\\
\simeq
&\delta_X^\Delta(\mathbb{D}_{X/k}(j_!g_*F)\boxtimes_kj_!g_*F).
\end{split}
\end{align}
\item
The map $(e)$ is the composition 
\begin{align}
\begin{split}
&\delta_X^\Delta\delta_{X*}j_!(\mathbb{D}_{U/k}(g_*F)\otimes g_*F)\\
\overset{\eqref{eq:dDl!}}{\simeq}
&j_!\delta_U^\Delta\delta_{U*}(\mathbb{D}_{U/k}(g_*F)\otimes g_*F)\\
\simeq
&j_!\delta_U^\Delta\delta_{U*}g_*(g^*g_*\mathbb{D}_{V/k}(F)\otimes F)\\
\to
&j_!\delta_U^\Delta\delta_{U*}g_*(\mathbb{D}_{V/k}(F)\otimes F)\\
\overset{\eqref{eq:dDl*}}{\simeq}
&j_!g_*\delta_V^\Delta\delta_{V*}(\mathbb{D}_{V/k}(F)\otimes F)\\
\simeq
&f_*l_!\delta_V^\Delta\delta_{V*}(\mathbb{D}_{V/k}(F)\otimes F)\\
\overset{\eqref{eq:dDl!}}{\simeq}
&f_*\delta_Y^\Delta\delta_{Y*}l_!(\mathbb{D}_{V/k}(F)\otimes F).
\end{split}
\end{align}
\end{itemize}

Since all maps in diagram~\eqref{eq:bigdiag} are deduced from adjunctions of the six functors, the commutativity of each square can be checked using compatibilities between these operations. For example, we now check the commutativity of the square $\Gamma$ in diagram~\eqref{eq:bigdiag}.
Using the diagram
\begin{align*}
\begin{split}
\resizebox{1\textwidth}{!}{
  \xymatrix@=10pt{
    f_*\delta_Y^\Delta(\mathbb{D}_{Y/k}(l_!F)\boxtimes_kl_!F) \ar[rr]^-{} \ar[d]^-{}  &   & f_*\delta_Y^\Delta\delta_{Y*}l_!(\mathbb{D}_{V/k}(F)\otimes F) \ar[ld]^-{\sim} \\
    f_*\delta_Y^!((\mathbb{D}_{Y/k}(l_!F)\boxtimes_kl_!F)\otimes\gamma_{Y*}\mathbbold{1}) \ar[r]^-{} \ar[d]_-{} \ar@{}[rd]|{\eqref{eq:delta*!com}} & f_*\delta_Y^*((\mathbb{D}_{Y/k}(l_!F)\boxtimes_kl_!F)\otimes\gamma_{Y*}\mathbbold{1}) & \\
    \delta_X^!(f\times f)_*((\mathbb{D}_{Y/k}(l_!F)\boxtimes_kl_!F)\otimes\gamma_{Y*}\mathbbold{1}) \ar[r]^-{} \ar[d]_-{} & \delta_X^*(f\times f)_*((\mathbb{D}_{Y/k}(l_!F)\boxtimes_kl_!F)\otimes\gamma_{Y*}\mathbbold{1}) \ar[u]_-{} \ar[d]_-{} \ar@{}[ru]|{\ \ \ \ \ \ \ \ \Gamma_1} & \\
    \delta_X^!(f\times f)_*((\mathbb{D}_{Y/k}(l_!F)\boxtimes_kl_!F)\otimes(f\times f)^*\gamma_{X*}\mathbbold{1}) \ar[r]^-{} & \delta_X^*(f\times f)_*((\mathbb{D}_{Y/k}(l_!F)\boxtimes_kl_!F)\otimes(f\times f)^*\gamma_{X*}\mathbbold{1}) \ar[r]^-{\sim} & \delta_X^\Delta\delta_{X*}j_!(\mathbb{D}_{U/k}(g_*F)\otimes g_*F) \ar[uuu]_-{(e)}
  }
}
\end{split}
\end{align*}
we see that is suffices to check the commutativity of the subdiagram $\Gamma_1$. This follows from the following commutative diagram:
\begin{align*}
\begin{split}
\resizebox{1\textwidth}{!}{
  \xymatrix@=10pt{
    \delta_X^*(f\times f)_*((\mathbb{D}_{Y/k}(l_!F)\boxtimes_kl_!F)\otimes\gamma_{Y*}\mathbbold{1}) \ar[r]^-{} \ar[d]_-{\wr} & f_*\delta_Y^*((\mathbb{D}_{Y/k}(l_!F)\boxtimes_kl_!F)\otimes\gamma_{Y*}\mathbbold{1}) \ar[d]^-{\wr} \\
    \delta_X^*(f\times f)_*(1_Y\times l)_!(l\times 1_V)_*((\mathbb{D}_{V/k}(F)\boxtimes_kF)\otimes\gamma_{V*}\mathbbold{1}) \ar[r]_-{} \ar[d]_-{} & f_*\delta_Y^*(1_Y\times l)_!(l\times 1_V)_*((\mathbb{D}_{V/k}(F)\boxtimes_kF)\otimes\gamma_{V*}\mathbbold{1}) \ar[d]^-{\wr} \\
    \delta_X^*(f\times f)_*(1_Y\times l)_!(l\times 1_V)_*((\mathbb{D}_{V/k}(F)\boxtimes_kF)\otimes (g\times g)^*\gamma_{U*}\mathbbold{1}) \ar[r]_-{}  & f_*\delta_Y^*(1_Y\times l)_!(l\times 1_V)_*((\mathbb{D}_{V/k}(F)\boxtimes_kF)\otimes (g\times g)^*\gamma_{U*}\mathbbold{1}) \ar[d]^-{\wr} \\
    \delta_X^*(f\times f)_*((\mathbb{D}_{Y/k}(l_!F)\boxtimes_kl_!F)\otimes (f\times f)^*\gamma_{X*}\mathbbold{1}) \ar[r]_-{}  \ar[d]_-{\wr} \ar[u]^-{\wr} & f_*\delta_Y^*((\mathbb{D}_{Y/k}(l_!F)\boxtimes_kl_!F)\otimes (f\times f)^*\gamma_{X*}\mathbbold{1}) \ar[d]^-{\wr} \ar[u]_-{\wr} \\
    j_!(\mathbb{D}_{U/k}(g_*F)\otimes g_*F)\otimes\delta_X^*\gamma_{X*}\mathbbold{1} \ar[r]_-{} \ar[d]_-{\wr} & j_!g_*(\mathbb{D}_{V/k}(F)\otimes F)\otimes\delta_X^*\gamma_{X*}\mathbbold{1} \ar[d]^-{\wr} \\
    \delta_X^\Delta\delta_{X*}j_!(\mathbb{D}_{U/k}(g_*F)\otimes g_*F) \ar[r]^-{(e)}  & f_*\delta_Y^\Delta\delta_{Y*}l_!(\mathbb{D}_{V/k}(F)\otimes F)
  }
}
\end{split}
\end{align*}
and the proof is finished.
\endproof

\begin{remark}
In the construction of the map~\eqref{eq:hardmap}, we use Lemma~\ref{lm:ppf} reduce to the open part $g:V\to U$ which is assumed \'etale in order to use the map~\eqref{eq:gammaXY}. In the case where $f$ is assumed \'etale, the construction of the map~\eqref{eq:hardmap} is much simpler.
\end{remark}

\noindent
\emph{Proof of Proposition~\ref{prop:CXZpf}.} 
Denote by $\nu$ the canonical map
\begin{align}
\nu:[i_*\mathbbold{1}_Z,j_!j^!f^!\mathbf{H}\mathbb{Z}]
\to
[i_*\mathbbold{1}_Z,\delta^\Delta_X\delta_{X*}j_!j^!f^!\mathbf{H}\mathbb{Z}].
\end{align}
By construction of the class $C_{X,U}^Z(F)$, we know that
\begin{itemize}
\item
The image of the class $h_*C_{Y,V}^W(G)$ by $\nu$ is the class
\begin{align}
h_*(B_{Y,V}^W(G,\mathbf{H}\mathbb{Z})
-
\operatorname{rk}(G)\cdot B_{Y,V}^W(\mathbbold{1}_V,\mathbf{H}\mathbb{Z})).
\end{align}
\item
The class $C_{X,U}^Z(g_*G)-\operatorname{rk}(G)\cdot C_{X,U}^Z(g_*\mathbbold{1}_V)$ is the unique class whose image by $\nu$ is
\begin{align}
\begin{split}
&B_{X,U}^Z(g_*G)-\operatorname{rk}(g_*G)\cdot B_{Y,V}^W(\mathbbold{1}_V)
-\operatorname{rk}(G)\cdot(B_{X,U}^Z(g_*\mathbbold{1}_V)-\operatorname{rk}(g_*\mathbbold{1}_V)\cdot B_{Y,V}^W(\mathbbold{1}_V))\\
=&B_{X,U}^Z(g_*G)-\operatorname{rk}(G)\cdot B_{X,U}^Z(g_*\mathbbold{1}_V).
\end{split}
\end{align}
\end{itemize}
It follows that the commutative diagram~\eqref{eq:diagpfB} implies formula~\eqref{eq:locpf2}.
By~\ref{num:lmCXUZ}, we finish the proof of Proposition~\ref{prop:CXZpf}.

\endproof


\begin{thebibliography}{}

\bibitem[AS07]{AS}
A. Abbes, T. Saito,
\emph{The characteristic class and ramification of an $l$-adic \'etale sheaf}, Invent. Math. \textbf{168} (2007), no. 3, 567-612. 

\bibitem[Abe23]{Abe}
T. Abe,
\emph{Trace formalism for motivic cohomology}, 
\'Epijournal G\'eom. Alg\'ebrique \textbf{7} (2023), Art. 7, 18 pp.

\bibitem[AKMW02]{AKMW}
D. Abramovich, K. Karu, K. Matsuki, J. W\l odarczyk,
\emph{Torification and factorization of birational maps}, 
J. Amer. Math. Soc. \textbf{15} (2002), no. 3, 531-572. 

\bibitem[Alu06]{Alu}
P. Aluffi,
\emph{Limits of Chow groups, and a new construction of Chern-Schwartz-MacPherson classes}, 
Pure Appl. Math. Q. \textbf{2} (2006), no. 4, Special Issue: In honor of Robert D. MacPherson. Part 2, 915-941. 

\bibitem[Ayo07]{Ayo}
J. Ayoub,
\emph{Les six op\'erations de Grothendieck et le formalisme des cycles \'evanescents dans le monde motivique}, Ast\'erisque No. \textbf{314-315} (2007).

\bibitem[Azo22]{Azo}
R. Azouri,
\emph{Motivic characteristic classes for singular spaces}, 
\href{https://arxiv.org/abs/2208.14440}{arXiv:2208.14440}.

\bibitem[Blo85]{Blo}
S. Bloch, 
\emph{Cycles on arithmetic schemes and Euler characteristics of curves}, 
Algebraic geometry, Bowdoin, 1985 (Brunswick, Maine, 1985), 421-450.
Proc. Sympos. Pure Math., \textbf{46}, Part 2.
American Mathematical Society, Providence, RI, 1987.

\bibitem[Blo86]{Blo1}
S. Bloch,
\emph{Algebraic cycles and higher $K$-theory}, 
Adv. Math. \textbf{61} (1986), no. 3, 267-304. 

\bibitem[BD17]{BD}
M. Bondarko, F. D\'eglise,
\emph{Dimensional homotopy t-structures in motivic homotopy theory}, Adv. Math. \textbf{311} (2017), 91-189.

\bibitem[CD16]{CD}
D.-C. Cisinski, F. D\'eglise, 
\emph{\'Etale motives},
Comp. Math. \textbf{152} (2016) 556-666.

\bibitem[CD19]{CD1}
D.-C. Cisinski, F. D\'eglise,
\emph{Triangulated categories of motives}, 
Springer Monographs in Mathematics. Springer, Cham, 2019.

\bibitem[Con07]{Con}
B. Conrad,
\emph{Deligne’s notes on Nagata compactifications}, 
J. Ramanujan Math. Soc. \textbf{22} (2007), no. 3, 205-257.

\bibitem[CP19]{CP}
V. Cossart, O. Piltant,
\emph{Resolution of singularities of arithmetical threefolds}, 
J. Algebra \textbf{529} (2019), 268-535.

\bibitem[DF20]{DF}
F. D\'eglise, J. Fasel,
\emph{The Milnor-Witt motivic ring spectrum and its associated theories}, in \emph{Milnor-Witt Motives},
to appear in Mem. Am. Math. Soc..

\bibitem[DFJK21]{DFJK}
F. D\'eglise, J. Fasel, F. Jin, A. Khan, 
\emph{On the rational motivic homotopy category}, 
J. Ec. Polytech. Math. \textbf{8} (2021), 533-583.

\bibitem[DJK21]{DJK}
F. D\'eglise, F. Jin, A. Khan,
\emph{Fundamental classes in motivic homotopy theory}, 
J. Eur. Math. Soc. \textbf{23} (2021), no. 12, 3935-3993.

\bibitem[DM15]{DM}
F. D\'eglise, N. Mazzari,
\emph{The rigid syntomic ring spectrum}, 
J. Inst. Math. Jussieu \textbf{14} (2015), no.4, 753-799.

\bibitem[Del70]{Del}
P. Deligne,
\emph{Equations diff\'erentielles \`a points singuliers r\'eguliers}, 
Lecture Notes in Mathematics, Vol. 163. Springer-Verlag, Berlin-New York, 1970.

\bibitem[EK20]{EK}
E. Elmanto, A. Khan,
\emph{Perfection in motivic homotopy theory},  Proc. Lond. Math. Soc. \textbf{120} (2020), no. 1, 28-38.

\bibitem[FHM03]{FHM}
H. Fausk, P. Hu, J. P. May,
\emph{Isomorphisms between left and right adjoints}, Theory Appl. Categ. \textbf{11} (2003), No. 4, 107-131.

\bibitem[Fel22]{Fel}
N. Feld,
\emph{Birational invariance of the Chow-Witt group of zero-cycles}, 
\href{https://arxiv.org/abs/2210.03995}{arXiv:2210.03995}.

\bibitem[Ful93]{Fult}
W. Fulton,
\emph{Introduction to toric varieties}, Ann. of Math. Stud., \textbf{131}.
William Roever Lectures Geom.
Princeton University Press, Princeton, NJ, 1993.

\bibitem[Ful98]{Ful}
W. Fulton,
\emph{Intersection theory}, Second edition, Ergebnisse der Mathematik und ihrer Grenzgebiete. 3. Folge. A Series of Modern Surveys in Mathematics, \textbf{2}. Springer-Verlag, Berlin, 1998.

\bibitem[G-S81]{Gon}
G. Gonz\'alez-Sprinberg,
\emph{L'obstruction locale d'Euler et le th\'eor\`eme de MacPherson}, in \emph{Caract\'eristique d'Euler-Poincar\'e: S\'eminaire E. N. S. 1978-1979} (Ast\'erisque 82-83), 7-32. Soc. Math. France, Paris, 1981.

\bibitem[GPS14]{GPS}
M. Groth, K. Ponto, M. Shulman,
\emph{The additivity of traces in monoidal derivators}, J. K-Theory \textbf{14} (2014), no. 3, 422-494. 

\bibitem[Gro58]{Gro}
A. Grothendieck, 
\emph{La th\'eorie des classes de Chern}, 
Bull. Soc. Math. France \textbf{86} (1958), 137-154. 

\bibitem[Gro86]{ReS}
A. Grothendieck,
\emph{R\'ecoltes et Semailles: r\'eflexions et t\'emoignages sur un pass\'e de math\'ematicien}, Universit\'e des sciences et techniques du Languedoc et Centre national de la recherche scientifique, Montpellier, 1986.

\bibitem[Ill15]{Ill}
L. Illusie, 
\emph{From Pierre Deligne’s secret garden: looking back at some of his letters}, 
Jpn. J. Math. (3) \textbf{10}, No. 2, 237-248 (2015). 

\bibitem[Jin16]{Jin}
F. Jin, 
\emph{Borel-Moore motivic homology and weight structure on mixed motives}, 
Math. Z. \textbf{283} (2016), no. 3-4, 1149-1183.

\bibitem[JY21]{JY}
F. Jin, E. Yang,
\emph{K\"unneth formulas for motives and additivity of traces}, 
Adv. Math. \textbf{376} (2021), Article ID 107446.

\bibitem[JY22]{JY2}
F. Jin, E. Yang,
\emph{The quadratic Artin conductor of a motivic spectrum}, 
\href{https://arxiv.org/abs/2211.10985}{arXiv:2211.10985}.

\bibitem[KS90]{KS}
M. Kashiwara, P. Schapira,
\emph{Sheaves on manifolds}, 
With a chapter in French by Christian Houzel. Grundlehren der Mathematischen Wissenschaften \textbf{292}. Springer-Verlag, Berlin, 1990.

\bibitem[KS08]{KSa}
K. Kato, T. Saito,
\emph{Ramification theory for varieties over a perfect field}, Ann. of Math. (2) \textbf{168} (2008), no. 1, 33-96. 

\bibitem[Kel14]{Kel}
S. Kelly, 
\emph{Vanishing of negative $K$-theory in positive characteristic}, 
Compos. Math. \textbf{150} (2014), no. 8, 1425-1434. 

\bibitem[Lau81]{Lau}
G. Laumon,
\emph{Comparaison de caractéristiques d'Euler-Poincaré en cohomologie l-adique}, C. R. Acad. Sci. Paris Sér. I Math.\textbf{292}(1981), no.3, 209-212.

\bibitem[LMS86]{LMS}
L. G. Lewis, J. P. May, M. Steinberger,
\emph{Equivariant stable homotopy theory}, 
With contributions by J. E. McClure. Lecture Notes in Mathematics, \textbf{1213}. Springer-Verlag, Berlin, 1986.

\bibitem[LZ22]{LZ}
Q. Lu, W. Zheng,
\emph{Categorical traces and a relative Lefschetz–Verdier formula},
Forum of Mathematics, Sigma, \textbf{10}, E10, 2022.

\bibitem[Mac74]{Mac}
R. MacPherson,
\emph{Chern classes for singular algebraic varieties}, 
Ann. of Math. (2) \textbf{100} (1974), 423-432.

\bibitem[May01]{May}
J. P. May,
\emph{The additivity of traces in triangulated categories}, 
Adv. Math. \textbf{163} (2001), no. 1, 34-73. 

\bibitem[MVW06]{MVW}
C. Mazza, V. Voevodsky, C. Weibel,
\emph{Lecture notes on motivic cohomology}, 
Clay Math. Monogr., \textbf{2}. American Mathematical Society, Providence, RI; Clay Mathematics Institute, Cambridge, MA, 2006.

\bibitem[Ols16]{Ols}
M. Olsson,
\emph{Motivic cohomology, localized Chern classes, and local terms}, Manuscripta Math. \textbf{149} (2016), no. 1-2, 1-43.

\bibitem[Rio05]{Rio}
J. Riou,
\emph{Dualit\'e de Spanier-Whitehead en g\'eom\'etrie alg\'ebrique}, 
C. R. Math. Acad. Sci. Paris \textbf{340} (2005), no.6, 431-436.

\bibitem[Sai80]{Sai1}
K. Saito,
\emph{Theory of logarithmic differential forms and logarithmic vector fields}, 
J. Fac. Sci. Univ. Tokyo Sect. IA Math. \textbf{27} (1980), no.2, 265-291. 

\bibitem[Sai93]{Sai2}
T. Saito,
\emph{$\epsilon$-factor of a tamely ramified sheaf on a variety}, Invent. Math.\textbf{113} (1993), no.2, 389-417.

\bibitem[Sai17]{Sai}
T. Saito,
\emph{The characteristic cycle and the singular support of a constructible sheaf}, Invent. Math. \textbf{207} (2017), no. 2, 597-695.

\bibitem[SGA4]{SGA4}
M. Artin, A. Grothendieck, J.-L. Verdier,
\emph{Th\'eorie des topos et cohomologie \'etale des sch\'emas}, 
in: S\'eminaire de G\'eom\'etrie Alg\'ebrique du Bois-Marie 1963–1964 (SGA 4).
Dirig\'e par M. Artin, A.Grothendieck, et J. L. Verdier. Avec la collaboration de N. Bourbaki, P. Deligne et B. Saint-Donat, in: Lecture Notes in Mathematics, vol. 269, 270, 305, Springer-Verlag, Berlin-New York, 1972-1973.

\bibitem[SGA5]{SGA5}
A. Grothendieck,
\emph{Cohomologie l-adique et fonctions L}, S{\'e}minaire de g{\'e}om{\'e}trie alg{\'e}brique du Bois-Marie 1965-66 (SGA 5). Avec la collaboration de I. Bucur, C. Houzel, L. Illusie, J.-P. Jouanolou, et J.-P. Serre. Springer Lecture Notes, Vol. 589. Springer-Verlag, Berlin-New York, 1977.

\bibitem[Sil96]{Sil}
R. Silvotti,
\emph{On a conjecture of Varchenko}, 
Invent. Math. \textbf{126} (1996), no. 2, 235-248. 

\bibitem[Spi18]{Spi}
M. Spitzweck,
\emph{A commutative $\mathbb{P}^1$-spectrum representing motivic cohomology over Dedekind domains}, 
M\'em. Soc. Math. Fr. (N.S.) No. 157 (2018), 110 pp.

\bibitem[SV00]{SV}
A. Suslin, V. Voevodsky, 
\emph{Bloch-Kato conjecture and motivic cohomology with finite coefficients}, 
in \emph{The arithmetic and geometry of algebraic cycles}
(Banff, AB, 1998), pp. 117-189, NATO Sci. Ser. CMath. Phys. Sci., 548, Kluwer Acad. Publ., Dordrecht (2000).

\bibitem[Stacks]{Stack}
\emph{The Stacks project}, 
available at https://stacks.math.columbia.edu/, 2020. 

\bibitem[Tem17]{Tem}
M. Temkin,
\emph{Tame distillation and desingularization by $p$-alterations}, 
Ann. Math. (2) \textbf{186}, No. 1, 97-126 (2017).

\bibitem[Tsu11]{Tsu}
T. Tsushima,
\emph{On localizations of the characteristic classes of $\ell$-adic sheaves and conductor formula in characteristic $p>0$}, 
Math. Z. \textbf{269} (2011), no.1-2, 411-447.


\bibitem[Wil06]{Wil}
J. Wildeshaus,
\emph{The boundary motive: definition and basic properties}, Compos. Math. \textbf{142} (2006), no.3, 631-656.

\bibitem[YZ25]{YZ}
E. Yang, Y. Zhao,
\emph{Cohomological Milnor formula and Saito's conjecture on characteristic classes}, 
Invent. Math. \textbf{240} (2025): 123-191.

\end{thebibliography}
\end{document}